\newcommand{\SPRINGER}[2]{
  \ifdefined\FORSPRINGER
    #1
  \else
    #2
  \fi
}
\SPRINGER{
  \RequirePackage{fix-cm}
  \documentclass[smallextended,numbook]{svjour3}
  \smartqed 
  \usepackage{mathptmx} 
}{
  \documentclass{article}
  \usepackage[a4paper]{geometry}
  \pdfoutput=1
}
\usepackage{graphicx}

\usepackage{pdfsync}

\usepackage{latexsym,amsfonts,amsmath,amssymb,bbm}
\SPRINGER{}{
\numberwithin{equation}{section}
\usepackage{amsthm}
\theoremstyle{plain}
  \newtheorem{theorem}{Theorem}[section]
  \newtheorem{lemma}{Lemma}[section]

\theoremstyle{definition}

\theoremstyle{remark}
  
\newcommand{\email}[1]{#1}
\newenvironment{acknowledgements}{\section*{Acknowledgements}}{}
\newcommand{\keywords}[1]{\par\addvspace\medskipamount{\def\and{\ifhmode\unskip\nobreak\fi\ $\cdot$
}\noindent\textbf{Keywords}\enspace\ignorespaces#1\par}}
\newcommand{\subclass}[1]{\par\addvspace\medskipamount{\def\and{\ifhmode\unskip\nobreak\fi\ $\cdot$
}\noindent\textbf{Mathematics Subject Classification
(2000)}\enspace\ignorespaces#1\par}}
}

\usepackage{multirow,url}
\usepackage{mathrsfs}   

\usepackage[usenames,dvipsnames,svgnames]{xcolor}
\usepackage[colorlinks=true,linkcolor=purple,citecolor=purple]{hyperref} 
\usepackage{enumitem} 

\renewcommand{\boldsymbol}[1]{\pmb{#1}} 
\newcommand{\sbinom}[2]{{\textstyle{\binom{#1}{#2}}}}
\newcommand{\satop}[2]{\stackrel{\scriptstyle{#1}}{\scriptstyle{#2}}}

\newcommand{\bsalpha}{{\boldsymbol{\alpha}}}
\newcommand{\bsbeta}{{\boldsymbol{\beta}}}
\newcommand{\bstau}{{\boldsymbol{\tau}}}
\newcommand{\bsgamma}{{\boldsymbol{\gamma}}}
\newcommand{\bsDelta}{{\boldsymbol{\Delta}}}

\newcommand{\bsnu}{{\boldsymbol{\nu}}}
\newcommand{\bszero}{{\boldsymbol{0}}}
\newcommand{\bshalf}{{\boldsymbol{\tfrac{1}{2}}}}

\newcommand{\bsb}{{\boldsymbol{b}}}

\newcommand{\bse}{{\boldsymbol{e}}}

\newcommand{\bsk}{{\boldsymbol{k}}}

\newcommand{\bsm}{{\boldsymbol{m}}}

\newcommand{\bst}{{\boldsymbol{t}}}

\newcommand{\bsw}{{\boldsymbol{w}}}
\newcommand{\bsx}{{\boldsymbol{x}}}
\newcommand{\bsy}{{\boldsymbol{y}}}
\newcommand{\bsz}{{\boldsymbol{z}}}

\newcommand{\rd}{{\mathrm{d}}}
\newcommand{\bbA}{{\mathbb{A}}}
\newcommand{\bbB}{{\mathbb{B}}}
\newcommand{\bbR}{{\mathbb{R}}}
\newcommand{\bbZ}{{\mathbb{Z}}}
\newcommand{\bbN}{{\mathbb{N}}}
\newcommand{\bbE}{{\mathbb{E}}}
\newcommand{\bbP}{{\mathbb{P}}}
\newcommand{\calH}{{\mathcal{H}}}
\newcommand{\calI}{{\mathcal{I}}}

\newcommand{\calO}{{\mathcal{O}}}
\newcommand{\calP}{{\mathcal{P}}}

\newcommand{\calX}{{\mathcal{X}}}
\newcommand{\calY}{{\mathcal{Y}}}
\newcommand{\setu}{{\mathrm{\mathfrak{u}}}}
\newcommand{\setv}{{\mathrm{\mathfrak{v}}}}

\newcommand{\scA}{{\mathscr{A}}}
\newcommand{\wor}{\mathrm{wor}}

\newcommand{\indx}{{\mathfrak F}}
\newcommand{\supp}{{\mathrm{supp}}}
\newcommand{\mask}[1]{{}}
\definecolor{darkred}{RGB}{139,0,0}
\definecolor{darkgreen}{RGB}{0,100,0}
\definecolor{darkmagenta}{RGB}{180,0,180}
\definecolor{darkblue}{RGB}{0,0,190}

\newcommand{\mybullet}{{\raise .5ex\hbox{\tiny$\bullet$}}\;}

\newcommand{\SL}{^{\rm SL}}
\newcommand{\ML}{^{\rm ML}}
\newcommand{\ran}{_{\rm ran}}
\newcommand{\determ}{_{\rm det}}

\newcommand{\X}{\varsigma}
\usepackage[g]{esvect}
\newcommand{\vvec}[1]{\vv{#1}} 

\newcommand{\Bj}{\mathscr{B}\!_j}

\newcommand{\vpsi}{\varpi}

\SPRINGER{\journalname{Foundations of Computational Mathematics}}{}
\begin{document}
\title{Application of quasi-Monte Carlo methods to elliptic PDEs with random diffusion coefficients
-- a survey of analysis and implementation}
\author{
Frances~Y.~Kuo\SPRINGER{}{\footnotemark[1]}
\and
Dirk Nuyens\SPRINGER{}{\footnotemark[2]}
}
\footnotetext[1]{
School of Mathematics and Statistics,
              University of New South Wales, Sydney NSW 2052, Australia\SPRINGER{\\}{,}
              \email{f.kuo@unsw.edu.au}
}
\footnotetext[2]{
Department of Computer Science, KU Leuven,
              Celestijnenlaan 200A, 3001 Leuven, Belgium\SPRINGER{\\}{,}
              \email{dirk.nuyens@cs.kuleuven.be}
}
\SPRINGER{\institute{Frances~Y.~Kuo \at
              School of Mathematics and Statistics,
              University of New South Wales, Sydney NSW 2052, Australia \\
              \email{f.kuo@unsw.edu.au}
           \and
           Dirk Nuyens \at
              Department of Computer Science, KU Leuven,
              Celestijnenlaan 200A, 3001 Leuven, Belgium \\
              \email{dirk.nuyens@cs.kuleuven.be}
}}{}
\date{June 2016}
\maketitle

\begin{abstract}
This article provides a survey of recent research efforts on the
application of quasi-Monte Carlo (QMC) methods to elliptic partial
differential equations (PDEs) with random diffusion coefficients. It
considers, and contrasts, the uniform case versus the lognormal case,
single-level algorithms versus multi-level algorithms, first order QMC
rules versus higher order QMC rules, and deterministic QMC methods versus
randomized QMC methods. It gives a summary of the error analysis and proof
techniques in a unified view, and provides a practical guide to the
software for constructing and generating QMC points tailored to the PDE
problems. The analysis for the uniform case can be generalized to cover a
range of affine parametric operator equations.

\keywords{Quasi-Monte Carlo methods \and
          Infinite dimensional integration \and
          Partial differential equations with random coefficients \and
          Uniform \and Lognormal \and
          Single-level \and Multi-level \and
          First order \and Higher order \and
          Deterministic \and Randomized
          }
\subclass{65D30 \and 65D32 \and 65N30}
\end{abstract}
\markboth{F.~Y. KUO AND D. NUYENS}{Application of QMC methods to PDEs with
random coefficients}

\setcounter{tocdepth}{2} \tableofcontents

\section{Introduction} \label{sec:intro}

In this article we provide a survey of recent research efforts on the
application of quasi-Monte Carlo (QMC) methods to elliptic partial
differential equations (PDEs) with random diffusion coefficients. Such PDE
problems occur in the area of uncertainty quantification. In recent years
many papers have been written on this topic, using a variety of methods
(see below). QMC methods are relatively new to this application area. This
article considers and contrasts different models for the randomness
(uniform versus lognormal) and different algorithms (single-level versus
multi-level, first order versus higher order, deterministic versus
randomized). It gives a summary of the QMC error analysis and proof
techniques in a unified view, and provides a practical guide to the
software for constructing and generating QMC points tailored to the PDE
problems. The analysis for the uniform case can be generalized to cover a
range of affine parametric operator equations.

\subsection{Motivating example} \label{sec:mot}

Many physical, biological or geological models involve spatially varying
input data which may be subject to uncertainty. This induces a
corresponding uncertainty in the outputs of the model and in any physical
quantities of interest which may be derived from these outputs. A common
way to deal with these uncertainties is by considering the input data to
be a \emph{random field}, in which case the derived quantity of interest
will in general also be a random variable or a random field. The
computational goal is usually to find the expected value, higher order
moments or other statistics of these derived quantities.

A prime example is the flow of water through a disordered porous medium.
Because of the near impossibility of modelling the microscopic channels
through which water can flow in a porous layer, it is common engineering
practice to model the porous medium as a random field. Mathematically, the
flow can be modeled by Darcy's law coupled with the mass conservation law,
i.e.,
\begin{align*}
 q(\bsx,\omega) + a(\bsx,\omega)\,\nabla p(\bsx,\omega) &\,=\, 0\,, \\
 \nabla\cdot q(\bsx,\omega) &\,=\, 0\,,
\end{align*}
for $\bsx$ in a bounded domain $D\subset\bbR^d$, $d\in\{1,2,3\}$, and for
almost all events $\omega$ in the probability space
$(\Omega,\mathcal{A},\bbP)$. Here $q(\bsx,\omega)$ is the velocity (also
called the specific discharge) and $p(\bsx,\omega)$ is the residual
pressure, while $a(\bsx,\omega)$ is the permeability (or more precisely,
the ratio of permeability to dynamic viscosity) which is modelled as a
random field. Uncertainty in $a(\bsx,\omega)$ leads to uncertainty in
$q(\bsx,\omega)$ and $p(\bsx,\omega)$. Quantities of interest include for
example the breakthrough time of a plume of pollution moving through the
medium.

To compute the expected value of some quantity of interest, one can
generate a number of realizations of the random permeability field, for
each realization solve the PDE numerically and compute the quantity of
interest, and then take the average of all solutions from different
realizations. This describes \emph{Monte Carlo
\textnormal{(}MC\textnormal{)} simulation}, and it is regularly employed
in such problems, see, e.g., \cite{Dag84,NHS98a,NHS98b,Zha02,Rub03,LZ04}.
By definition, expected values are integrals, with the dimensionality as
high as the number of parameters needed to describe the randomness. This
leads to the consideration of \emph{quasi-Monte Carlo
\textnormal{(}QMC\textnormal{)} methods}, which are quadrature methods for
tackling high dimensional integrals, with the hope to improve upon the
slow convergence of MC simulation.

Throughout this article we refer to the number of integration variables
$s$ as the ``stochastic dimension'', which can be in the hundreds or
thousands or more, in contrast to the ``spatial dimension'' $d$ which is
just $1$, $2$ or $3$.

\subsection{The QMC story}

QMC methods \cite{Nie92,SJ94}, including the families of ``lattice rules"
and ``digital nets", are equal-weight quadrature rules where the
quadrature points are chosen deterministically and designed cleverly to
beat the random sampling of MC. They date back to the 1960s from number
theorists, but the theories of that era were \emph{not} adequate for very
high dimensional problems because either that fast convergence was
obtained by assuming periodicity which is unrealistic in high dimensions,
or that the typical error bounds for a QMC method requiring $n$ function
evaluations in $s$ variables were of the form $c_s(\log n)^{s-1}/n$. In
the latter case although the convergence rate appears faster than the MC
rate of $1/\sqrt{n}$, the fatal flaw is that for fixed $s$, the function
$(\log n)^{s-1}/n$ {\em increases} with increasing $n$ until
$n\approx\exp(s)$, a number that is truly astronomical if $s$ is large.

So until perhaps the middle of the 1990s it was generally thought that QMC
methods would not be effective in dimensions of more than say $20$ or
$30$. But then a dramatically successful computational experiment of
treating a $360$-dimensional integral coming from Wall Street \cite{PT95}
changed people's perceptions of what might be possible. At that time there
was nothing in the available QMC theory that could explain the success.
This led to many theoretical developments, as researchers struggled to
understand how such high dimensionality could be handled successfully by
QMC methods.

Modern QMC analysis takes into account that integrands from practical
applications can have ``low effective dimension'' \cite{CMO97}, meaning
that the problem although having a very high nominal dimension may in fact
depend mainly on a number of leading variables, or may be mainly affected
by the interaction of a small number of variables at a time. This concept
was formalized in theory through the introduction of ``weighted'' function
spaces \cite{SW98}: a set of parameters called ``weights'' are built into
the function space norm to model the relative importance between different
subsets of variables. Then ``tractability'' \cite{NW08,NW10,NW12} analysis
has been conducted in these settings to obtain, for example, a necessary
and sufficient condition on the weights under which the integration error
in the so-called ``worst case'' sense is \emph{independent of the
dimension~$s$}. Thus, rather than saying that all high dimensional
problems can be successfully tackled by QMC methods, we now know how to
recognize and analyze mathematically the particular features that make
some high dimensional problems manageable.

Following the introduction of weighted spaces, a new class of constructive
methods known as the ``component-by-component (CBC) construction'' has
flourished \cite{SR02,SKJ02b}. The algorithm constructs the components of
the ``generating vector'' for lattice rules one at a time: the $(j+1)$-th
component is obtained by successive $1$-dimensional searches, with the
previous $j$ components kept unchanged. The inductive nature of such
algorithms provides the means to obtain new methods for arbitrarily high
dimensions. It has been established that this ``greedy'' algorithm yields
lattice rules which achieve the optimal rate of convergence close to order
$n^{-1}$ in the underlying weighted function spaces, with the implied
constant independent of $s$, see \cite{Kuo03,Dic04}. This was followed by
many further works, most notably the use of fast Fourier transform (FFT)
to speed up the computation \cite{NC06a,NC06b}, the construction of
``extensible lattice sequences'' \cite{CKN06,DPW08,HKKN12}, the use of
``tent transform'' to achieve close to order $n^{-2}$ convergence
\cite{Hic02,DNP14}, the carrying over of the lattice technology to digital
nets and sequences \cite{DP05,DKPS05,NC06c}, and the revolutionary
invention of ``higher order digital nets'' which allow a convergence rate
of order $n^{-\alpha}$, $\alpha>1$, for sufficiently smooth integrands
\cite{Dic07,Dic08,Dic09,BD09,BDGP11,BDLNP12,GSY16}. For surveys of these
recent QMC developments see \cite{DP10,KSS11,DKS13,Nuy14,Slo15}.

By now weights of many forms have been considered in the literature. At
the unrealistic extreme we have ``general weights'' which allow for a
different weight parameter to be attached to each of the $2^s$ subsets of
the indices from $1$ to $s$. The original and simplest weights from
\cite{SW98} are now called ``product weights'', and in between we have
``finite order'' weights, ``finite diameter'' weights, and ``order
dependent'' weights \cite{DSWW04,SWW04}. Furthermore, there are the more
recent additions in the context of applying QMC to PDE problems: ``POD''
(for ``product and order dependent'') weights \cite{KSS12} and ``SPOD''
(for ``smoothness-driven product and order dependent'') weights
\cite{DKLNS14}. The driving motivation for this flowering of possibilities
has been the desire to describe in a more precise way the influence of
particular combinations of the variables. For the CBC construction of
lattice rules mentioned above, the weights definitely matter, since they
appear as parameters in the algorithm that determine the integration rule.
The philosophy is therefore to choose the weights according to the
dimension structure of the practical integrands and then construct QMC
methods that are tailored to the given application.

\subsection{Progress on PDEs with random coefficients}

Returning to the motivating example in Subsection~\ref{sec:mot}, because
the permeability $a(\bsx,\omega)$ is physically positive, it is popular
and natural to assume that $a(\bsx,\omega)$ is a ``lognormal'' random
field, that is, $\log (a(\bsx,\omega))$ is a Gaussian random field on the
spatial domain $D$ with a specified mean and a covariance function. There
is some evidence from field data that lognormality gives a reasonable
representation of reality in certain cases \cite{Fre75,HK85}. See below
for many recent references which considered lognormal random fields.

A common approach to represent the random field $a(\bsx,\omega)$ is to use
the \emph{Karhunen-Lo\`eve} (\emph{KL}) \emph{expansion} \cite{Loe78} to
write $\log (a(\bsx,\omega))$ as an infinite series involving the
eigenvalues and eigenfunctions of the integral operator associated with
the covariance function, where the series is parametrised and depends
linearly on a sequence $y_j = y_j(\omega)$, $j\ge 1$, of i.i.d.\ standard
normal random numbers from $\bbR$. In practical computations the infinite
sum is truncated to, say, $s$ terms, giving rise to a truncation error to
be managed. While this approach can be very effective when the KL
expansion converges rapidly, it faces the serious challenges of high cost
combined with large truncation error when the convergence of the KL
expansion is slow.

Instead of sampling the continuous random field \emph{everywhere} in the
spatial domain by a truncated KL expansion, an alternative approach is to
sample the random field only at a \emph{discrete set of grid points} with
respect to the covariance matrix inherited from the given covariance
function of the continuous field. The random field is then represented
exactly at these grid points, thus eliminating completely the truncation
error. (However, interpolation would be required at the time of assembling
the stiffness matrix for solving the PDE numerically, and care is needed
to ensure that the interpolation error is no worse than, e.g., the
quadrature error in integrating the finite element basis functions.) The
resulting large matrix factorization problem can be handled by the
``circulant embedding'' technique, see, e.g., \cite{DN97,GKNSS11,KSSSU}.

The unbounded parameters $\bsy = (y_j)_{j\ge 1}$ from the lognormal model
present some challenges in the theoretical analysis. One major challenge
is that the random coefficient is not uniformly bounded from above and
below and so the Lax-Milgram lemma cannot be applied directly. Many
researchers therefore consider the simpler ``uniform'' model where
$a(\bsx,\omega)$ is written as an infinite series that depends linearly on
a sequence $y_j = y_j(\omega)$, $j\ge 1$, of i.i.d.\ uniform random
numbers from a bounded interval of $[-1,1]$ or
$[-\frac{1}{2},\frac{1}{2}]$.

There is a huge body of literature on treating these PDEs with random
coefficients. Some methods apply simultaneous approximation in both
physical and probability space; these go under names such as ``stochastic
Galerkin'', ``stochastic collation'', ``polynomial chaos'', or
``generalized polynomial chaos'', see, e.g.,
\cite{GS91,XK03,BTZ04,ST06,ST07,BNT07,NTW08a,NTW08b,CDS10,CD15,SG11,HS12,HS13,KunSch13,GWZ14,CCS15}.
In recent times these methods are also combined with ``multi-level''
schemes to reduce the computational cost without loss of accuracy, see,
e.g., \cite{BSZ11,CGST11,CST13,TSGU13,CHNST15,HPS15a,HPS15b,TJWG}. The
methods are also being applied to the area of Bayesian inversion, see,
e.g., \cite{HSS12,SS13,SS14}.

Table~\ref{tab:papers} provides a summary of some recent progress on the
application of QMC to PDEs with random coefficients. Firstly,
comprehensive numerical experiments were carried out in \cite{GKNSS11}
showing promising results for the lognormal case with the circulant
embedding strategy, but without any theoretical justification. The first
theoretical analysis was done in \cite{KSS12} for the simpler uniform case
under the KL expansion framework. The analysis was then generalized to the
lognormal case in \cite{GKNSSS15} and extended to a multi-level scheme for
the uniform case in \cite{KSS15}. The use of first order QMC methods in
the uniform case was then replaced by higher order QMC methods in
\cite{DKLNS14}, and the corresponding multi-level analysis was done in
\cite{DKLS}. The multi-level analysis for the lognormal case was done in
\cite{KSSSU}, while the analysis for the discrete sampling method combined
with circulant embedding technique is being considered in \cite{GKNSS}.
There are also other QMC related works, see, e.g.,
\cite{Leg12,Sch12,GS14,HPS15c,RNV15,GH15,DLS}.

\begin{table}[t]
 \caption{Application of QMC to PDEs with random coefficients.
 The * indicates that there are accompanying numerical results.}\label{tab:papers}
 \begin{center}
 \begin{tabular}{|c|c|c|c|}
 \hline
 & Uniform & \multicolumn{2}{c|}{Lognormal} \\
 & KL expansion & KL expansion & Circulant embedding \\
 \hline
  Numerical experiments only & & &  \hphantom{*}\cite{GKNSS11}* \\
  First order single-level analysis & \cite{KSS12} & \hphantom{*}\cite{GKNSSS15}* & \hphantom{*}\cite{GKNSS}* \\
  First order multi-level analysis & \cite{KSS15} & \hphantom{*}\cite{KSSSU}* & \\
  Higher order single-level analysis & \cite{DKLNS14} & & \\
  Higher order multi-level analysis & \hphantom{*}\cite{DKLS}* & & \\
 \hline
 \end{tabular}
 \end{center}
\end{table}

\subsection{Overview of this article}

This article surveys the results from
\cite{KSS12,KSS15,GKNSSS15,KSSSU,DKLNS14,DKLS} in a unified view. The
first order results \cite{KSS12,KSS15,GKNSSS15,KSSSU} are based on
\emph{randomly shifted lattice rules} and are accompanied by probabilistic
error bounds. The higher order results \cite{DKLNS14,DKLS} are based on
\emph{interlaced polynomial lattice rules} and are accompanied by
deterministic error bounds. Part of the gain in the improved convergence
rates arises because of the switch from an $\ell^2$ norm in the function
space setting to an $\ell^\infty$ norm. The lognormal results
\cite{GKNSSS15,KSSSU} require a non-standard function space setting for
integrands with domain $\bbR^s$ for some $s$. There is as yet no QMC
theory that can give higher order convergence for the lognormal case.

We will \emph{not} discuss the theory \cite{GKNSS} for the case of
circulant embedding. Also, the results in \cite{DKLNS14,DKLS} were
formulated for general affine parametric operator equations, but we will
only touch on this aspect very briefly in the article.

This article includes a practical guide on how to make use of the QMC
technology for this and other PDE problems. Computer programs are
available from the website \\
\url{http://people.cs.kuleuven.be/~dirk.nuyens/qmc4pde/}.

The structure of this article is as follows. In Section~\ref{sec:para} we
give a parametric formulation of the PDE problem and discuss both the
uniform and lognormal cases. In Section~\ref{sec:appr} we provide some
minimum background on QMC and finite element (FE) methods as well as the
dimension truncation analysis. In Section~\ref{sec:algo} we introduce the
single-level and multi-level algorithms for both the deterministic and
randomized variants. In Section~\ref{sec:qmc} we introduce three weighted
function space settings for QMC error analysis. In Section~\ref{sec:error}
we summarize the error analysis from
\cite{KSS12,KSS15,GKNSSS15,KSSSU,DKLNS14,DKLS} in a high-level unified
way. We explain the strategies and proof techniques, deferring proofs to
Section~\ref{sec:proof}, the Appendix. In Section~\ref{sec:comp} we
provide a practical guide on how to use the software from the website for
constructing and generating QMC points. In Section~\ref{sec:conc} we end
the article with some concluding remarks.

\section{Uniform versus lognormal coefficients} \label{sec:para}

With the motivating example from Subsection~\ref{sec:mot} in mind, we
consider throughout this article a model parametric elliptic problem with
homogeneous Dirichlet condition. This model problem has been considered in
many papers. It is simple enough for illustrating the kind of results that
we can obtain as well as the proof techniques. The strategy can be
extended to more general domains and boundary conditions as well as other
PDE problems.

We take the view that the random coefficient $a(\bsx,\omega)$ has been
parameterized by a vector $\bsy(\omega) =
(y_1(\omega),y_2(\omega),\ldots)$, and, for fixed $\omega$, we denote the corresponding
deterministic parametric coefficient by $a(\bsx,\bsy)$. Specifically, for a given parameter $\bsy$ we
consider the parametric elliptic Dirichlet problem
\begin{equation}\label{eq:strong}
 - \nabla \cdot (a(\bsx,\bsy)\,\nabla u(\bsx,\bsy))
 \,=\, f(\bsx) \quad\mbox{for $\bsx$ in $D$}\,,
 \quad
 u(\bsx,\bsy) \,=\, 0 \quad\mbox{for $\bsx$ on $\partial D$}\,,
\end{equation}
for domain $D\subset\bbR^d$ a bounded, convex, Lipschitz polyhedron with
boundary $\partial D$, where the spatial dimension $d=1,2$, or $3$ is
assumed given and fixed. The differential operators in~\eqref{eq:strong}
are understood to be with respect to the physical variable $\bsx$ which
belongs to $D$. The parametric variable $\bsy = (y_j)_{j\geq 1}$ belongs
to either a bounded or unbounded domain, depending on which of the two
popular formulations of the parametric coefficient $a(\bsx,\bsy)$ is being
considered: the ``uniform'' case or the ``lognormal'' case; see below.

\subsection{Uniform case}

In the ``uniform'' case, we assume that the parameter $\bsy$ is
distributed on
\begin{align*}
  U \,:=\, [-\tfrac{1}{2},\tfrac{1}{2}]^\bbN\,
\end{align*}
with the uniform probability measure $\mu(\rd\bsy) = \bigotimes_{j\geq 1}
\rd y_j = \rd\bsy$. Throughout the article $\bbN = \{1,2,3,\ldots\}$
denotes the set of positive integers. The parametric coefficient
$a(\bsx,\bsy)$ is further assumed to depend linearly on the parameters
$y_j$ as follows:
\begin{align} \label{eq:axy-unif}
  a(\bsx,\bsy)
  \,=\, a_0(\bsx) + \sum_{j\geq 1} y_j\, \psi_j(\bsx)\,,
  \qquad
  \bsx\in D\,, \quad\bsy\in U\,.
\end{align}
The functions $\psi_j$ can arise from either the eigensystem of a
covariance operator or other suitable function systems in $L^2(D)$.

We will impose a number of assumptions on $a_0$ and $\psi_j$ as required.
In the following, the $L^\infty(D)$ norm is defined as the essential
supremum in $D$ as per usual, and for $\|\nabla v\|_{L^\infty}$ we take
the essential supremum of the Euclidean norm of $\nabla v$. The
$W^{1,\infty}(D)$ norm is defined by $\|v\|_{W^{1,\infty}} :=
\max\{\|v\|_{L^\infty},\|\nabla v\|_{L^\infty}\}$.

\begin{enumerate}[label=\textbf{\textup{(U\arabic*)}},leftmargin=2.5em]
\item\label{U1}%
We have $a_0 \in L^\infty(D)$ and $\sum_{j\geq 1}\|\psi_j\|_{L^\infty}
< \infty$.
\item\label{U2}%
There exist $a_{\max}$ and $a_{\min}$ such that $0 < a_{\min}\le
a(\bsx,\bsy)\le a_{\max}<\infty$ for all $\bsx\in D$ and $\bsy\in U$.
\item\label{U3}%
We have $a_0 \in W^{1,\infty}(D)$ and $\sum_{j\geq
1}\|\psi_j\|_{W^{1,\infty}} < \infty$.
\item\label{U4}%
The sequence $\psi_j$ is ordered so that $\|\psi_1\|_{L^\infty}\ge
\|\psi_2\|_{L^\infty} \ge\cdots$.
\item\label{U5}%
There exists $p_0\in (0,1)$ such that $\sum_{j\ge 1}
\|\psi_j\|^{p_0}_{L^\infty}< \infty$.
\item\label{U6}%
There exists $p_1 \in (0,1)$ such that $\sum_{j\ge 1}
\|\psi_j\|^{p_1}_{W^{1,\infty}}< \infty$.
\item\label{U7}%
For a non-negative integer $t$, there exists $p_t \in (0,1)$ such that
$\sum_{j\ge 1} \|\psi_j\|^{p_t}_{\calX_t}< \infty$, where $\calX_t$
denotes a Sobolev space of functions on $D$ with smoothness scale $t$,
with $t$ roughly corresponding to the number of derivatives that exist
in $\bsx$.
\end{enumerate}
Assumption~\ref{U1} ensures that the coefficient $a(\bsx,\bsy)$ is
well-defined for all parameters $\bsy\in U$. Assumption~\ref{U2} yields
the continuity and coercivity needed for the standard FE analysis, so that
a unique solution exists. Assumption~\ref{U3} guarantees that the FE
solutions converge to the solution of~\eqref{eq:strong}.
Assumption~\ref{U4} enables the analysis for truncating the infinite sum
in~\eqref{eq:axy-unif}. Assumption~\ref{U5} implies decay of the
fluctuation coefficients $\psi_j$, with faster decay for smaller $p_0$;
Assumptions~\ref{U6} and~\ref{U7} have similar implications. Although not
explicitly specified, we are interested in those values of $p_0$, $p_1$,
and $p_t$ which are as small as possible. Typically we have $p_0 < p_1 <
p_2 < \cdots$. The values of $p_0$, $p_1$, and $p_t$ will determine our
QMC convergence rates for different algorithms.

For convenience of later analysis, we define the sequence $\bsb =
(b_j)_{j\ge 1}$ by
\begin{align} \label{eq:bj-unif}
  b_j \,:=\, \frac{\|\psi_j \|_{L^\infty}}{a_{\min}} \,,
  \quad j\ge 1\,,
\end{align}
and define the sequence $\overline\bsb = (\overline{b}_j)_{j\ge 1}$ by
\begin{align} \label{eq:barbj-unif}
  \overline{b}_j
  \,:=\,
  \frac{\|\psi_j\|_{W^{1,\infty}}}{a_{\min}}
  \,=\,
  \frac{\max(\|\psi_j\|_{L^\infty},\|\nabla \psi_j\|_{L^\infty})}{a_{\min}}
  \,\ge\, b_j\,,
  \quad j\ge 1\,.
\end{align}

Our goal is to compute the integral, i.e., the expected value, with
respect to $\bsy\in U$, of a bounded linear functional $G$ applied to the
solution $u(\cdot,\bsy)$ of the PDE \eqref{eq:strong}
\begin{align} \label{eq:int-unif}
  I(G(u))
  &\,\hphantom{:}=\,
  \int_{[-\tfrac{1}{2},\tfrac{1}{2}]^\bbN} G(u(\cdot,\bsy))\,\rd\bsy \\
  &\,:=\, \lim_{s\to\infty} \int_{[-\tfrac{1}{2},\tfrac{1}{2}]^s}
        G(u(\cdot,(y_1,\ldots,y_s,0,0,\ldots)))\,\rd y_1\cdots\rd y_s\,. \nonumber
\end{align}
We remark that our analysis relies heavily on the boundedness and
linearity of $G$, especially in the duality arguments.

The uniform framework can be extended to the general framework of
``affine'' parametric operator equations, see \cite{Sch12} as well as
\cite{DKLNS14,DKLS}. Let $\{A_j\}_{j\ge 1}$ denote a sequence of bounded
linear operators in $L(\calX,\calY^*)$, between suitably defined spaces
$\calX$ and $\calY^*$. For every $f\in \calY^*$ and every $\bsy\in U$, the
task is to seek $u(\bsy)\in\calX$ such that $A(\bsy)\,u(\bsy) = f$, where
$A(\bsy) = A_0 + \sum_{j\ge 1} y_j\, A_j$. In this general setting,
coercivity is replaced by inf-sup conditions, and the results depend,
e.g., on the summability of $\sum_{j\ge
1}\|A_j\|_{L(\calX,\calY^*)}^{p_0}$ for $p_0\in (0,1)$. We will not
discuss this general framework further in this article, other than to
summarize some results from \cite{DKLNS14,DKLS} at the end of
Section~\ref{sec:error}.

\subsection{Lognormal case}

In the ``lognormal'' case, we assume that the parameter $\bsy$ is
distributed on $\bbR^\bbN$ according to the product Gaussian measure
$\mu_G = \bigotimes_{j\geq 1} N(0,1)$. The parametric coefficient
$a(\bsx,\bsy)$ now takes the form
\begin{align} \label{eq:axy-logn}
  a(\bsx,\bsy) \,=\,
  a_0(\bsx)\exp\bigg(\sum_{j\ge 1}
  y_j\,\sqrt{\mu_j}\,\xi_j(\bsx)\bigg)\,,
  \qquad
  \bsx\in D\,, \quad\bsy\in \bbR^\bbN\,,
\end{align}
where $a_0(\bsx)>0$.

The coefficient $a(\bsx,\bsy)$ of the form~\eqref{eq:axy-logn} can arise
from the \emph{Karhunen-Lo\`eve} (KL) \emph{expansion} in the case where
$\log(a)$ is a stationary Gaussian random field with a specified mean and
a covariance function. As an example we focus on the isotropic
\emph{Mat\'ern} covariance $\rho_\nu(|\bsx-\bsx'|)$, with
\begin{equation} \label{eq:matern}
  \rho_\nu (r) \,:=\, \sigma^2 \frac{2^{1-\nu}}{\Gamma(\nu)}
  \left(2\sqrt{\nu} \frac{r}{\lambda_C}\right)^\nu
  K_\nu\left(2\sqrt{\nu} \frac{r}{\lambda_C} \right)\,,
\end{equation}
where $\Gamma$ is the gamma function and $K_\nu$ is the modified Bessel
function of the second kind. The parameter $\nu>1/2$ is a smoothness
parameter, $\sigma^2$ is the variance and $\lambda_C$ is the correlation
length scale. Then $\{(\mu_j,\xi_j)\}_{j\ge1}$ is the sequence of
eigenvalues and eigenfunctions of the integral operator $(Rw)(\bsx) =
\int_D \rho_\nu(|\bsx-\bsx'|)\,w(\bsx')\,\rd\bsx'$, with eigenvalues
$\mu_j$ enumerated in nonincreasing order and with eigenfunctions $\xi_j$
normalized in $L^2(D)$. Moreover, the sequence $\{\xi_j\}_{j\ge 1}$ form
an orthonormal basis in $L^2(D)$.

We define the sequence $\bsbeta = (\beta_j)_{j\ge 1}$ by
\begin{align} \label{eq:bj-logn}
  \beta_j \,:=\, \sqrt{\mu_j}\,\|\xi_j\|_{L^\infty}\,,
  \qquad j\ge 1\,,
\end{align}
and define the set of admissible parameters
\begin{align*}
  U_\bsbeta \,:=\,
  \bigg\{ \bsy\in \bbR^\bbN \;:\; \sum_{j\ge 1} \beta_j\, |y_j| < \infty \bigg\} \,\subseteq\, \bbR^\bbN\,.
\end{align*}
We also define the sequence $\overline\bsbeta = (\overline{\beta}_j)_{j\ge
1}$ by
\begin{align} \label{eq:barbj-logn}
  \overline{\beta}_j \,:=\, \sqrt{\mu_j}\,\|\xi_j\|_{W^{1,\infty}}
  \,=\, \sqrt{\mu_j}\,\max\{\|\xi_j\|_{L^\infty},\|\nabla \xi_j\|_{L^\infty}\}
  \,\ge\, \beta_j\,,
  \qquad j\ge 1\,,
\end{align}
and define analogously $U_{\overline\bsbeta} \subseteq U_\bsbeta \subseteq
\bbR^\bbN$.

Similarly to the uniform case, we will impose a number of assumptions in
the lognormal case as required. We follow the setting of \cite{KSSSU},
with the exception of Assumption~\ref{L4} which came from \cite{GKNSSS15}.

\begin{enumerate}[label=\textbf{\textup{(L\arabic*)}},leftmargin=2.5em]
\item\label{L1}%
We have $a_0 \in L^\infty(D)$ and $\sum_{j\geq 1} \beta_j < \infty$.
\item\label{L2}%
For every $\bsy\in U_\bsbeta$, the expressions $a_{\max}(\bsy) :=
\max_{\bsx\in \overline{D}} a(\bsx,\bsy)$ and $a_{\min}(\bsy) :=
\min_{\bsx\in \overline{D}} a(\bsx,\bsy)$ are well defined and satisfy
$0 < a_{\min}(\bsy) \le a(\bsx,\bsy) \le a_{\max}(\bsy) < \infty$.
\item\label{L3}%
We have $a_0 \in W^{1,\infty}(D)$ and $\sum_{j\geq 1} \overline\beta_j
< \infty$.
\item\label{L4}%
There exist $C_1,C_2>0$, $\Theta>1$, and $\varepsilon\in
[0,\frac{\Theta-1}{2\Theta})$ such that $\mu_j \le C_1 j^{-\Theta}$
and $\|\xi_j\|_{C^0(\overline{D})} +
\mu_j\,\|\nabla\xi_j\|_{C^0(\overline{D})} \le C_2\,
\mu_j^{-\varepsilon}$ for $j\ge 1$.
\item\label{L5}%
There exists $p_0\in (0,1)$ such that $\sum_{j\ge 1} \beta_j^{p_0}<
\infty$.
\item\label{L6}%
There exists $p_1 \in (0,1)$ such that $\sum_{j\ge 1}
\overline{\beta}_j^{p_1}< \infty$.
\end{enumerate}
Assumption~\ref{L1} ensures that the series~\eqref{eq:axy-logn} converges
in $L^\infty(D)$ for every $\bsy\in U_\bsbeta$ and that $\mu_G(U_\bsbeta)
= 1$, see \cite[Lemma 2.28]{SG11}. Moreover, Assumption~\ref{L1} implies
Assumption~\ref{L2}, which in turn yields the continuity and coercivity of
the bilinear form (see~\eqref{eq:bilinear} below) for every $\bsy\in
U_\bsbeta$. Assumption~\ref{L3} ensures that the
series~\eqref{eq:axy-logn} converges in $W^{1,\infty}(D)$ for every
$\bsy\in U_{\overline\bsbeta}$ and that $\mu_G(U_{\overline\bsbeta}) = 1$;
it also guarantees that for every $\bsy\in U_{\overline\bsbeta}$ the FE
solutions converge to the solution of~\eqref{eq:strong}.
Assumption~\ref{L4} is needed for the dimension truncation result.
Assumptions~\ref{L5} and~\ref{L6} play analogous roles to
Assumptions~\ref{U5} and~\ref{U6} in the uniform case. Typically we have
$p_0 < p_1$.

Our goal is again to compute the integral of a bounded linear functional
$G$ applied to the solution $u(\cdot,\bsy)$ of the PDE, but now the
integral is over $\bsy\in \bbR^\bbN$ with a countable product Gaussian
measure $\mu_G(\rd\bsy)$. Recall that we restrict ourselves to $\bsy\in
U_\bsbeta$ with full measure $\mu_G(U_\bsbeta) = 1$. Abusing the standard
notation in measure and integration theory, we write this integral simply
as
\begin{align} \label{eq:int-logn}
 I(G(u))
 &\,=\, \int_{\bbR^\bbN} G(u(\cdot,\bsy))\,\prod_{j\ge 1} \phi(y_j)\,\rd\bsy
 \,=\, \int_{[0,1]^\bbN} G(u(\cdot,\Phi^{\mbox{-}1}(\bsw)))\,\rd\bsw\,.
\end{align}
Here $\phi(y) := \exp(-y^2/2)/\sqrt{2\pi}$ is the univariate standard
normal probability density function. Denoting the corresponding cumulative
distribution function by $\Phi(y) := \int_{-\infty}^y \phi(t)\,\rd t$ and
its inverse by $\Phi^{\mbox{-}1}$, we can apply the change of variables
componentwise as follows:
\[
  \bsy \,=\, \Phi^{\mbox{-}1}(\bsw) \,=\,
  (\Phi^{\mbox{-}1}(w_1),\Phi^{\mbox{-}1}(w_2),\ldots) \,\in\, \bbR^\bbN
  \qquad\mbox{for}\qquad
  \bsw\,\in\, (0,1)^\bbN\,.
\]
This leads to the transformed integral over the unit cube on the
right-hand side of~\eqref{eq:int-logn}. The equivalence of the integrals
in~\eqref{eq:int-logn} follows from Kakutani's theorem on the equivalence
of infinite product measures (see, e.g., \cite{Bog98}).

\section{Quadrature, spatial discretization, dimension truncation} \label{sec:appr}

\subsection{QMC quadrature}

For $F$ a general real-valued function defined over the $s$-dimensional
unit cube $[0,1]^s$, with $s$ finite and fixed, we consider the integral
\begin{align*}
  I(F)
  \,=\, \int_{[0,1]^s} F(\bsy)\, \rd\bsy\,.
\end{align*}
An $n$-point \emph{quasi-Monte Carlo} (QMC) approximation to this integral
is an equal-weight quadrature rule of the form
\begin{align} \label{eq:qmc}
  Q(F)
  \,=\, \frac{1}{n} \sum_{i=1}^n F(\bst_i)\,,
\end{align}
with carefully chosen points $\bst_1,\ldots,\bst_n\in [0,1]^s$. See
\cite{DKS13} for a comprehensive survey of recent developments on QMC
methods. In this article we will consider two families of QMC methods:
\emph{randomly shifted lattice rules} and (deterministic) \emph{interlaced
polynomial lattice rules}. More details about these QMC methods and their
error analysis will be given in Section~\ref{sec:qmc}. Here we provide
only a brief overview of the general framework.

We define the \emph{worst case error} for QMC integration in some Banach
space of functions $\calH$ to be $e^\wor(\bst_1,\ldots,\bst_n) :=
  \sup_{F\in\calH,\,\|F\|_{\calH}\le 1} |I(F) - Q(F)|$.
Then we have the error bound
\begin{align*}
  |I(F) - Q(F)| \,\le\, e^\wor(\bst_1,\ldots,\bst_n)\,\|F\|_{\calH}
  \qquad\mbox{for all}\quad F\in\calH\,.
\end{align*}
An error bound of this form conveniently separates the dependence on the
QMC point set from the dependence on the integrand. For a given integrand
$F$, the general idea is to choose a suitable function space $\calH$ so
that the norm $\|F\|_\calH$ is finite, and then to construct QMC points
$\bst_1,\ldots,\bst_n$ to make the worst case error
$e^\wor(\bst_1,\ldots,\bst_n)$ as small as possible. There could be a
trade-off between these two quantities, but the ultimate goal is to make
the product of the two quantities as small as possible.

The advantages of deterministic QMC methods include the exact
reproducibility and the fully deterministic theoretical error bound; these
properties might be favored by practitioners. However, one may also argue
that deterministic QMC methods have the drawback of being biased and
lacking a practical error estimate. In contrast, randomized QMC methods
are unbiased and a practical error estimate can be easily obtained. Here
we discuss only the simplest kind of randomization, namely, \emph{random
shifting}. For a fixed shift $\bsDelta\in [0,1]^s$, the $\bsDelta$-shift
of the QMC rule~\eqref{eq:qmc} is
\begin{align} \label{eq:qmc-sh}
  Q(F;\bsDelta)
  \,=\, \frac{1}{n} \sum_{i=1}^n F(\{\bst_i + \bsDelta\})\,,
\end{align}
where the braces around a vector indicate that we take the fractional part
of each component in the vector. Essentially, we move all the QMC points
by the same amount, and if any point falls outside of the unit cube it is
simply ``wrapped'' back in from the opposite side. If the shift $\bsDelta$
is generated randomly from the uniform distribution on $[0,1]^s$, then it
is easy to verify that $\bbE[Q(F;\cdot)] = \int_{[0,1]^s}
Q(F;\bsDelta)\,\rd\bsDelta = I(F)$, that is, a randomly shifted QMC rule
provides an unbiased estimate of the integral, and in turn, the variance
of $Q(F;\cdot)$ is precisely the mean-square error $\bbE[|I(F) -
Q(F;\cdot)|^2]$. A probabilistic error bound for a randomly shifted QMC
rule in $\calH$ is
\begin{align} \label{eq:sh-avg-wce}
  \sqrt{\bbE[|I(F) - Q(F;\cdot)|^2]}
  \,\le\, e^\wor_{\rm sh}(\bst_1,\ldots,\bst_n)\,\|F\|_{\calH}
  \qquad\mbox{for all}\quad F\in\calH\,,
\end{align}
where the quantity
$e^\wor_{\rm sh}(\bst_1,\ldots,\bst_n) :=
(\int_{[0,1]^s}
(e^\wor(\{\bst_1+\bsDelta\},\ldots,\{\bst_n+\bsDelta\}))^2\,\rd\bsDelta)^{1/2}$ is known as
the \emph{shift-averaged worst case error}.

The idea is then to construct QMC points $\bst_1,\ldots,\bst_n$ to make
the shift-averaged worst case error as small as possible. In practice, we
can take a number of independent random shifts
$\bsDelta_1,\ldots,\bsDelta_r$ drawn from the uniform distribution on
$[0,1]^s$ and use the average
\[
 Q\ran(F)
 \,=\,
 \frac1r \sum_{k=1}^r Q(F;\bsDelta_k)
\]
as the approximation to the integral. Notice our use of the subscript
``ran'' to denote that this is a randomized rule. Since $Q(F;\bsDelta_1),
\ldots, Q(F;\bsDelta_r)$ are i.i.d.\ random variables with mean $I(F)$,
the variance of their average $Q\ran(F)$ is precisely the variance of a
single $Q(F;\bsDelta_k)$ divided by~$r$. This together with
\eqref{eq:sh-avg-wce} gives
\[
  \sqrt{\bbE[|I(F) - Q\ran(F)|^2]}
  \,\le\, r^{-1/2} \, e^\wor_{\rm sh}(\bst_1,\ldots,\bst_n)\,\|F\|_{\calH}
 \qquad\mbox{for all}\quad F\in\calH\,.
\]
A practical estimate of the standard error can be obtained by calculating
\begin{align*}
 \sqrt{\frac{1}{r(r-1)} \sum_{k=1}^r (Q(F;\bsDelta_k) - Q\ran(F))^2}\,,
\end{align*}
from which a confidence interval for $Q\ran(F)$ can be deduced.

A QMC approximation to an integral which is formulated over the Euclidean
space $\bbR^s$ can be done by first mapping the integral to the unit cube
as follows:
\begin{align} \label{eq:int-Rs}
  I(F) \,=\, \int_{\bbR^s} F(\bsy)\,\prod_{j=1}^s \phi(y_j)\,\rd\bsy
  &\,=\, \int_{[0,1]^s} F(\Phi^{-1}(\bsw))\,\rd\bsw \\
  &\,\approx\, \frac{1}{n} \sum_{i=1}^n F(\Phi^{-1}(\bst_i)) \,=\, Q(F) \,. \nonumber
\end{align}
Here $\phi$ can be any general univariate probability density function,
and $\Phi^{-1}$ denotes the component-wise application of the inverse of
the cumulative distribution function corresponding to $\phi$. Note that in
many practical applications we need to first apply some clever
transformation to formulate the integral in the above form; some examples
are discussed in \cite{NK14}.

\subsection{FE discretization}

In the variational setting, we consider the Sobolev space $V = H^1_0(D)$
of functions with vanishing trace on the boundary, with norm $\|v\|_V :=
\|\nabla v\|_{L^2}$, together with its dual space $V^* = H^{-1}(D)$ and
pivot space $L^2(D)$. We now discuss the weak formulation
of~\eqref{eq:strong}. For $f\in V^*$ and $\bsy\in U$ (or $\bsy\in
U_\bsbeta$ in the lognormal case), find $u(\cdot,\bsy)\in V$ such that
\begin{align} \label{eq:weak}
  \scA(\bsy;u(\cdot,\bsy), v) \,=\, \langle f,v\rangle
  \quad\mbox{for all}\quad v \in V\,,
\end{align}
where the parametric bilinear form is given by
\begin{align} \label{eq:bilinear}
  \scA(\bsy; w,v) \,:=\,
  \int_D a(\bsx,\bsy)\,\nabla w(\bsx)\cdot\nabla v(\bsx)\,\rd\bsx
  \quad\mbox{for all}\quad w, v\in V\,,
\end{align}
and $\langle\cdot,\cdot\rangle$ denotes the duality pairing between $V$
and $V^*$.

In the uniform case under Assumptions~\ref{U1} and~\ref{U2}, it
follows that for all $\bsy \in U$ the bilinear form is continuous and
coercive on $V\times V$, and we may infer from the Lax--Milgram lemma that
for every $f\in V^*$ and every $\bsy\in U$ there exists a unique solution
$u(\cdot,\bsy)\in V$ to~\eqref{eq:weak} satisfying the standard \emph{a
priori} estimate
\begin{align} \label{eq:apriori-unif}
 \|u(\cdot,\bsy)\|_V \,\le\, \frac{\|f\|_{V^*}}{a_{\min}}\,.
\end{align}
In addition, if Assumption~\ref{U3} holds and if we assume that the
representer $f$ of the right-hand side of \eqref{eq:weak} is in $L^2(D)$,
then we have
\begin{align} \label{eq:delta-unif}
  \|\Delta u(\cdot,\bsy)\|_{L^2} \,\lesssim\, \|f\|_{L^2}\,.
\end{align}
Throughout this article, the notation $P \lesssim Q$ indicates $P\le C\,Q$
for some constant $C>0$ which is independent of all relevant parameters.

We denote by $\{V_h\}_{h>0}$ a family of subspaces $V_h \subset V$ of
finite dimension $M_h$. For example, $V_h$ can be the space of continuous
piecewise linear finite elements on a sequence of regular triangulations
of $D$ with meshwidth $h>0$. We define the parametric finite element (FE)
approximation as follows: for $f\in V^*$ and $\bsy\in U$ (or $\bsy\in
U_{\overline\bsbeta}$ in the lognormal case), find $u_h(\cdot,\bsy)\in
V_h$ such that
\[
  \scA(\bsy;u_h(\cdot,\bsy), v_h) \,=\, \langle f,v_h\rangle
  \quad\mbox{for all}\quad v_h\in V_h\,.
\]
Then, in the uniform case under Assumptions~\ref{U1} and~\ref{U2}, it is
known that the FE approximation $u_h(\cdot,\bsy)$ of $u(\cdot,\bsy)$ is
stable, that is, \eqref{eq:apriori-unif} holds with $u(\cdot,\bsy)$
replaced by $u_h(\cdot,\bsy)$. Recall that $G\in V^*$ is the linear
functional considered in \eqref{eq:int-unif} and \eqref{eq:int-logn}. In
the same way as we did with $f\in V^*$, we use the same notation to denote
the representer of $G$. In addition, if Assumption~\ref{U3} holds and
$f\in L^2(D)$ and $G\in L^2(D)$, then as $h\to 0$ we have
\begin{align}
  \|u(\cdot,\bsy) - u_h(\cdot,\bsy)\|_V
  &\,\lesssim\, h\,\|\Delta u(\cdot,\bsy)\|_{L^2} \,\lesssim\, h\,\|f\|_{L^2}\,, \label{eq:FEerror-unif}
  \\
  |G(u(\cdot,\bsy)) - G(u_h(\cdot,\bsy))| &\,\lesssim\, h^2\, \|f\|_{L^2}\, \|G\|_{L^2}\,, \label{eq:duality-unif}
  \\
  |I(G(u)) - I(G(u_h))| &\,\lesssim\, h^2\, \|f\|_{L^2}\, \|G\|_{L^2}\,. \label{eq:IGuh-unif}
\end{align}

In the lognormal case under Assumptions~\ref{L1} and~\ref{L2}, the \emph{a
priori} estimate~\eqref{eq:apriori-unif} is replaced by
\begin{align} \label{eq:apriori-logn}
  \|u(\cdot,\bsy)\|_V \,\le\, \frac{\|f\|_{V^*}}{a_{\min}(\bsy)}
  \quad\mbox{for all}\quad
  \bsy\in U_\bsbeta\,.
\end{align}
Adding also Assumption~\ref{L3} and $f\in L^2(D)$, the
bound~\eqref{eq:delta-unif} is replaced by
\[
  \|\Delta u(\cdot,\bsy)\|_{L^2} \,\lesssim\, T(\bsy)\,\|f\|_{L^2}
  \quad\mbox{for all}\quad
  \bsy\in U_{\overline\bsbeta}\,,
\]
where
\begin{align} \label{eq:defT}
  T(\bsy) \,:=\, \frac{\|\nabla a(\cdot,\bsy)\|_{L^\infty}}{a_{\min}^2(\bsy)}
  + \frac{1}{a_{\min}(\bsy)}\,,
\end{align}
while~\eqref{eq:FEerror-unif} is generalized to
\begin{align} \label{eq:FEerror-logn}
  \|a^{1/2}(\cdot,\bsy)\,\nabla (u- u_h)(\cdot,\bsy)\|_{L^2}
  \,\lesssim\, h\,\, a_{\max}^{1/2}(\bsy)\, \|\Delta u(\cdot,\bsy)\|_{L^2}
  \quad\mbox{for all}\quad
  \bsy\in U_{\overline\bsbeta}\,,
\end{align}
see, e.g., \cite{KSSSU}.
Furthermore, if $G\in L^2(D)$ then we also obtain an analogous result
to~\eqref{eq:IGuh-unif} for the lognormal case, see, e.g.,
\cite[Theorem~6]{GKNSSS15}.

To allow for the analysis of higher order methods, we need to assume that
we are given scales of smoothness spaces $\{\calX_t\}_{t\ge 0}$ in the
spatial domain. For example, $\calX_0 = H^1_0(D)$ and $\calX_1 = (H^2\cup
H^1_0)(D)$. For higher order regularity we may consider ``weighted Sobolev
spaces of Kondratiev type'', see, e.g., \cite{NS13}. Then we may consider
families of finite dimensional subspaces $\{\calX^h\}_{h>0} \subset
\calX_0$ and use higher order FE methods to achieve, in the uniform case,
\begin{align}
  \|u(\cdot,\bsy) - u_h(\cdot,\bsy)\|_{\calX_0}
  &\,\lesssim\, h^t \|f\|_{\calX_t^*}\,, \nonumber
  \\
  |G(u(\cdot,\bsy)) - G(u_h(\cdot,\bsy))|
  &\,\lesssim\, h^{t+t'}\, \|f\|_{\calX_t^*}\, \|G\|_{\calX_{t'}^*}\,. \label{eq:FE-ho}
\end{align}

\subsection{Dimension truncation}

For both the uniform and lognormal cases, we observe that truncating the
infinite sum in~\eqref{eq:axy-unif} and~\eqref{eq:axy-logn} to $s$ terms
is the same as setting $y_j = 0$ for $j>s$. We denote the corresponding
weak solution for the truncated case of~\eqref{eq:weak} by $u^s(\bsx,\bsy)
:= u(\bsx,(\bsy_{\{1:s\}};\bszero))$. Throughout this article we refer to
the value of $s$ as the ``truncation dimension''.

In the uniform case under Assumptions~\ref{U1} and~\ref{U2}, for every
$f\in V^*$, every $G\in V^*$, every $\bsy\in U$ and every $s\in\bbN$, we
have from \cite[Theorem~5.1]{KSS12} that, with $b_j$ defined in
\eqref{eq:bj-unif},
\begin{align}
  \| u(\cdot,\bsy) - u^s(\cdot,\bsy) \|_V
  &\,\lesssim\, \|f\|_{V^*}\,\sum_{j\ge s+1} b_j\,, \nonumber
  \\
  |I(G(u)) - I(G(u^s))|
  &\,\lesssim\, \|f\|_{V^*}\|G\|_{V^*}
  \bigg(\sum_{j\ge s+1}b_j\bigg)^2\,. \label{eq:trunc-unif}
\end{align}
In addition, if Assumptions~\textnormal{\ref{U4}}
and~\textnormal{\ref{U5}} hold, then
\begin{align} \label{eq:stechkin}
  \sum_{j\ge s+1} b_j
  \,\le\,
  \min\left(\frac{1}{1/p_0-1},1\right)
  \bigg(\sum_{j\ge1} b_j^{p_0} \bigg)^{1/p_0}
  s^{-(1/p_0-1)}\,.
\end{align}

Dimension truncation analysis in the lognormal case is more complicated.
We summarize here the results from \cite{GKNSSS15} which make use of
\cite{Cha10}. Under Assumption~\ref{L4}, we know from
\cite[Theorem~8]{GKNSSS15} that for $f\in V^*$, $G\in V^*$, $s\in\bbN$ and
$h>0$,
\begin{align}\label{eq:trunc-logn}
  |I(G(u_h)) - I(G(u^s_h))|
  \,\lesssim\, \|f\|_{V^*}\|G\|_{V^*}\,
  s^{-\chi}\,, \quad
  0 < \chi < (\tfrac{1}{2} - \varepsilon)\Theta - \tfrac{1}{2}\,.
\end{align}
For the Mat\'ern covariance~\eqref{eq:matern} with $\nu>d/2$, we know from
\cite[Proposition~9]{GKNSSS15} that Assumption~\ref{L4} holds with $\Theta
= 1 + 2\nu/d$ and $\varepsilon \in
(\frac{1}{2\Theta},\frac{\Theta-1}{2\Theta})$. This implies
that~\eqref{eq:trunc-logn} holds for all $0 < \chi < \nu/d - 1/2$.

\section{Single-level versus multi-level algorithms} \label{sec:algo}

\subsection{Single-level algorithms}

We approximate the integral~\eqref{eq:int-unif} or~\eqref{eq:int-logn} in
three steps:
\begin{enumerate}
\item Dimension truncation: the infinite sum in~\eqref{eq:axy-unif}
or~\eqref{eq:axy-logn} is truncated to $s$ terms.
\item FE discretization: the PDE in weak formulation~\eqref{eq:weak}
    is solved using the piecewise linear FE method with meshwidth $h$.
\item QMC quadrature: the integral of the FE solution for the
    truncated problem is estimated using a deterministic or randomized
    QMC method.
\end{enumerate}

The deterministic version of this algorithm is therefore
\begin{align*}
  A\SL\determ(G(u))
  &\,:=\,
  Q(G(u_h^s))
  =
  \frac{1}{n} \sum_{i=1}^{n} G(u^s_h(\cdot,\bsy_i))\,,
  &
  \bsy_i
  &\,=\,
  \begin{cases}
  \bst_i - \bshalf & \mbox{for uniform}, \\
  \Phi^{\mbox{-}1}(\bst_i) & \mbox{for lognormal},
  \end{cases}
\end{align*}
where $\bst_1,\ldots,\bst_n \in [0,1]^s$ are $n$ QMC points from the
$s$-dimensional standard unit cube. In the uniform case, these points are
translated to the unit cube $[-\tfrac{1}{2},\frac{1}{2}]^s$. In the
lognormal case, these points are mapped to the Euclidean space $\bbR^s$ by
applying the inverse of the cumulative normal distribution function
component-wise.

A randomized version of this algorithm with random shifting is then given by
\begin{multline*}
 A\SL\ran(G(u))
 \,:=\,
 Q\ran(G(u_h^s))
 \,=\,
  \frac{1}{r} \sum_{k=1}^{r}
 \frac{1}{n} \sum_{i=1}^{n} G(u^s_h(\cdot,\bsy_{i,k}))\,,
  \\
  \bsy_{i,k}
  \,=\,
  \begin{cases}
  \{\bst_i + \bsDelta_k\} - \bshalf & \mbox{for uniform}, \\
  \Phi^{\mbox{-}1}(\{\bst_i + \bsDelta_k\}) & \mbox{for lognormal},
  \end{cases}
\end{multline*}
where $\bst_1,\ldots,\bst_n \in [0,1]^s$ are $n$ QMC points as above, and
$\bsDelta_1,\ldots,\bsDelta_r \in [0,1]^s$ are $r$ independent random
shifts generated from the uniform distribution on $[0,1]^s$. Recall that
the braces in $\{\bst_i + \bsDelta_k\}$ mean that we take the fractional
part of each component in the vector $\bst_i + \bsDelta_k$. The total
number of evaluations of the integrand is $r\,n$.

We assume that the cost for assembling the stiffness matrix is $\calO(s \,
h^{-d})$ operations, and further assume that this is higher than the FE
solve itself. Thus the overall cost for the deterministic algorithm is
$\calO(n\,s\,h^{-d})$ operations, while for the randomized algorithm it is
$\calO(r\,n\,s\,h^{-d})$ operations. In practice we assume that $r$ is a
fixed small constant, e.g., $r = 10$ or $20$.

We sometimes refer to these algorithms as ``single-level'' algorithms, in
contrast to ``multi-level'' algorithms to be discussed next.

\subsection{Multi-level algorithms}

The idea of multi-level algorithms in MC simulation originated from
\cite{Hei98,Hei01} and was reinvented in \cite{Gil07,Gil08}, see also
\cite{Gil15}. The general concept is quite easy to explain: if we denote
the integral~\eqref{eq:int-unif} or~\eqref{eq:int-logn} by $I_\infty$ and
define a sequence $I_0, I_1, \ldots$ of approximations converging to
$I_\infty$, then we can write $I_\infty$ as a telescoping sum
\begin{align*}
  I_\infty
  \,=\, (I_\infty - I_L) + \sum_{\ell=0}^L (I_\ell - I_{\ell-1}),
  \qquad I_{-1}:=0\,,
\end{align*}
and then apply different quadrature rules to the differences $I_\ell -
I_{\ell-1}$, which we anticipate to get smaller as $\ell$ increases.

In our case, for each $\ell\ge 0$ we define $I_\ell$ to be the integral of
$G(u^{s_\ell}_{h_\ell})$ corresponding to the FE solution with meshwidth
$h_\ell$ of the truncation problem with $s_\ell$ terms. We assume that
\begin{align*}
  1 \le s_0 \le s_1 \le s_2 \le \cdots \le s_L \le \cdots
  \quad\mbox{and}\quad
  h_0 \ge h_1 \ge h_2 \ge \cdots \ge h_L \ge \cdots > 0\,,
\end{align*}
so that $I_\ell$ becomes a better approximation to $I_\infty$ as $\ell$
increases. For convenience we take $h_\ell \asymp 2^{-\ell}$. Throughout
this article, the notation $P \asymp Q$ means that we have $P\lesssim Q$
and $Q\lesssim P$.

The deterministic version of our multi-level algorithm takes the form
(remembering the linearity of $G$)
\begin{multline*}
  A\ML\determ(G(u))
  \,:=\,
  \sum_{\ell=0}^L \bigg(\frac{1}{n_\ell} \sum_{i=1}^{n_\ell}
                 G((u^{s_\ell}_{h_\ell}-u^{s_{\ell-1}}_{h_{\ell-1}})
                 (\cdot,\bsy_i^\ell))\bigg)\,,
  \quad 
  \bsy_i^\ell \,=\,
  \begin{cases}
  \bst_i^\ell - \bshalf & \mbox{for uniform}, \\
  \Phi^{\mbox{-}1}(\bst_i^\ell) & \mbox{for lognormal},
  \end{cases}
\end{multline*}
where we apply an $s_\ell$-dimensional QMC rule with $n_\ell$ points
$\bst_1^\ell,\ldots,\bst_{n_\ell}^\ell \in [0,1]^{s_\ell}$ to the
integrand $G(u^{s_\ell}_{h_\ell}-u^{s_{\ell-1}}_{h_{\ell-1}})$, and we
define $u^{s_{-1}}_{h_{-1}} :=0$. The total number of evaluations of the
integrand is $\calO(\sum_{\ell=0}^L n_\ell)$.

We can also use $r_\ell$ random shifts at each level, noting that the
shifts should all be independent. Then a randomized version of our
multi-level algorithm with random shifting takes the form
\begin{multline*}
 A\ML\ran(G(u))
 \,:=\,
 \sum_{\ell=0}^L \bigg(\frac{1}{r_\ell} \sum_{k=1}^{r_\ell}
 \frac{1}{n_\ell} \sum_{i=1}^{n_\ell} G((u^{s_\ell}_{h_\ell}-u^{s_{\ell-1}}_{h_{\ell-1}})
 (\cdot,\bsy_{i,k}^\ell))\bigg)\,,
 \\
 \bsy_{i,k}^\ell
 \,=\,
  \begin{cases}
  \{\bst_i^\ell + \bsDelta_k^\ell\} - \bshalf & \mbox{for uniform}, \\
  \Phi^{\mbox{-}1}(\{\bst_i^\ell + \bsDelta_k^\ell\}) & \mbox{for lognormal}.
  \end{cases}
\end{multline*}
In this case the total number of evaluations of the integrand is
$\calO(\sum_{\ell=0}^L r_\ell\, n_\ell)$.

We assume that the overall cost for the multi-level algorithm is
$\calO(\sum_{\ell=0}^L n_\ell\, s_\ell\, h_\ell^{-d})$ operations for the
deterministic version and $\calO(\sum_{\ell=0}^L r_\ell\, n_\ell\,
s_\ell\, h_\ell^{-d})$ operations for the randomized version.

\section{First order versus higher order methods} \label{sec:qmc}

Contemporary analysis of QMC methods is often carried out in the setting
of \emph{weighted} spaces following \cite{SW98,DSWW04,SWW04}. The general
concept is the observation that in many practical examples, not all the
integration variables are of equal importance, and furthermore, there
could be a difference in importance associated with each different subset
of variables. Under appropriate conditions, it is known that QMC methods
can be constructed with error bounds that are independent of the
dimension. In this section we briefly summarize known results for
\emph{randomly shifted lattice rules} (first order) and \emph{interlaced
polynomial lattice rules} (higher order) in suitably weighted spaces.
Construction of these point sets are surveyed in \cite{Nuy14} and we
provide practical pointers to a software implementation in
Section~\ref{sec:comp}. We note that there usually is a tight bond between
the QMC method and the chosen function space setting.

\subsection{Weighted Sobolev spaces over $[0,1]^s$ and randomly shifted lattice
rules} \label{sec:sob}

An $n$-point lattice rule in the unit cube $[0,1]^s$ is a QMC
method~\eqref{eq:qmc} with points
\begin{align}\label{eq:r1lr}
  \bst_i \,=\, \left\{\frac{i\,\bsz}{n} \right\}
  \,=\, \frac{i\,\bsz\bmod n}{n}
  \,,
  \qquad i = 1,\ldots,n\,,
\end{align}
where $\bsz\in \bbZ^s$ is known as the \emph{generating vector} and the
braces indicate that we take the fractional parts of a vector, as
in~\eqref{eq:qmc-sh}. The quality of a lattice rule is determined by the
choice of the generating vector.

We analyze \emph{randomly shifted lattice rules} in a \emph{weighted} and
\emph{unanchored} Sobolev space which is a Hilbert space containing
functions defined over the unit cube $[0,1]^s$, with square integrable
mixed first (weak) derivatives. The norm is given by
\begin{align} \label{eq:norm1}
  \|F\|_{s,\bsgamma}
  \,=\,
  \Bigg[
  \sum_{\setu\subseteq\{1:s\}}
  \frac{1}{\gamma_\setu}
  \int_{[0,1]^{|\setu|}}
  \bigg(\int_{[0,1]^{s-|\setu|}}
  \frac{\partial^{|\setu|}F}{\partial \bsy_\setu}(\bsy_\setu;\bsy_{\{1:s\}\setminus\setu})
  \,\rd\bsy_{\{1:s\}\setminus\setu}
  \bigg)^2
  \rd\bsy_\setu
  \Bigg]^{1/2},
\end{align}
where $\{1:s\}$ is a shorthand notation for the set of indices
$\{1,2,\ldots,s\}$, $(\partial^{|\setu|}F)/(\partial \bsy_\setu)$ denotes
the mixed first derivative of $F$ with respect to the ``active'' variables
$\bsy_\setu = (y_j)_{j\in\setu}$, while $\bsy_{\{1:s\}\setminus\setu} =
(y_j)_{j\in\{1:s\}\setminus\setu}$ denotes the ``inactive'' variables.
This norm is said to be ``unanchored'' because the inactive variables are
integrated out, as opposed to being ``anchored'' at some fixed value, say,
$0$.

There is a weight parameter $\gamma_\setu\ge 0$ associated with each
subset of variables $\bsy_\setu$. A small $\gamma_\setu$ means that $F$
depends weakly on the set of variables $\bsy_\setu$. If $\gamma_\setu = 0$
then it is understood that the corresponding mixed first derivative is
also zero, and then the convention $0/0=0$ is used. We denote by
$\bsgamma$ the set of all weights $\gamma_\setu$, and we take
$\gamma_\emptyset = 1$. There are in general $2^s$ weights in $s$
dimensions, far too many for practical purposes. Special forms of weights
have been considered in the literature, including the so-called ``product
weights'' and ``order-dependent weights'', see, e.g.,
\cite{SW98,DSWW04,SWW04}. Later we will show that a hybrid of these two
forms of weights, called ``product and order dependent weights" or ``POD
weights" for short, naturally arise in the context of PDE applications.
They take the form
\begin{equation} \label{eq:POD}
  \gamma_\setu \,=\, \Gamma_{|\setu|}\,\prod_{j\in\setu} \Upsilon_j\,,
\end{equation}
which is specified by two sequences $\Upsilon_1\ge\Upsilon_2\ge \cdots
> 0$ and $\Gamma_0= \Gamma_1=1, \Gamma_2,\Gamma_3,\ldots\ge 0$.
In this context the cardinality $|\setu|$ of the set $\setu$ is commonly
referred to as the ``order''. Here the factor $\Gamma_{|\setu|}$ is said
to be order dependent because it is determined solely by the cardinality
of $\setu$ and not the precise indices in $\setu$. The dependence of the
weight $\gamma_\setu$ on the indices $j\in\setu$ is controlled by the
product of terms $\Upsilon_j$. Each term $\Upsilon_j$ in the sequence
corresponds to one coordinate direction; the sequence being non-increasing
indicates that successive coordinate directions become less important.

For randomly shifted lattice rules in the unanchored Sobolev space, we
have the root-mean-square error bound~\eqref{eq:sh-avg-wce} where an
explicit expression for the shift-averaged worst case error is known,
allowing it to be analyzed in theory and computed in practice. It has been
proved that a good generating vector $\bsz$ for an $n$-point rule can be
constructed to achieve the optimal convergence rate of
$\calO(n^{-1+\delta})$, $\delta>0$, and the implied constant can be
independent of the dimension $s$ under appropriate conditions on the
weights~$\bsgamma$. The construction is by a \emph{component-by-component}
(\emph{CBC}) \emph{algorithm}: the components of the generating vector
$\bsz$ are obtained one at a time while keeping previously chosen
components fixed. \emph{Fast} CBC algorithms (using FFT) can construct an
$n$-point rule in $s$ dimensions in $\calO(s\,n\log n)$ operations in the
case of product weights, and in $\calO(s\,n\log n + s^2\,n)$ operations in
the case of POD weights.

We summarize the error bound in the theorem below. In the following,
$\zeta(x) := \sum_{k=1}^\infty k^{-x}$ denotes the Riemann zeta function.

\begin{theorem} \label{thm:QMC1}
Let $F$ belong to the unanchored Sobolev space defined over $[0,1]^s$
with weights $\bsgamma$. A randomly shifted lattice rule
with $n=2^m$~points in $s$~dimensions can be constructed by a CBC algorithm
such that for $r$ independent shifts and for all $\lambda\in (1/2,1]$,
\begin{align*}
  \sqrt{\bbE\left[ |I(F) - Q\ran(F)|^2 \right]}
  &\,\le\,
  \frac{1}{\sqrt{r}}
  \Bigg(
  \frac2n
  \sum_{\emptyset\ne\setu\subseteq\{1:s\}} \gamma_\setu^\lambda\,
  [\varrho(\lambda)]^{|\setu|}
  \Bigg)^{1/(2\lambda)}
  \,\|F\|_{s,\bsgamma}\,,
\end{align*}
where
\[
  \varrho(\lambda)
  \,=\,\frac{2\zeta(2\lambda)}{(2\pi^2)^\lambda}
  \,.
\]
\end{theorem}

In the theorem above we have restricted ourselves to the case where $n$ is
a power of $2$, as this is the most convenient setting for generating the
points on a computer. For general $n$, the factor $2/n$ should be replaced
by $1/\varphi_{\rm tot}(n)$, where $\varphi_{\rm tot}(n)$ is the Euler
totient function, i.e., the number of integers between $1$ and $n-1$ that
are relatively prime to $n$. When $n$ is a power of a prime, we have
$1/\varphi_{\rm tot}(n)\le 2/n$ and hence the theorem as stated also holds
in this case.

The best rate of convergence clearly comes from choosing $\lambda$ close
to $1/2$, but the advantage is offset by the fact that
$\zeta(2\lambda)\to\infty$ as $\lambda\to (1/2)_+$.

The CBC construction yields a lattice rule which is ``extensible'' in
dimension $s$, meaning that a generating vector constructed for dimension
$s$ can be used in lower dimensions by taking only the initial components,
and that components for higher dimensions can be appended to existing
components by continuing with the construction.

It is also possible to construct ``lattice sequences'' which are
extensible or embedded in the number of points~$n$, meaning that the same
generating vector can be used to generate more points or less points
without having to construct the existing points anew. This extensibility
in $n$ can be achieved at the expense of increasing the implied constant
in the error bound, see, e.g., \cite{CKN06,DPW08}.

\subsection{Weighted space setting in $\bbR^s$ and randomly shifted lattice
rules} \label{sec:Rs}

For an integral formulated over the Euclidean space $\bbR^s$ as
in~\eqref{eq:int-Rs}, the transformed integrand $F\circ \Phi^{-1}$ arising
from practical applications typically does not belong to the Sobolev space
defined over the unit cube due to the integrand being unbounded near the
boundary of the cube, or because the mixed derivatives of the transformed
integrand do not exist or are unbounded. Thus most QMC theories, including
Theorem~\ref{thm:QMC1}, generally do not apply in practice. Here we
summarize a special weighted space setting in $\bbR^s$ for which randomly
shifted lattice rules have been shown to achieve nearly the optimal
convergence rate of order one, see \cite{KSWWat10,NK14}. The norm in this
setting is given by
\begin{align} \label{eq:norm2}
  \|F\|_{s,\bsgamma}
  &\,=\,
  \Bigg[
  \sum_{\setu\subseteq\{1:s\}} \frac{1}{\gamma_\setu}
  \int_{\bbR^{|\setu|}}
  \bigg(
  \int_{\bbR^{s-|\setu|}}
  \frac{\partial^{|\setu|} F}{\partial \bsy_\setu}(\bsy_\setu;\bsy_{\{1:s\}\setminus\setu})
  \bigg(\prod_{j\in\{1:s\}\setminus\setu} \phi(y_j)\bigg)
  \,\rd\bsy_{\{1:s\}\setminus\setu}
  \bigg)^2 \nonumber\\
  &\qquad\qquad\qquad\qquad\qquad\qquad\qquad\qquad\;\;\times
  \bigg(\prod_{j\in\setu} \vpsi_j^2(y_j)\bigg)
  \,\rd\bsy_\setu \Bigg]^{1/2} \,.
\end{align}
Comparing~\eqref{eq:norm2} with~\eqref{eq:norm1}, apart from the
difference that the integrals are now over the unbounded domain, there is
a probability density function $\phi$ as well as additional \emph{weight
functions} $\vpsi_j$ which can be chosen to reflect the boundary behavior
of the mixed derivatives of $F$.

In the context of PDEs with lognormal random coefficients, $\phi$ is the
standard normal density. To ensure that the integrands from the PDE
problems belong to the function space setting, we may restrict ourselves
to the choice (see \eqref{eq:hindsight} ahead)
\begin{equation} \label{eq:psi}
  \vpsi_j^2(y_j) \,=\, \exp (-2\,\alpha_j\, |y_j|)\;,
  \qquad\alpha_j > 0\,.
\end{equation}
We have the following result from \cite[Theorem~15]{GKNSSS15} for randomly
shifted lattice rules in this setting; see \cite[Theorem~8]{NK14} for
results on general functions $\phi$ and $\vpsi_j$. We state the result
again for $n$ a power of~$2$, but it holds when $n$ is a power of any
prime.

\begin{theorem} \label{thm:QMC2}
Let $F$ belong to the weighted function space over $\bbR^s$ with
weights~$\bsgamma$, with $\phi$ being the standard normal density, and
with weight functions $\vpsi_j$ given by~\eqref{eq:psi}. A randomly
shifted lattice rule with $n=2^m$~points in $s$~dimensions can be
constructed by a CBC algorithm such that for $r$ independent shifts and
for all $\lambda\in (1/2,1]$,
\begin{align*}
  \sqrt{\bbE\left[ |I(F) - Q\ran(F)|^2 \right]}
  \,\le\,
  \frac{1}{\sqrt{r}}
  \Bigg(
  \frac2n
  \sum_{\emptyset\ne \setu\subseteq\{1:s\}}
  \gamma_\setu^\lambda\, \prod_{j\in\setu} \varrho_j(\lambda)\Bigg)^{1/(2\lambda)}\,
   \,
   \|F\|_{s,\bsgamma}\,,
\end{align*}
where
\[
  \varrho_j(\lambda)
  \,=\, 2\left(\frac{\sqrt{2\pi}\,\exp(\alpha_j^2/\eta_*)}{\pi^{2-2\eta_*}(1-\eta_*)\eta_*}\right)^\lambda\,
  \zeta\big(\lambda + \tfrac{1}{2}\big),
  \quad\mbox{and}\quad
  \eta_* \,=\, \frac{2\lambda-1}{4\lambda}\,.
\]
\end{theorem}

Alternatively, for the lognormal case it is also possible to choose the
weight functions $\vpsi_j$ to take the form of a normal density that
decays much slower than $\phi$. The corresponding result is given below,
which can be obtained from \cite[Theorem~8]{NK14} together with
\cite[Example~4]{KSWWat10}.

\begin{theorem} \label{thm:QMC2b}
If we replace the weight functions $\vpsi_j$ in Theorem~\ref{thm:QMC2} by
$\vpsi_{\rm alt}$, given by $\vpsi^2_{\rm alt}(y_j) = \exp (- \alpha\,
y_j^2)$ with $\alpha < 1/2$. Then the root-mean-square error bound in
Theorem~\ref{thm:QMC2} holds for all $\lambda\in (1/(2-2\alpha),1]$, but
with $\rho_j(\lambda)$ replaced by
\[
  \varrho_{\rm alt}(\lambda)
  \,=\, 2\left(\frac{\sqrt{2\pi}}{\pi^{2-2\alpha}(1-\alpha)\alpha}\right)^\lambda\,
  \zeta\big(2(1-\alpha)\lambda\big)\,.
\]
\end{theorem}

\subsection{Weighted space of smooth functions over $[0,1]^s$ and interlaced polynomial lattice rules}

We now introduce \emph{interlaced polynomial lattice rules} which are a
special family of higher order QMC rules, see, e.g.,
\cite{God15,GD15,DKLNS14} and \cite{Dic07,Dic08}. We limit ourselves to
the case where the polynomials are over the finite field $\bbZ_2$ with two
elements. This has major advantages in the computer implementations for
constructing these rules and generating the points, and simplifies the
presentation.

Let $\bbN_0 = \{0,1,2,\ldots\}$. We will identify numbers $x \in [0,1)$
and $i \in \bbN_0$ that have finite base 2 representations (i.e., having
only a finite number of bits being $1$) by a vector enumerating the bits,
denoted by $\vvec{x}$ and $\vvec{i}$, or by a polynomial in the formal
variable $\X$, denoted by $x(\X)$ and $i(\X)$.
We will use the common notation $\bbZ_2[\X]$ for formal
polynomials in $\X$ and $\bbZ_2((\X^{-1}))$ for formal
polynomials in $\X^{-1}$ as well as $\bbZ_2((\X))$ for
formal Laurent series (with powers going in both
directions).
Then, for $x \in [0,1)$, with
$x_\ell \in \bbZ_2$ for $\ell \ge 1$ denoting the base~$2$ digits of~$x$,
we have in this multitude of notations
\begin{align*}
  x \,=\, \sum_{\ell \ge 1} x_\ell \, 2^{-\ell}
    \,=\, ( 0.x_1 x_2 \ldots )_2
  \in \bbR
  &\;\simeq\;
  ( x_1, x_2, \ldots )^\top
  &&\,=:\,\vvec{x}
  \in \bbZ_2^\infty
  \\&\;\simeq\;
  x_1 \, \X^{-1} + x_2 \, \X^{-2} + \cdots
  &&\,=:\,
  x(\X)
  \in \bbZ_2((\X^{-1}))
  \,.
\end{align*}
Similarly for integers $i \in \bbN_0$ where the bit expansion goes in the
other direction, with $i_\ell \in \bbZ_2$ for $\ell \ge 0$ denoting the
base~$2$ digits of~$i$, we have
\begin{align*}
  i \,=\, \sum_{\ell \ge 0} i_\ell \, 2^\ell
    \,=\, ( \ldots i_2 i_1 i_0 )_2
  \in \bbN_0
  &\;\simeq\;
  ( i_0, i_1, i_2, \ldots )^\top
  &&\,=:\,
  \vvec{i}
  \in \bbZ_2^\infty
  \\&\;\simeq\;
  i_0 \, \X^0 + i_1 \, \X^1 + i_2 \, \X^2 + \cdots
  &&\,=:\,
  i(\X)
  \in \bbZ_2[\X]
  \,.
\end{align*}
We need one further notation to limit the number of bits to~$m$; this will
be denoted by $[x(\X)]_m$ and $[\vvec{x}]_m$, or $[i(\X)]_m$ and
$[\vvec{i}]_m$.

With this notation in place, the points of a \emph{polynomial lattice
rule} are given by
\begin{align}\label{eq:pr1lr}
  \bst_i
  &\,\simeq\,
  [ \bstau_i(\X) ]_m
  \,,
  &
  \bstau_i(\X)
  &\,=\,
  \frac{ i(\X) \, \bsz(\X) \bmod{P(\X)}}{P(\X)}
  \in \bigl(\bbZ_2((\X^{-1}))\bigr)^s
  \,,
  &
  i
  &\,=\,
  0,\ldots,2^m-1
  \,,
\end{align}
where $P(\X)$ is an irreducible polynomial of degree $m$ over $\bbZ_2$,
known as the \emph{modulus}, and the vector of polynomials $\bsz(\X) =
(z_1(\X), \ldots, z_s(\X)) \in \bbZ[\X]^s$ is known as the
\emph{generating vector}. We see that polynomial lattice rules take the
same form as lattice rules in \eqref{eq:r1lr}, but the integers are
replaced by polynomials, and thus the multiplication and division are
replaced by their polynomial equivalents over $\bbZ_2((\X))$. We remark
that the points are here indexed from $0$ to $2^m-1$ which is different
from the convention we used in the rest of this article.

Alternatively, the polynomial multiplication and division can be written
in matrix-vector notation over $\bbZ_2$ by identifying a Hankel matrix
$C_j = ( a_{j,r+t-1} )_{r,t\ge1} $ with the division $z_j(\X) / P(\X) =
\sum_{\ell} a_{j,\ell} \, \X^{-\ell}$ for $j=1,\ldots,s$. In principle
these matrices are in $\bbZ_2^{\infty \times \infty}$ but we restrict them
to be of size $m \times m$. In terms of these ``generating matrices''
$C_1,\ldots,C_s \in \bbZ_2^{m\times m}$ we can then write the points of
the polynomial lattice rule as
\begin{align}\label{eq:pr1lr-matrix}
  \bst_i
  &\,\simeq\,
  ( C_1 \, [\vvec{i}]_m , \ldots, C_s \, [\vvec{i}]_m )^\top
  \,,
  &
  i
  &\,=\,
  0,\ldots,2^m-1
  \,.
\end{align}
The matrix-vector form is to be preferred when generating the actual
points as point generation then boils down to simple bit operations. In
fact, if the order in which the points are iterated is changed into Gray
code ordering so that only one bit changes in the index of the point at a
time, then the coordinate of the next point can be obtained by simply
applying the \textsf{xor}-operation between the previous value of that
coordinate and the particular column from the generator matrix
corresponding to the single bit change in the Gray code.

An \emph{interlaced polynomial lattice rule} in $s$ dimensions with
interlacing factor $\alpha$ is now obtained by taking a polynomial lattice
rule in $\alpha\,s$ dimensions and then interlacing the bits from every
successive $\alpha$ dimensions to yield one dimension. More explicitly,
for $\alpha = 3$, given three coordinates $x = (0.x_1 x_2 \ldots x_m)_2$,
$y = (0.y_1 y_2 \ldots y_m)_2$ and $z = (0.z_1 z_2 \ldots z_m)_2$ we
interlace their bits to obtain $w = (0.x_1 y_1 z_1 x_2 y_2 z_2 \ldots x_m
y_m z_m)_2$. Note that interlacing can be applied to any existing point
set, but the construction of the polynomial lattice rules in this section
actually takes the interlacing into account.

The interlacing operation can also be done directly by constructing new
generator matrices which are obtained by interlacing the rows of $\alpha$
successive generator matrices. That is, if we have generating matrices
$C_1,\ldots,C_{\alpha\,s} \in \bbZ_2^{m \times m}$ then we obtain new
generating matrices $B_1,\ldots,B_s \in \bbZ_2^{\alpha m \times m}$ and
the points of the interlaced polynomial lattice rule are given by
\begin{align}\label{eq:ipr1lr}
  \bst_i
  &\,\simeq\,
  ( B_1 \, [\vvec{i}]_m , \ldots, B_s \, [\vvec{i}]_m )^\top
  \,,
  &
  i
  &\,=\,
  0,\ldots,2^m-1
  \,.
\end{align}
We remark that for efficiency reasons the points are normally iterated
in Gray code ordering such that each coordinate can be obtained by
a single \textsf{xor} operation.

A function space setting for smooth integrands defined over the unit cube
was introduced in \cite{DKLNS14}. Let $\alpha, s\in\bbN$, and $1\le q,r
\le \infty$, and let $\bsgamma = (\gamma_\setu)_{\setu\subset\bbN}$ be a
collection of weights as in the previous subsections. The norm in this
setting for $1\le q,r<\infty$ is given by
\begin{multline*}
 \|F\|_{s,\alpha,\bsgamma,q,r}
 \,:=\,
 \Bigg[ \sum_{\setu\subseteq\{1:s\}} \Bigg( \frac{1}{\gamma_\setu^q}
 \sum_{\setv\subseteq\setu}\,
 \sum_{\bstau_{\setu\setminus\setv} \in \{1:\alpha\}^{|\setu\setminus\setv|}}
 \\
 \int_{[0,1]^{|\setv|}} \bigg|\int_{[0,1]^{s-|\setv|}}
 (\partial^{(\bsalpha_\setv,\bstau_{\setu\setminus\setv},\bszero)} F)(\bsy) \,
 \rd \bsy_{\{1:s\} \setminus\setv}
 \bigg|^q \rd \bsy_\setv \Bigg)^{r/q} \Bigg]^{1/r},
\end{multline*}
with the obvious modifications if $q$ or $r$ is infinite, see
\eqref{eq:norm3} below. Here
$(\bsalpha_\setv,\bstau_{\setu\setminus\setv},\bszero)$ denotes a vector
$\bsnu$ with $\nu_j = \alpha$ for $j\in\setv$, $\nu_j = \tau_j$ for
$j\in\setu\setminus\setv$, and $\nu_j = 0$ for $j\notin\setu$. We denote
the $\bsnu$-th partial derivative of $F$ by $\partial^\bsnu F \,=\,
(\partial^{|\bsnu|}F)/(\partial^{\nu_1}_{y_1}\partial^{\nu_2}_{y_2}\cdots\partial^{\nu_s}_{y_s})$.
We remark that the norm is stated incorrectly in
\cite[Equation~(3.7)]{DKLNS14}.

A new form of weights arose from the analysis of interlaced polynomial
lattice rules in this setting in the context of PDE applications, see
\cite{DKLNS14}. They are called ``smoothness-driven product and order
dependent weights'' or ``SPOD weights'' for short, and take the form
\begin{equation} \label{eq:SPOD}
 \gamma_\setu
 \,=\,
 \sum_{\bsnu_\setu \in \{1:\alpha\}^{|\setu|}}
 \Gamma_{|\bsnu_\setu|}  \prod_{j\in\setu} \Upsilon_j(\nu_j)\,.
\end{equation}

For PDE applications, it was shown in \cite{DKLNS14} that good theoretical
results can be obtained by taking $r=\infty$ while $q$ can be arbitrary.
Therefore, for simplicity, here we take $q=r=\infty$ and denote the
corresponding norm by
\begin{multline}\label{eq:norm3}
  \|F\|_{s,\alpha,\bsgamma}
  \,:=\,
  \|F\|_{s,\alpha,\bsgamma,\infty,\infty}
  \\\,=\,
 \sup_{\setu\subseteq\{1:s\}}
 \sup_{\bsy_\setv \in [0,1]^{|\setv|}}
 \frac{1}{\gamma_\setu}
 \sum_{\setv\subseteq\setu} \,
 \sum_{\bstau_{\setu\setminus\setv} \in \{1:\alpha\}^{|\setu\setminus\setv|}}
 \bigg|\int_{[0,1]^{s-|\setv|}}
 (\partial^{(\bsalpha_\setv,\bstau_{\setu\setminus\setv},\bszero)} F)(\bsy) \,
 \rd \bsy_{\{1:s\} \setminus\setv}
 \bigg|\,
  .
\end{multline}
The following theorem is adjusted from \cite[Theorem~3.10]{DKLNS14} for
$b=2$. We made use of a new result in \cite{Yos15} that a constant
$C_{\alpha,b}$ which usually appears is exactly~$1$ when $b=2$.
Furthermore, we followed the proof of \cite[Theorem~5.1]{Nuy14} to obtain
the factor $2/n$ instead of $2/(n-1)$. We remark that the interlacing
factor $\alpha$ needs to be at least~$2$ for the theorem to hold.

\begin{theorem} \label{thm:QMC3}
Let $F$ belong to the weighted space of smooth functions over $[0,1]^s$
with $\alpha\ge 2$ and weights $\bsgamma$. An interlaced polynomial
lattice rule with interlacing factor $\alpha$, with irreducible modulus
polynomial of degree $m$, and with $n=2^m$ points in $s$ dimensions, can
be constructed by a CBC algorithm such that, for all $\lambda \in
(1/\alpha,1]$,
\begin{align*}
  |I(F) - Q(F)|
  \,\le\,
  \left(\frac{2}{n} \sum_{\emptyset\ne\setu\subseteq{\{1:s\}}}
  \gamma_\setu^\lambda\, [\rho_{\alpha}(\lambda)]^{|\setu|}\right)^{1/\lambda}\,
  \|F\|_{s,\alpha,\bsgamma}\,,
\end{align*}
where
\[
  \rho_{\alpha}(\lambda) \,=\,
  2^{\alpha\lambda(\alpha-1)/2}
  \bigg(\bigg(1+\frac{1}{2^{\alpha\lambda}-2}\bigg)^\alpha-1\bigg)\,.
\]
\end{theorem}

If the weights $\bsgamma$ are SPOD weights, then the CBC algorithm has
cost $\calO(\alpha\,s\, n\log n + \alpha^2\,s^2 n)$ operations. If the
weights $\bsgamma$ are product weights, then the CBC algorithm has cost
$\calO(\alpha\,s\, n\log n)$ operations.


\section{Error analysis} \label{sec:error}

In this section we summarize the error analysis for various algorithms
based on QMC methods in the uniform and lognormal cases. Not surprisingly,
the error is a combination of dimension truncation error, FE
discretization error, and QMC quadrature error. The errors are additive in
the case of single-level algorithms, while in the case of multi-level
algorithms the overall error includes some multiplicative effects between
FE and QMC errors. Recall that bounds on truncation and FE errors are
summarized in Section~\ref{sec:appr}, while bounds on QMC errors for
different spaces are summarized in Section~\ref{sec:qmc}.

The convergence of the QMC method (deterministic or randomized) can be
independent of the truncation dimension both in the rate (i.e., the
exponent of $1/n$ in the error estimate) and in the asymptotic constant.
This is achieved by working in weighted spaces with strategically chosen
weights $\gamma_\setu$. To bound the QMC error we will use
Theorem~\ref{thm:QMC1}, \ref{thm:QMC2}, or~\ref{thm:QMC3} depending on the
setting. From these theorems we see that the key step of our analysis is
to obtain bounds on the particular weighted norm of the integrand, which
depends on the PDE solution $u^s_h(\bsx,\bsy)$. Specifically, this means
that we need to obtain bounds on the mixed partial derivatives of
$u^s_h(\bsx,\bsy)$ with respect to $\bsy$. Once we obtain estimates on the
norm of the integrand, we then choose suitable weights $\gamma_\setu$ to
reduce the error bound, optimizing on the theoretical QMC convergence rate
while ensuring (under some conditions) that the implied constant is
independent of the truncation dimension. The chosen weights $\gamma_\setu$
then enter the fast CBC construction to produce tailored QMC methods for
the given setting.

In the next subsections we outline the error analysis for different
algorithms under different settings. First order results are based on
randomly shifted lattice rules and we obtain probabilistic error bounds.
Higher order results are based on interlaced polynomial lattice rules and
we obtain deterministic error bounds. In the latter case we also replace
the $\ell^2$ norm in the typical definition of the function space norm by
the $\ell^\infty$ norm, and this enables us to gain an extra factor of
$n^{-1/2}$ in the QMC convergence rate in the context of PDE applications.
However, this analysis only applies in the uniform case.

To obtain multiplicative error bounds in the case of multi-level
algorithms, we need to assume a stronger regularity of $u^s_h(\bsx,\bsy)$
in $\bsx$, and we need to establish regularity results simultaneously in
$\bsx$ and $\bsy$. This makes the analysis more challenging. The resulting
weights $\bsgamma_\setu$ are also different from those in the single-level
algorithms.

The error versus cost analysis depends crucially on the cost model
assumption. In the single-level algorithms we choose $n,s,h$ to balance
the three sources of errors. In the multi-level algorithms we choose
$n_\ell,s_\ell,h_\ell$, for each level, to minimize the error for a fixed
cost using Lagrange multiplier arguments, and we choose $L$ such that the
combined error meets the required threshold. We assume that $r$ and
$r_\ell$ are fixed constants.

Before we proceed we introduce some notation. For a multi-index $\bsnu =
(\nu_j)_{j\ge1}$ with $\nu_j\in \{0,1,2,\ldots\}$, we write its ``order''
as $|\bsnu| := \sum_{j\ge 1} \nu_j$ and its ``support'' as $\supp(\bsnu)
:= \{j\ge 1: \nu_j\ge 1\}$. Furthermore, we write $\bsnu! := \prod_{j\ge
1} \nu_j!$, which is different from $|\bsnu|! = (\sum_{j\ge 1} \nu_j)!$.
We denote by $\indx$ the (countable) set of all ``finitely supported''
multi-indices:
\[
  \indx
  \,:=\,
  \{ \bsnu \in \bbN_0^\bbN : \supp(\bsnu) < \infty \}
  .
\]
For $\bsnu\in\indx$, we denote the $\bsnu$-th partial derivative with
respect to the parametric variables $\bsy$ by
\begin{align*}
  \partial^{\bsnu}
  \,=\, \frac{\partial^{|\bsnu|}}{\partial y_1^{\nu_1}\partial y_2^{\nu_2}\cdots} \,.
\end{align*}
For any sequence of real numbers $\bsb = (b_j)_{j\ge 1}$, we write
$\bsb^\bsnu := \prod_{j\ge 1} b_j^{\nu_j}$. By $\bsm\le\bsnu$ we mean that
the multi-index $\bsm$ satisfies $m_j\le \nu_j$ for all $j$. Moreover,
$\bsnu-\bsm$ denotes a multi-index with the elements $\nu_j-m_j$, and
$\binom{\bsnu}{\bsm} := \prod_{j\ge 1} \binom{\nu_j}{m_j}$. We denote by
$\bse_j$ the multi-index whose $j$th component is $1$ and all other
components are $0$.

We remind the reader that throughout this article, the gradient~$\nabla$
and the Laplacian~$\Delta$ are to be taken with respect to the spatial
variables $\bsx$, while the partial derivatives~$\partial^{\bsnu}$ are to
be taken with respect to the parametric variables $\bsy$.

Note that in Section~\ref{sec:qmc} we have used two different notations
for the mixed derivatives with respect to $\bsy$. In the norms
\eqref{eq:norm1} and \eqref{eq:norm2} we restrict only to mixed first
derivatives, and the subsets $\setu\subseteq\{1:s\}$ are used to identify
the indices of the variables with respect to which we differentiate. For
example, if $\setu = \{1,2,5\}$ then
\[
  \frac{\partial^{|\setu|}F}{\partial\bsy_\setu}
  \,=\, \frac{\partial^3F}{\partial y_1\partial y_2\partial y_5}
  \,=\, \partial^{\bsnu}F\,,
\]
with $\bsnu=(1,1,0,0,1,0,0,\ldots)$ in the multi-index notation.


\subsection{First order, single-level, uniform} \label{sec:fo-sl-unif}

The mean-square error for our single-level algorithm with randomly shifted
lattice rules can be expressed as
\begin{align} \label{eq:err-SL-ran}
  \bbE \,[ |(I-A\SL\ran)(G(u))|^2]
  \,=\, | I(G(u-u^s_h)) |^2 \,+\, \bbE\,[ |(I- Q\ran)(G(u^s_h))|^2 ]\,,
\end{align}
where we used the linearity of $G$ and the unbiased property of randomly
shifted QMC rules, i.e., $\bbE\,[Q\ran(F)] = I(F)$. The first term on the
right-hand side of \eqref{eq:err-SL-ran} can be split trivially using
linearity of $I$ and $G$ into
\begin{align} \label{eq:split-s}
  I(G(u-u^s_h)) \,=\, I(G(u-u^s)) \,+\, I(G(u^s-u^s_h))\,.
\end{align}

In the uniform case, we estimate the truncation error using
\eqref{eq:trunc-unif}--\eqref{eq:stechkin}, and estimate the FE error
using \eqref{eq:IGuh-unif} but adapt it to the truncated solutions (this
is valid since \eqref{eq:duality-unif} holds for all $\bsy\in U$,
including those with $y_j = 0$ for $j>s$). For the QMC error we use
Theorem~\ref{thm:QMC1} which requires a bound on the norm
$\|G(u^s_h)\|_{s,\bsgamma}$, and we see from the definition
\eqref{eq:norm1} that we need to bound the mixed first partial derivatives
of $G(u^s_h(\cdot,\bsy))$. Due to linearity and boundedness of~$G$, we
have
\begin{align} \label{eq:lin-G}
 \bigg|\frac{\partial^{|\setu|}}{\partial \bsy_\setu} G(u^{s}_h(\cdot,\bsy))\bigg|
 \,=\, \bigg|G\bigg(\frac{\partial^{|\setu|}}{\partial \bsy_\setu} u^{s}_h(\cdot,\bsy)\bigg)\bigg|
 \,\le\, \|G\|_{V^*}\,\bigg\|\frac{\partial^{|\setu|}}{\partial \bsy_\setu} u^{s}_h(\cdot,\bsy)\bigg\|_V\,.
\end{align}
Hence we need a regularity result on the PDE solution $u^s_h(\bsx,\bsy)$
with respect to the parameters $\bsy$. In Lemma~\ref{lem:regu1} below we
state such a result for general mixed derivatives using the multi-index
notation.

\begin{lemma} \label{lem:regu1}
In the uniform case under Assumptions~\ref{U1} and~\ref{U2}, for every
$f\in V^*$, every $\bsy\in U$, and every $\bsnu\in\indx$, we have
\begin{align*}
  \|\partial^{\bsnu}u(\cdot,\bsy)\|_{V}
  \,=\, \|\nabla (\partial^{\bsnu}u (\cdot,\bsy))\|_{L^2}
 \,\le\, |\bsnu|!\, \bsb^\bsnu\, \frac{\|f\|_{V^*}}{a_{\min}}\,,
\end{align*}
where the sequence $\bsb$ is defined in \eqref{eq:bj-unif}. The same
estimate holds when the exact solution $u$ is replaced by $u^s_h$.
\end{lemma}

This result was proved in \cite{CDS10}. The proof is by induction on
$|\bsnu|$. We take the mixed partial derivative $\partial^{\bsnu}$ on both
sides of the weak formulation \eqref{eq:weak} and then substitute the test
function $v = \partial^{\bsnu}u(\cdot,\bsy)$. Rearranging and estimating
the terms then yields the required bound. Since the same proof strategy is
used repeatedly in subsequent more complicated proofs, we include this
relatively simple proof in Section~\ref{sec:proof} as a first
illustration.

For bounds on mixed first derivatives in the norm \eqref{eq:norm1}, we
restrict Lemma~\ref{lem:regu1} to multi-indices $\bsnu$ with $\nu_j\le 1$
for all $j$. Using also \eqref{eq:lin-G}, we obtain the estimate
\begin{equation} \label{eq:Gu-norm}
  \|G(u^s_h)\|_{s,\bsgamma}
  \,\le\, \frac{\|f\|_{V^*}\,\|G\|_{V^*}}{a_{\min}}
  \Bigg(
  \sum_{\setu\subseteq\{1:s\}} \frac{(|\setu|!)^2 \prod_{j\in\setu} b_j^2}{\gamma_\setu}
  \Bigg)^{1/2}\,.
\end{equation}
Combining this with Theorem~\ref{thm:QMC1} gives
\begin{align*}
  \bbE\left[ |(I-Q_{\rm ran})(G(u^s_h))|^2 \right]
  \,\le\,
  \frac1r
  \frac{C_\bsgamma^2(\lambda)}{n^{1/\lambda}}\, \frac{\|f\|_{V^*}^2\,\|G\|_{V^*}^2}{a_{\min}^2}\,,
\end{align*}
where
\begin{align*}
  C_\bsgamma^2(\lambda)
  \,=\, \Bigg( 2
  \sum_{\setu\subseteq\{1:s\}} \gamma_\setu^\lambda\,
  [\rho(\lambda)]^{|\setu|}
  \Bigg)^{1/\lambda}
  \,
  \Bigg(
  \sum_{\setu\subseteq\{1:s\}} \frac{(|\setu|!)^2 \prod_{j\in\setu} b_j^2}{\gamma_\setu}
  \Bigg)\,.
\end{align*}
Elementary calculus (see \cite[Lemma~6.2]{KSS12}) yields the result that
for any given $\lambda$, $C_\bsgamma^2(\lambda)$ is minimized as a
function of $\gamma_\setu$ by taking
\begin{equation} \label{eq:weight1}
  \gamma_\setu \,=\,
  \bigg(|\setu|!\, \prod_{j\in\setu} \frac{b_j}{\sqrt{\rho(\lambda)}} \bigg)^{2/(1+\lambda)}\,,
\end{equation}
which are of the POD form \eqref{eq:POD}. Under Assumption~\ref{U5}, it is
proved in \cite[Theorem~6.4]{KSS12} that if we take
\begin{equation} \label{eq:lambda-p0}
 \lambda \,=\,
 \begin{cases}
 \displaystyle\frac{1}{2-2\delta} \quad\mbox{for some } \delta\in (0,1/2)
 & \mbox{when } p_0\in (0,2/3]\,, \vspace{0.1cm} \\
 \displaystyle\frac{p_0}{2-p_0} & \mbox{when } p_0\in (2/3,1)\,,
 \end{cases}
\end{equation}
then the constant $C_\bsgamma^2(\lambda)$ can be bounded independently of
$s$. The case $p_0=1$ can also be analyzed but it requires an additional
assumption which can be quite restrictive; we omit this case here.

We summarize the final result in the theorem below.

\begin{theorem} \label{thm:main1}
In the uniform case under Assumptions~\ref{U1}--\ref{U5}, for every $f\in
L^2(D)$ and every $G\in L^2(D)$, the single-level algorithm
$A\SL\ran(G(u))$ using a randomly shifted lattice rule with $n=2^m$ points
constructed from a CBC algorithm with POD weights
\eqref{eq:weight1}--\eqref{eq:lambda-p0}, at the pre-computation cost of
$\calO(s\, n\log n + s^2n)$ operations, achieves the mean-square error
bound
\begin{align*}
  \bbE\left[ |(I-A\SL\ran)(G(u))|^2 \right]
  \,\lesssim\,
  s^{-4(1/p_0-1)} + h^{4} + r^{-1}\,n^{-2 \min(1/p_0-1/2,1-\delta)}\,,
  \; \delta\in (0,\tfrac{1}{2})\,,
\end{align*}
where the implied constant is independent of $s$, $h$, $r$ and $n$.
\end{theorem}

If we treat $r$ as a fixed constant and choose $s$, $h$, $n$ to balance
the terms so that the mean-square error is $\calO(\varepsilon^2)$, then
the cost of the algorithm $A\SL\ran(G(u))$ is $\calO(s\, n\, h^{-d}) =
\calO(\varepsilon^{-\tau})$ with
\[
  \tau \,=\, \frac{p_0}{2-2p_0} + \frac{d}{2}
  + \max\bigg(\frac{2p_0}{2-p_0}, \frac{1}{1-\delta}\bigg) \,,
  \quad \delta\in (0,\tfrac{1}{2})\,.
\]


\subsection{First order, multi-level, uniform} \label{sec:fo-ml-unif}

For our multi-level algorithm with randomly shifted lattice rules, we can
write the mean-square error as
\begin{multline} \label{eq:err-ML-ran}
  \bbE \,[ |(I- A\ML\ran)(G(u))|^2]
  \,=\, | I(G(u-u^{s_L}_{h_L})) |^2
  \,+\, \sum_{\ell=0}^{L}
  \bbE\,[ |(I - Q\ran^\ell)(G(u^{s_\ell}_{h_\ell}-u^{s_{\ell-1}}_{h_{\ell-1}}))|^2 ]\,,
\end{multline}
where $Q\ran^\ell$ denotes a randomly shifted lattice rule in $s_\ell$
dimensions with $n_\ell$ points and $r_\ell$ independent shifts. The first
term on the right-hand side of \eqref{eq:err-ML-ran} can be estimated in
exactly the same way as for the single-level algorithm in the previous
subsection. Each term in the sum over $\ell$ in \eqref{eq:err-ML-ran} can
be estimated using Theorem~\ref{thm:QMC1}. For simplicity we restrict our
discussion here to the case where $s_\ell = s$ for all $\ell$. The case
$s_\ell\ne s_{\ell-1}$ was analyzed in \cite[Theorem~8]{KSS15}. Thus we
need to estimate the norm
\begin{align*}
 \| G(u^s_{h_\ell} - u^s_{h_{\ell-1}})\|_{s,\bsgamma}
 \,\le\, \| G(u^s - u^s_{h_{\ell}})\|_{s,\bsgamma}
 \,+\, \| G(u^s- u^s_{h_{\ell-1}})\|_{s,\bsgamma}\,.
\end{align*}
For this we need to estimate the mixed first derivatives of $G(u-u_h)$
with respect to~$\bsy$. We state a result for general mixed derivatives in
Lemma~\ref{lem:regu2c} below.

To establish Lemma~\ref{lem:regu2c} we need some intermediate results in
Lemmas~\ref{lem:regu2a} and~\ref{lem:regu2b}. These three lemmas together
correspond to \cite[Theorems~6 and~7]{KSS15}, but in addition to a
generalization of \cite[Theorem~7]{KSS15} from first derivatives to
general derivatives, the results are different because the sequence
$\overline{\bsb}$ defined here is simpler compared to \cite{KSS15}, and we
take a different (arguably more direct) route with the proofs. However, we
get bigger factorials here. For simplicity we restrict our discussion to
the case $f,G\in L^2(D)$. The results can be generalized to cover $f\in
H^{-1+t}(D)$ and $G\in H^{-1+t'}(D)$ for $t,t'\in [0,1]$ as
in~\cite{KSS15}.

\begin{lemma} \label{lem:regu2a}
In the uniform case under Assumptions~\ref{U1}--\ref{U3}, for every $f\in
L^2(D)$, every $\bsy\in U$, and every $\bsnu\in\indx$, we have
\begin{align*}
  \|\Delta (\partial^{\bsnu}u (\cdot,\bsy))\|_{L^2}
  \,\lesssim\,  (|\bsnu|+1)!\,\, \overline\bsb^{\bsnu}\, \|f\|_{L^2}\,,
\end{align*}
where the sequence $\overline\bsb$ is defined in \eqref{eq:barbj-unif}.
\end{lemma}

\begin{lemma} \label{lem:regu2b}
In the uniform case under Assumptions~\ref{U1}--\ref{U3}, for every $f\in
L^2(D)$, every $\bsy\in U$, every $\bsnu\in\indx$, and every $h>0$, we
have
\begin{align*}
  \|\partial^\bsnu (u-u_h)(\cdot,\bsy)\|_V
  \,=\, \|\nabla (\partial^\bsnu (u-u_h)(\cdot,\bsy))\|_{L^2}
  \,\lesssim\, h\,
  (|\bsnu|+2)!\,\,\overline\bsb^\bsnu \,\|f\|_{L^2}\,,
\end{align*}
where the sequence $\overline\bsb$ is defined in \eqref{eq:barbj-unif}.
\end{lemma}

\begin{lemma} \label{lem:regu2c}
In the uniform case under Assumptions~\ref{U1}--\ref{U3}, for every $f\in
L^2(D)$, every $G\in L^2(D)$, every $\bsy\in U$, every $\bsnu\in\indx$,
and every $h>0$, we have
\begin{align*}
  |\partial^{\bsnu} G((u-u_h)(\cdot,\bsy))|
  \,\lesssim\, h^2\,
  (|\bsnu|+5)!\,\,\overline\bsb^\bsnu\, \|f\|_{L^2}\,\|G\|_{L^2}\,,
\end{align*}
where the sequence $\overline\bsb$ is defined in \eqref{eq:barbj-unif}.
\end{lemma}

The proofs of these three lemmas are given in Section~\ref{sec:proof}.
Lemma~\ref{lem:regu2a} is proved by induction on $|\bsnu|$, similarly to
the proof of Lemma~\ref{lem:regu1}, but this time by differentiating the
strong form \eqref{eq:strong} of the PDE and obtaining estimates involving
the Laplacian of the mixed derivatives of $u$; the proof makes use of
Lemma~\ref{lem:regu1}. Lemma~\ref{lem:regu2b} is also proved by induction
on $|\bsnu|$, but by differentiating the equation representing Galerkin
orthogonality; the proof makes use of Lemma~\ref{lem:regu2a}.
Lemma~\ref{lem:regu2c} is proved using a duality argument and it makes use
of Lemma~\ref{lem:regu2b}.

For $\ell\ge 1$ we obtain from Lemma~\ref{lem:regu2c} the norm estimate
\begin{align*}
 \| G(u^s_{h_\ell} - u^s_{h_{\ell-1}})\|_{s,\bsgamma}
 \,\lesssim\, h_{\ell-1}^2\, \|f\|_{L^2}\,\|G\|_{L^2}\,
 \Bigg(\sum_{\setu\subseteq\{1:s\}}
 \frac{[(|\setu|+5)!]^2\,\prod_{j\in\setu} \overline{b}_j^2}{\gamma_\setu}
 \Bigg)^{1/2}\,,
\end{align*}
which also holds for the case $\ell = 0$, see \eqref{eq:Gu-norm}, if we
define $h_{-1}:=1$. Combining these with Theorem~\ref{thm:QMC1} gives
\begin{align*}
  \sum_{\ell=0}^{L}
  \bbE\,[ |(I - Q\ran^\ell)(G(u^{s_\ell}_{h_\ell}-u^{s_{\ell-1}}_{h_{\ell-1}}))|^2 ]
  \,\lesssim\, C_\bsgamma^2(\lambda)\, \|f\|_{L^2}^2\,\|G\|_{L^2}^2
  \sum_{\ell=0}^L r_\ell^{-1}\, n_\ell^{-1/\lambda}\, h_{\ell-1}^4\,,
\end{align*}
where
\begin{align*}
  C_\bsgamma^2(\lambda)
  \,=\, \Bigg(2
  \sum_{\setu\subseteq\{1:s\}} \gamma_\setu^\lambda\,
  [\rho(\lambda)]^{|\setu|}
  \Bigg)^{1/\lambda}
  \,
  \Bigg(
  \sum_{\setu\subseteq\{1:s\}} \frac{[(|\setu|+5)!]^2 \prod_{j\in\setu} \overline{b}_j^2}{\gamma_\setu}
  \Bigg)\,.
\end{align*}
Here $C_\bsgamma^2(\lambda)$ is minimized by taking a different set of POD
weights
\begin{align} \label{eq:weight2}
  \gamma_\setu \,=\,
  \bigg(\frac{(|\setu|+5)!}{120}\, \prod_{j\in\setu} \frac{\overline{b}_j}{\sqrt{\rho(\lambda)}} \bigg)^{2/(1+\lambda)}\,.
\end{align}
Under Assumption~\ref{U6}, we can prove that if we take
\begin{equation} \label{eq:lambda-p1}
 \lambda \,=\,
 \begin{cases}
 \displaystyle\frac{1}{2-2\delta} \quad\mbox{for some } \delta\in (0,1/2)
 & \mbox{when } p_1\in (0,2/3]\,, \vspace{0.1cm} \\
 \displaystyle\frac{p_1}{2-p_1} & \mbox{when } p_1\in (2/3,1)\,,
 \end{cases}
\end{equation}
then the constant $C_\bsgamma^2(\lambda)$ can be bounded independently of
$s$. We point out again that the weights \eqref{eq:weight2} have a larger
factorial factor than those in \cite[Theorem~10]{KSS15} due to differences
in the definition of $\overline{b}_j$ and the proof argument.

We summarize the final result in the theorem below.

\begin{theorem} \label{thm:main2}
In the uniform case under Assumptions~\ref{U1}--\ref{U6}, for every $f\in
L^2(D)$ and every $G\in L^2(D)$, the multi-level algorithm
$A\ML\ran(G(u))$ using randomly shifted lattice rules with $s_\ell = s$
for all $\ell$, and with $n_\ell=2^{m_\ell}$ points constructed from a CBC
algorithm with POD weights \eqref{eq:weight2}--\eqref{eq:lambda-p1}, at
the pre-computation cost of $\calO(s\, n_\ell\log n_\ell + s^2n_\ell)$
operations, achieves the mean-square error bound
\begin{multline*}
  \bbE\left[ |(I-A\ML\ran)(G(u))|^2 \right]
  \,\lesssim\,
  s^{-4(1/p_0-1)} + h_L^{4} +
  \sum_{\ell=0}^L r_\ell^{-1} n_\ell^{-2\min(1/p_1-1/2,1-\delta)} h_{\ell-1}^{4}\,,
  \quad \delta\in (0,\tfrac{1}{2}) \,,
\end{multline*}
where the implied constant is independent of $s$, $h_\ell$, $r_\ell$ and
$n_\ell$.
\end{theorem}

In \cite{KSS15} it is assumed more generally that $f\in H^{-1+t}(D)$ and
$G\in H^{-1+t'}$ for $t,t'\in [0,1]$. Moreover, the analysis allowed for
different $s_\ell$ at different levels to arrive at the mean-square error
bound
\begin{align*}
  s_L^{-4(1/p_0-1)} \,+\, h_L^{2(t+t')}
  \,+\, \sum_{\ell=0}^L r_\ell^{-1} n_\ell^{-2\min(1/p_1-1/2,1-\delta)}
  \big( \theta_{\ell-1}\,s_{\ell-1}^{-2(1/p_0-1/p_1)} \,+\, h_{\ell-1}^{2(t+t')} \big)\,,
\end{align*}
for $\delta\in (0,1/2)$, where $\theta_{\ell-1}$ is $0$ if
$s_\ell=s_{\ell-1}$ and is $1$ otherwise. For this analysis a modified
sequence $\overline{\overline{b}}_j := \max (\overline{b}_j,
b_j^{p_0/p_1})$ is needed in place of $\overline{b}_j$ in the choice of
weights $\gamma_\setu$, see \cite[Theorem~11]{KSS15}.

Since the $s_\ell$ are potentially different for different $\ell$, the
strategy in \cite{KSS15} is to first choose $s_\ell$ in relation to
$h_\ell \asymp 2^{-\ell}$ to balance the truncation error and FE error on
each level. Observe from the error bound that the truncation error between
the levels decays significantly more slowly than the truncation error at
the highest level, i.e., comparing the exponent $2(1/p_0 - 1/p_1)$ with
$4(1/p_0-1)$. For this reason one ends up with a sequence of $s_\ell$ that
is strictly increasing initially but then the remaining $s_\ell$ become
constant and equal to $s_L$.

Then, with $r_\ell=r$ assumed to be a fixed constant, and with the
assumption that the cost of the algorithm $A\ML\ran(G(u))$ is
$\calO(\sum_{\ell=0}^L s_\ell\, n_\ell\, h_\ell^{-d})$ operations, a
Lagrange multiplier argument is used in \cite{KSS15} to choose $n_\ell$ in
relation to $h_\ell$ to minimize the mean-square error subject to a fixed
cost. Finally the value of $L$ is chosen so that the combined error meets
the required threshold. If the mean-square error is
$\calO(\varepsilon^2)$, then the cost of the multi-level algorithm can be
expressed as $\calO(\varepsilon^{-\tau}\,(\log \varepsilon^{-1})^\eta)$,
with $\tau$ much smaller than the corresponding exponent in the
single-level algorithm in most cases, see \cite[Theorem~12]{KSS15}.


\subsection{First order, single-level, lognormal} \label{sec:fo-sl-logn}

In the lognormal case, the mean-square error for our single-level
algorithm with randomly shifted lattice rules can be expressed as in
\eqref{eq:err-SL-ran}, but instead of using \eqref{eq:split-s}, here we
split the first term into
\[
  I(G(u-u^s_h)) \,=\, I(G(u-u_h)) \,+\, I(G(u_h -u^s_h))\,,
\]
where we estimate the FE error using an analogous result to
\eqref{eq:IGuh-unif} and we estimate the truncation error using
\eqref{eq:trunc-logn}. We use Theorem~\ref{thm:QMC2} for the QMC error,
namely, the second term in \eqref{eq:err-SL-ran}. As in the uniform case,
we need an estimate of $\|G(u^s_h)\|_{s,\bsgamma}$ but now with the norm
defined by \eqref{eq:norm2}. Using again \eqref{eq:lin-G}, we conclude
that we need a regularity result analogous to Lemma~\ref{lem:regu1}.

\begin{lemma} \label{lem:regu3}
In the lognormal case under Assumptions~\ref{L1} and~\ref{L2}, with the
sequence~$\bsbeta$ defined in \eqref{eq:bj-logn}, for every $f\in V^*$,
every $\bsy\in U_\bsbeta$, and every $\bsnu\in\indx$, we have
\begin{align*}
  \|\partial^{\bsnu}u(\cdot,\bsy)\|_{V}
  \,=\, \|\nabla (\partial^{\bsnu}u (\cdot,\bsy))\|_{L^2}
 \,\le\, \frac{|\bsnu|!}{(\ln 2)^{|\bsnu|}} \, \bsbeta^\bsnu\, \frac{\|f\|_{V^*}}{a_{\min}(\bsy)}\,.
\end{align*}
The same estimate holds when the exact solution $u$ is replaced by
$u^s_h$.
\end{lemma}

The proof is given in Section~\ref{sec:proof}; see also
\cite[Theorem~14]{GKNSSS15}. It is proved by obtaining an estimate on
$\|a^{1/2}(\cdot,\bsy)\,\nabla(\partial^{\bsnu}u(\cdot,\bsy))\|_{L^2}$
using induction on $|\bsnu|$.

Comparing with the corresponding result for the uniform case, the critical
difference in the lognormal case is that $a_{\min}(\bsy)$ depends on
$\bsy$, thus making the estimation of $\|G(u^s_h)\|_{s,\bsgamma}$ more
complicated. In the following, we argue along the lines of the proofs of
\cite[Theorems~16 and~20 and Corollary~21]{GKNSSS15}, however, the
estimates are slightly different since there was an omission of a factor
in \cite{GKNSSS15} which meant that the formula for the weights
$\gamma_\setu$ is different and the implied constant in
Theorem~\ref{thm:main3} is smaller here compared with \cite{GKNSSS15}.

We can use \eqref{eq:lin-G} and Lemma~\ref{lem:regu3}, together with the
bound
\[
  \frac{1}{a_{\min}(\bsy)} \,\le\,
  \frac{1}{\min_{\bsx\in \overline{D}} a_0(\bsx)}\,\prod_{j\ge 1} \exp(\beta_j\,|y_j|)\,,
\]
as well as
\begin{equation} \label{eq:beta-int1}
  \int_{-\infty}^\infty \exp(\beta_j\,|y_j|)\,\phi(y_j)\,\rd y_j
  \,=\, 2\exp(\beta_j^2/2)\Phi(\beta_j)
  \,\le\, \exp(\beta_j^2/2+2\beta_j/\sqrt{2\pi})
  \,,
\end{equation}
to obtain an estimate of the norm \eqref{eq:norm2}
\begin{multline} \label{eq:hindsight}
  \|G(u^s_h)\|_{s,\bsgamma}^2
  \,\le\,
  \bigg(\frac{\|f\|_{V^*}\,\|G\|_{V^*}}{\min_{\bsx\in \overline{D}} a_0
  (\bsx)}\bigg)^2
  \prod_{j=1}^s [2\exp(\beta_j^2/2)\Phi(\beta_j)]^2 \\
  \times
  \sum_{\setu\subseteq\{1:s\}} \frac{(|\setu|!)^2}{\gamma_\setu (\ln 2)^{2|\setu|}}
  \prod_{j\in\setu} \frac{\beta_j^2}{[2\exp(\beta_j^2/2)\Phi(\beta_j)]^2}\,
  \int_{-\infty}^\infty \exp(2\beta_j |y_j|)\,\vpsi_j^2(y_j)\,\rd y_j\,.
\end{multline}
This leads us to choose the weight functions $\vpsi_j$ to be exponential
functions given by \eqref{eq:psi}, with $\alpha_j> \beta_j$, so that
\begin{equation} \label{eq:beta-int2}
  \int_{-\infty}^\infty \exp(2\beta_j\,|y_j|)\,\vpsi_j^2(y_j)\,\rd y_j
  \,=\, \frac{1}{\alpha_j-\beta_j}\,,
\end{equation}
and thus
\begin{equation} \label{eq:messy-norm}
  \|G(u^s_h)\|_{s,\bsgamma}^2
  \,\lesssim\, \sum_{\setu\subseteq\{1:s\}} \frac{(|\setu|!)^2}{\gamma_\setu (\ln 2)^{2|\setu|}}
  \prod_{j\in\setu} \frac{\beta_j^2}{[2\exp(\beta_j^2/2)\Phi(\beta_j)]^2\, (\alpha_j-\beta_j)}\,,
\end{equation}
where the implied constant is independent of $s$ under
Assumption~\ref{L1}.

Combining this with Theorem~\ref{thm:QMC2} and following the strategy for
choosing weights $\gamma_\setu$ in the uniform case, we obtain
\begin{equation} \label{eq:weight3}
  \gamma_\setu \,=\,
  \Bigg( |\setu|!
  \prod_{j\in\setu} \frac{\beta_j}{2(\ln 2)\,\exp(\beta_j^2/2)\Phi(\beta_j)\, \sqrt{(\alpha_j-\beta_j)\,\rho_j(\lambda)}}
  \Bigg)^{2/(1+\lambda)}\,,
\end{equation}
with $\lambda$ chosen as in \eqref{eq:lambda-p0} but with $p_0$ as in
Assumption~\ref{L5}, where $\rho_j(\lambda)$ is as given in
Theorem~\ref{thm:QMC2}. We substitute this choice of weights into
Theorem~\ref{thm:QMC2} and \eqref{eq:messy-norm}, and then minimize the
resulting bound with respect to the parameters $\alpha_j$. As argued in
\cite[Corollary~21]{GKNSSS15}, this corresponds to minimizing
$[\varrho_j(\lambda)]^{1/\lambda}/(\alpha_j-\beta_j)$ with respect to
$\alpha_j$, and yields
\begin{equation} \label{eq:alpha3}
  \alpha_j \,=\, \frac{1}{2} \bigg(
  \beta_j + \sqrt{\beta_j^2+1-\frac{1}{2\lambda}}\bigg)\,.
\end{equation}

Note that $2\exp(\beta_j^2/2)\Phi(\beta_j)\to 1$ as $\beta_j\to 0$, and
$\rho_j(\lambda)$ is also bounded away from zero as $\beta_j\to 0$.
Moreover, $\alpha_j-\beta_j>0$ is minimized by the largest $\beta_j$. Thus
the summability of the product factors in \eqref{eq:weight3} is determined
by the summability of the numerator $\beta_j$. Arguing as in the proof of
\cite[Theorem~20]{GKNSSS15}, we conclude that
\[
  \bbE\,[ |(I-Q\ran)(G(u^s_h))|^2 ] \,\lesssim\, r^{-1} \, n^{-1/\lambda}\,,
\]
where the implied constant is independent of $s$ under
Assumption~\ref{L5}.

We summarize the final result in the theorem below.

\begin{theorem} \label{thm:main3}
In the lognormal case under Assumptions~\ref{L1}--\ref{L5}, for every
$f\in L^2(D)$ and every $G\in L^2(D)$, the single-level algorithm
$A\SL\ran(G(u))$ using a randomly shifted lattice rule with $n=2^m$ points
constructed from a CBC algorithm with POD weights
\eqref{eq:weight3}--\eqref{eq:alpha3} together with \eqref{eq:lambda-p0},
at the pre-computation cost of $\calO(s\, n\log n + s^2n)$ operations,
achieves the mean-square error bound
\begin{align*}
  \bbE\left[ |(I-A\SL\ran)(G(u))|^2 \right]
  \,\lesssim\,
  s^{-2\chi} + h^{4} + r^{-1}\,n^{-2\min(1/p_0-1/2,1-\delta)}\,, \quad
  \delta\in (0,\tfrac{1}{2})\,,
\end{align*}
with $\chi$ as in \eqref{eq:trunc-logn}, where the implied constant is
independent of $s$, $h$, $r$ and $n$.
\end{theorem}

Similarly to the uniform case, we can treat $r$ as a fixed constant and
choose $s$, $h$, and $n$ to achieve $\calO(\varepsilon^2)$ mean-square
error.

We remark that an alternative analysis can be carried out following
Theorem~\ref{thm:QMC2b} instead of Theorem~\ref{thm:QMC2}, by choosing
different weight functions $\vpsi_j$ in \eqref{eq:hindsight}. This would
yield a different formula \eqref{eq:beta-int2}, a different bound for the
norm \eqref{eq:messy-norm}, and a different choice of weights
\eqref{eq:weight3}. This is work in progress and could potentially lead to
better constants in the bounds and therefore better numerical results.


\subsection{First order, multi-level, lognormal} \label{sec:fo-ml-logn}

For our multi-level algorithm in the lognormal case, the mean-square error
with randomly shifted lattice rules can again be expressed as
\eqref{eq:err-ML-ran}. We now need to estimate $\|G(u^s_{h_\ell} -
u^s_{h_{\ell-1}})\|_{s,\bsgamma}$ with the norm defined by
\eqref{eq:norm2}. We need regularity results similar to
Lemmas~\ref{lem:regu2a}, \ref{lem:regu2b}, and \ref{lem:regu2c} in the
uniform case.

\begin{lemma} \label{lem:regu4a}
In the lognormal case under Assumptions~\ref{L1}--\ref{L3}, with the
sequence~$\overline\bsbeta$ defined in \eqref{eq:barbj-logn}, for every
$f\in L^2(D)$, every $\bsy\in U_{\overline\bsbeta}$, and every
$\bsnu\in\indx$, we have
\begin{align*}
  \|\Delta (\partial^{\bsnu}u (\cdot,\bsy))\|_{L^2}
  \,\lesssim\, T(\bsy)\, (|\bsnu|+1)!\, 2^{|\bsnu|}\,\overline\bsbeta^\bsnu\, \|f\|_{L^2}\,,
\end{align*}
where $T(\bsy)$ is defined in \eqref{eq:defT}.
\end{lemma}

\begin{lemma} \label{lem:regu4b}
In the lognormal case under Assumptions~\ref{L1}--\ref{L3}, with the
sequence~$\overline\bsbeta$ defined in \eqref{eq:barbj-logn}, for every
$f\in L^2(D)$, every $\bsy\in U_{\overline\bsbeta}$, every
$\bsnu\in\indx$, and every $h>0$, we have
\begin{align*}
  \|a^{1/2}(\cdot,\bsy)\nabla (\partial^{\bsnu}(u-u_h)(\cdot,\bsy))\|_{L^2}
  \,\lesssim\, h\, T(\bsy)\,a_{\max}^{1/2}(\bsy)\,
  \frac{(|\bsnu|+2)!}{2}\, 2^{|\bsnu|}\,\overline\bsbeta^\bsnu\, \|f\|_{L^2}\,,
\end{align*}
where $T(\bsy)$ is defined in \eqref{eq:defT}.
\end{lemma}

\begin{lemma} \label{lem:regu4c}
In the lognormal case under Assumptions~\ref{L1}--\ref{L3}, with the
sequence~$\overline\bsbeta$ defined in \eqref{eq:barbj-logn}, for every
$f\in L^2(D)$, every $G\in L^2(D)$, every $\bsy\in U_{\overline\bsbeta}$,
every $\bsnu\in\indx$, and every $h>0$, we have
\begin{align*}
  |\partial^{\bsnu} G((u - u_h)(\cdot,\bsy))|
  \,\lesssim\, h^2\, T^2(\bsy)\,a_{\max}(\bsy)\,
  (|\bsnu|+5)!\, 2^{|\bsnu|}\,\overline\bsbeta^\bsnu\, \|f\|_{L^2}\,\|G\|_{L^2}\,,
\end{align*}
where $T(\bsy)$ is defined in \eqref{eq:defT}.
\end{lemma}

The proof of these three lemmas are given in Section~\ref{sec:proof}; see
also \cite{KSSSU}.
These proofs follow the same general steps as in the proofs for
Lemmas~\ref{lem:regu2a}, \ref{lem:regu2b}, and \ref{lem:regu2c} in the
uniform case. The main challenge is in identifying which power of
$a(\cdot,\bsy)$ to be included in the $L^2$ norm estimate so that the
recursion works. Lemma~\ref{lem:regu4a} is proved by obtaining an estimate
on $\|a^{-1/2}(\cdot,\bsy)\nabla
\cdot(a(\cdot,\bsy)\,\nabla(\partial^{\bsnu}u(\cdot,\bsy)))\|_{L^2}$ using
induction on $|\bsnu|$; the proof makes use of an estimate on
$\|a^{1/2}(\cdot,\bsy)\,\nabla(\partial^{\bsnu}u(\cdot,\bsy))\|_{L^2}$
which was established in the proof of Lemma~\ref{lem:regu3}.
Lemma~\ref{lem:regu4b} is also obtained by induction on $|\bsnu|$; the
proof makes use of the estimate established in proof of
Lemma~\ref{lem:regu4a}. The proof of Lemma~\ref{lem:regu4c} makes use of
Lemma~\ref{lem:regu4b}.

As in \cite{KSSSU} we can show that
\begin{align*}
  T^2(\bsy) \, a_{\max}(\bsy)
  &\le
  \left(1 + \frac{\|\nabla a_0\|_{L^\infty}}{\min_{\bsx\in \overline{D}}a_0(\bsx)}\right)^2
  \frac{\|a_0\|_{L^\infty}^3}{(\min_{\bsx\in \overline{D}}a_0(\bsx))^4}
  \prod_{j\ge1} \exp(9 \overline{\beta}_j \, |y_j|)\,.
\end{align*}
Thus, with the weight functions $\vpsi_j$ in the norm \eqref{eq:norm2}
given by \eqref{eq:psi} where $\alpha_j > 9 \overline\beta_j$, and using
\eqref{eq:beta-int1} and \eqref{eq:beta-int2} with $\beta_j$ replaced by
$9\overline\beta_j$, we obtain from Lemma~\ref{lem:regu4c} the estimate
for $\ell\ge 1$
\begin{align*}
 \| G(u^s_{h_\ell} - u^s_{h_{\ell-1}})\|_{s,\bsgamma}^2 
 &\,\lesssim\, h_{\ell-1}^4\,
 \sum_{\setu\subseteq\{1:s\}} \frac{[(|\setu|+5)!]^2\, 2^{2|\setu|}}{\gamma_\setu}
 \prod_{j\in\setu} \frac{\overline{\beta}_j^2}
 {[2\exp(81\overline{\beta}_j^2/2)\Phi(9\overline{\beta}_j)]^2\, (\alpha_j-9\overline{\beta}_j)}\,,
\end{align*}
where the implied constant is independent of $s$ under
Assumption~\ref{L1}. This also holds for $\ell=0$ with $h_{-1}:=1$, see
\eqref{eq:messy-norm}. Combining this with Theorem~\ref{thm:QMC2} and
following the same strategy for choosing weights as in the previous
subsections, we conclude that
\begin{align*}
  \sum_{\ell=0}^{L}
  \bbE\,[ |(I - Q\ran^\ell)(G(u^s_{h_\ell}-u^s_{h_{\ell-1}}))|^2 ]
  \,\lesssim\, \sum_{\ell=0}^L r_\ell^{-1}\, n_\ell^{-1/\lambda}\, h_{\ell-1}^4\,,
\end{align*}
where the implied constant is independent of $s$ under
Assumption~\ref{L6}, if we take $\lambda$ as in \eqref{eq:lambda-p1} and
if we choose the weights
\begin{equation} \label{eq:weight4}
  \gamma_\setu \,=\,
  \bigg(\frac{(|\setu|+5)!}{120}\, \prod_{j\in\setu} \frac{\overline{\beta}_j}
  {\exp(81\overline{\beta}_j^2/2)\Phi(9\overline{\beta}_j)\, \sqrt{(\alpha_j-9\overline{\beta}_j)\,\rho_j(\lambda)}}
  \bigg)^{2/(1+\lambda)}\,.
\end{equation}
A generalization of \cite[Corollary~21]{GKNSSS15} yields the choice
\begin{equation} \label{eq:alpha4}
  \alpha_j \,=\, \frac{1}{2} \bigg(
  9\overline\beta_j + \sqrt{81\overline\beta_j^2+1-\frac{1}{2\lambda}}\bigg)\,.
\end{equation}

We summarize the final result in the theorem below.

\begin{theorem} \label{thm:main4}
In the lognormal case under Assumptions~\ref{L1}--\ref{L6}, for every
$f\in L^2(D)$ and every $G\in L^2(D)$, the multi-level algorithm
$A\ML\ran(G(u))$ using randomly shifted lattice rules with $s_\ell = s$
for all $\ell$, and with $n_\ell=2^{m_\ell}$ points constructed from a CBC
algorithm with POD weights \eqref{eq:weight4}--\eqref{eq:alpha4} together
with \eqref{eq:lambda-p1}, at the pre-computation cost of $\calO(s\,
n_\ell\log n_\ell + s^2n_\ell)$ operations, achieves the mean-square error
bound
\begin{align*}
  \bbE\left[ |(I-A\ML\ran)(G(u))|^2 \right]
  \,\lesssim\,
  s^{-\chi} + h_L^{4} +
  \sum_{\ell=0}^L r_\ell^{-1} n_\ell^{-2\min(1/p_1-1/2,1-\delta)} h_{\ell-1}^{4},
   \; \delta\in (0,\tfrac{1}{2}) \,,
\end{align*}
where the implied constant is independent of $s$, $h_\ell$, $r_\ell$, and
$n_\ell$.
\end{theorem}

We can carry out a cost versus error analysis as in the uniform case to
obtain the optimal choice of values $s$, $h_\ell$, and $n_\ell$. We could
also adjust the levels adaptively in practice, as described in
\cite{KSSSU}, see also \cite{GW09}.

As in the single-level algorithm, we could choose the weight functions
$\vpsi_j$ differently by following Theorem~\ref{thm:QMC2b} and this would
yield a different choice of weights \eqref{eq:weight4}.


\subsection{Higher order, single-level, uniform} \label{sec:ho-sl-unif}

The error for our single-level algorithm with interlaced polynomial
lattice rules can be expressed as
\[
  (I-A\SL\determ)(G(u))
  \,=\, I(G(u-u^s)) \,+\, I(G(u^s-u^s_h)) \,+\, (I - Q)(G(u^s_h))\,.
\]
The truncation error can be estimated using \eqref{eq:trunc-unif} as
before. The FE error of higher order can be estimated using
\eqref{eq:FE-ho}. For the QMC error we use Theorem~\ref{thm:QMC3} and
therefore we need an estimate on the norm
$\|G(u^s_h)\|_{s,\alpha,\bsgamma}$ defined by \eqref{eq:norm3}. Due to
linearity and boundedness of $G$ we have $|\partial^\bsnu
G(u^s_h(\cdot,\bsy))| \le \|G\|_{V^*}\,\|\partial^\bsnu
u^s_h(\cdot,\bsy)\|_V$, where the last norm can be estimated as in
Lemma~\ref{lem:regu1}. Thus we have
\begin{align*}
  \|G(u^s_h)\|_{s,\alpha,\bsgamma}
 &\,\le\, \frac{\|f\|_{V^*}\,\|G\|_{V^*}}{a_{\min}}
 \sup_{\setu\subseteq\{1:s\}}
 \frac{1}{\gamma_\setu}
 \sum_{\setv\subseteq\setu} \sum_{\bstau_{\setu\setminus\setv} \in \{1:\alpha\}^{|\setu\setminus\setv|}}\!
 |(\bsalpha_\setv,\bstau_{\setu\setminus\setv},\bszero)|!\,
 \bsb^{(\bsalpha_\setv,\bstau_{\setu\setminus\setv},\bszero)}
 \\
 &\,=\, \frac{\|f\|_{V^*}\,\|G\|_{V^*}}{a_{\min}}
 \sup_{\setu\subseteq\{1:s\}}
 \frac{1}{\gamma_\setu}
 \sum_{\bsnu_\setu \in \{1:\alpha\}^{|\setu|}}
 |\bsnu_\setu|!\,\prod_{j\in\setu} \big(2^{\delta(\nu_j,\alpha)} b_j^{\nu_j}\big)\,,
\end{align*}
where $\delta(\nu_j,\alpha)$ is $1$ if $\nu_j=\alpha$ and is $0$
otherwise.
We choose SPOD weights \eqref{eq:SPOD}
\begin{equation} \label{eq:weight5}
 \gamma_\setu
 \,=\,
 \sum_{\bsnu_\setu \in \{1:\alpha\}^{|\setu|}}
 |\bsnu_\setu|!\,\prod_{j\in\setu} \big(2^{\delta(\nu_j,\alpha)} b_j^{\nu_j}\big)\,,
\end{equation}
so that $\|G(u^s_h)\|_{s,\alpha,\bsgamma} \le
\|f\|_{V^*}\,\|G\|_{V^*}/a_{\min}$. Substituting these weights into
Theorem~\ref{thm:QMC3} and following the arguments in
\cite[pp.~2694--2695]{DKLNS14}, we take $\lambda = p_0$ and $\alpha =
\lfloor 1/p_0\rfloor + 1$, and conclude that
\[
  |(I - Q)(G(u^s_h))| \,\lesssim\, n^{-1/p_0}\,,
\]
with the implied constant independent of $s$.

We summarize the final result in the theorem below.

\begin{theorem} \label{thm:main5}
In the uniform case under Assumptions~\ref{U1}--\ref{U5}, for every $f\in
\calX_t^*$ and every $G\in \calX_{t'}^*$ with $t,t'\ge 0$, the
single-level algorithm $A\SL\determ(G(u))$ using an interlaced polynomial
lattice rule with interlacing factor $\alpha = \lfloor 1/p_0\rfloor+1\ge
2$ and with $n=2^m$ points constructed from a CBC algorithm with SPOD
weights \eqref{eq:weight5}, at the pre-computation cost of
$\calO(\alpha\,s\, n\log n + \alpha^2 s^2n)$ operations, achieves the
error bound
\begin{align*}
  |(I-A\SL\determ)(G(u))|
  \,\lesssim\,
  s^{-2(1/p_0-1)} + h^{t+t'} + n^{-1/p_0}\,,
\end{align*}
where the implied constant is independent of $s$, $h$, and $n$.
\end{theorem}

Again we can choose $s$, $h$, $n$ to achieve $\calO(\varepsilon)$ error.


\subsection{Higher order, multi-level, uniform} \label{sec:ho-ml-unif}

The error for our multi-level algorithm with interlaced polynomial lattice
rules can be expressed as
\begin{align*}
  (I - A\ML\determ)(G(u)) 
  &\,=\, I(G(u-u^{s_L})) \,+\, I(G(u^{s_L}-u^{s_L}_{h_L}))
  \,+\, \sum_{\ell=0}^L (I - Q^\ell)(G(u^{s_{\ell}}_{h_\ell} - u^{s_{\ell-1}}_{h_{\ell-1}}))\,,
\end{align*}
where $Q^\ell$ denotes an interlaced polynomial lattice rule in $s_\ell$
dimensions with $n_\ell$ points. Again for simplicity we restrict our
discussion to the case $s_\ell = s$ for all $\ell$. So we need an estimate
on the norm $\|G(u^{s}_{h_\ell} -
u^{s}_{h_{\ell-1}})\|_{s,\alpha,\bsgamma}$. From Lemma~\ref{lem:regu2c} we
conclude that for $\ell\ge 1$ we have
\begin{align*}
 &\|G(u^{s}_{h_\ell} - u^{s}_{h_{\ell-1}})\|_{s,\alpha,\bsgamma} \\
 &\,\lesssim\, h_{\ell-1}^2\, \|f\|_{L^2}\,\|G\|_{L^2}\,
 \sup_{\setu\subseteq\{1:s\}}
 \frac{1}{\gamma_\setu}
 \sum_{\setv\subseteq\setu} \sum_{\bstau_{\setu\setminus\setv} \in \{1:\alpha\}^{|\setu\setminus\setv|}}\!
 (|(\bsalpha_\setv,\bstau_{\setu\setminus\setv},\bszero)|+5)!\,
 \overline\bsb^{(\bsalpha_\setv,\bstau_{\setu\setminus\setv},\bszero)}
 \\
 &\,=\, h_{\ell-1}^2\, \|f\|_{L^2}\,\|G\|_{L^2}\,
 \sup_{\setu\subseteq\{1:s\}}
 \frac{1}{\gamma_\setu}
 \sum_{\bsnu_\setu \in \{1:\alpha\}^{|\setu|}}
 (|\bsnu_\setu|+5)!\,\prod_{j\in\setu} \big(2^{\delta(\nu_j,\alpha)} \overline{b}_j^{\nu_j}\big)\,.
\end{align*}
We therefore choose SPOD weights
\begin{equation} \label{eq:weight6}
 \gamma_\setu
 \,=\,
 \sum_{\bsnu_\setu \in \{1:\alpha\}^{|\setu|}}
 \frac{(|\bsnu_\setu|+5)!}{120}\,\prod_{j\in\setu} \big(2^{\delta(\nu_j,\alpha)} \overline{b}_j^{\nu_j}\big)\,.
\end{equation}
Substituting these weights into Theorem~\ref{thm:QMC3} and following
similar arguments as those before, we take now $\lambda = p_1$ and $\alpha
= \lfloor 1/p_1\rfloor + 1$, and conclude that the convergence rate is
$n^{-1/p_1}$, with the implied constant independent of $s$.

\begin{theorem} \label{thm:main6}
In the uniform case under Assumptions~\ref{U1}--\ref{U6}, for every $f\in
L^2(D)$ and every $G\in L^2(D)$, the multi-level algorithm
$A\ML\determ(G(u))$ using interlaced polynomial lattice rules with
interlacing factor $\alpha = \lfloor 1/p_1 \rfloor + 1\ge 2$, with $s_\ell
= s$ for all $\ell$, and with $n_\ell=2^{m_\ell}$ points constructed from
a CBC algorithm with SPOD weights \eqref{eq:weight6}, at the
pre-computation cost of $\calO(\alpha\,s\, n_\ell\log n_\ell +
\alpha^2\,s^2n_\ell)$ operations, achieves the error bound
\begin{align*}
  |(I-A\ML\determ)(G(u))|
  \,\lesssim\,
  s^{-2(1/p_0-1)} + h_L^{2} +
  \sum_{\ell=0}^L n_\ell^{-1/p_1} h_{\ell-1}^{2}\,,
\end{align*}
where the implied constant is independent of $s$, $h_\ell$, and $n_\ell$.
\end{theorem}

We could consider higher order FE methods here, but then we would need
stronger regularity on $f$ and $G$, as well as analogous results for
Lemmas~\ref{lem:regu2a}, \ref{lem:regu2b}, and \ref{lem:regu2c}. This is
where Assumption~\ref{U7} would be needed. We could also allow different
$s_\ell$ at different levels. These generalizations are considered in
\cite{DKLS} where a comprehensive error versus cost analysis of the
multi-level algorithm is provided for the general setting of affine
parametric operator equations.

Now we summarize and compare the results from
\cite{KSS12,KSS15,DKLNS14,DKLS} for the uniform case:
\begin{align*}
 & \mbox{First-order single-level \cite{KSS12}} \\
 & \quad s^{-2(1/p_0-1)} + h^{t+t'} + n^{-\min(1/p_0-1/2,1-\delta)} \quad\mbox{(rms)}. \\
 & \mbox{First-order multi-level \cite{KSS15}} \\
 & \quad s_L^{-2(1/p_0-1)} + h_L^{t+t'}
 + \sum_{\ell=0}^L n_\ell^{-\min(1/p_1-1/2,1-\delta)}
                   \big(\theta_{\ell-1}\,s_{\ell-1}^{-(1/p_0-1/p_1)}
                   + h_{\ell-1}^{t+t'}\big) \quad\mbox{(rms)}. \\
 & \mbox{Higher-order single-level \cite{DKLNS14}} \\
 & \quad s^{-2(1/p_0-1)} + h^{t+t'} + n^{-1/p_0}. \\
 & \mbox{Higher-order multi-level \cite{DKLS}} \\
 & \quad s_L^{-2(1/p_0-1)} + h_L^{t+t'}
 + \sum_{\ell=0}^L n_\ell^{-1/p_t}
 \big(\theta_{\ell-1}\,s_{\ell-1}^{-(1/p_0-1/p_t)} +  h_{\ell-1}^{t+t'} \big).
\end{align*}
For the first-order results, ``rms'' indicates that the error is in the
root-mean-square sense since we use a randomized QMC method. The
higher-order results are deterministic. The results include general
parameters $t,t'$ for the regularity of $f$ and $G$: in the first order
results we have $f\in H^{-1+t}(D)$ and $G\in H^{-1+t'}(D)$ for $t,t'\in
[0,1]$, while in the higher order results we have $f\in \calX_t^*$ and
$G\in \calX_{t'}^*$ for integers $t,t'\ge 0$. For the multi-level results
we include the analysis for potentially taking different $s_\ell$ at each
level. Recall that $\delta\in (0,1/2)$, and $\theta_{\ell-1}$ is $0$ if
$s_{\ell} = s_{\ell-1}$ and is $1$ otherwise.

Note that in many applications $p_t$ in Assumption~\ref{U7} satisfies
\[
  p_t \,=\, \frac{p_0}{1-tp_0/d}\,,
\]
which means it can be much bigger than $p_0$. So the higher order
multi-level algorithm does not necessarily lead to improved error bounds.

\section{A practical guide to the software for constructing QMC points} \label{sec:comp}

In this section we explain how to use the code which accompanies this
article to construct QMC rules for the different settings which have been
discussed in this article. The construction algorithms are all fast
component-by-component constructions, using results from
\cite{NC06a,NC06b,CKN06,DKLNS14} and \cite{Nuy14}. For lattice rules the
construction algorithm will output the generating vector and for the
interlaced polynomial lattice rules the algorithm will output the
generating matrices. These generating vectors and matrices can then be
used in the provided sample point generators.

For randomly shifted lattice rules we have to construct a good generating
vector~$\bsz$ for the lattice rule~\eqref{eq:r1lr}. The results in
Sections~\ref{sec:qmc} and~\ref{sec:error} were stated for $n=2^m$ points,
but as noted before they hold for any prime power (including $n$ being a
prime). The construction script expects the prime $p$ (which defaults
to~$2$) and the power $m$ to be given such that the (maximum) number of
points is $n=p^m$. The natural thing to do is to construct such rules
to be good for any intermediate power of $p$ as an embedded sequence of
lattice rules. Such \emph{lattice sequences} can then be used as a
sequence of QMC rule approximations. A construction of such lattice
sequences was proposed in \cite{CKN06} and this is the approach followed
in the lattice rule construction code.

For interlaced polynomial lattice rules~\eqref{eq:ipr1lr} we have to
construct the associated generating matrices. We fix the base of the
finite field to be~$2$ for practicalities in the construction, and in the
generation of the points. The number of points is $n = 2^m$, but, in
contrast to the lattice rules, we currently do not provide these as
embedded sequences. The construction script will automatically choose a
default irreducible polynomial of degree~$m$ as the modulus polynomial
(which can be overridden by the user). The specific choice of modulus
polynomial does not influence the error bound in Theorem~\ref{thm:QMC3}.
The output of the script will be both the generating matrices
$C_1,\ldots,C_{\alpha \, s} \in \bbZ_2^{m \times m}$ of the polynomial
lattice rule~\eqref{eq:pr1lr-matrix} as well as the generating matrices
$B_1,\ldots,B_s \in \bbZ_2^{\alpha m \times m}$ of the interlaced
polynomial lattice rule~\eqref{eq:ipr1lr}.

Once the generating vector or the generating matrices have been
constructed, they are used as input to the corresponding point generator.
These point generators are relatively straightforward to program, and
their computational cost to generate a point is really minor and so can be
neglected for practical purposes; in fact they are comparable to the
fastest random number generators.
To generate lattice points, one only requires an integer multiplication, a
modulus operation, and a fixed float multiplication/division per
dimension, see~\eqref{eq:r1lr}. This is comparable to the cost of a simple
LCG (linear congruential generator) per dimension.
To generate (interlaced) polynomial lattice points in Gray code ordering,
one needs an \texttt{xor} instruction and a fixed float
multiplication/division per dimension, see~\eqref{eq:pr1lr-matrix}
and~\eqref{eq:ipr1lr}. Additionally, a CTZ (count trailing zeros)
algorithm is used to determine the column number of the generating
matrices to perform each \texttt{xor} instruction, which is available on
most CPUs as a machine instruction or can be implemented with a simple
algorithm having a fixed low arithmetic cost. This is comparable to the
cost of a LFSR (linear feedback shift register) generator per dimension.
In the case of randomly shifted lattice rules the points still have to be
randomly shifted before being used as quadrature points. As QMC points are
naturally enumerated, it is straightforward to parallelize the solving of
the different PDE problems and we therefore equip the point generators
with an option to start at any offset in the enumeration of the points.

As the theory in Section~\ref{sec:error} is often quite involved, we
extract the essential properties of the analysis here and allow them to be
applied to any general integrand which shares similar characteristics. In
essence, the analysis in Section~\ref{sec:error} made use of bounds on the
mixed derivatives, $\partial^{\bsnu}$ for $\bsnu \in \{0, ...,
\alpha\}^s$, of the integrands to derive suitable weights $\gamma_\setu$.
In the uniform case these bounds can be stated as follows.
\begin{itemize}
\item Uniform, single-level,
      with $F(\bsy) \,=\, G(u_h^s(\cdot,\bsy))$:
\begin{align*}
  \left| \partial^{\bsnu} F(\bsy) \right|
  &\,\lesssim\,
  |\bsnu|! \, \prod_{j=1}^s b_j^{\nu_j}
  \,.
\end{align*}
\item Uniform, multi-level, with $F_\ell(\bsy) \,=\,
    G((u_{h_\ell}^s-u_{h_{\ell-1}}^s)(\cdot,\bsy))$:
\begin{align*}
   \left| \partial^{\bsnu} F_\ell(\bsy) \right|
  &\,\lesssim\,
  (|\bsnu|+5)! \, \prod_{j=1}^s \overline{b}_j^{\nu_j}
  \,.
\end{align*}
\end{itemize}
For first order methods the multi-index $\bsnu$ satisfies $\nu_j \le 1$
for all~$j$. In the lognormal case the bounds are not uniformly bounded in
$\bsy$ and can be stated as follows.
\begin{itemize}
\item Lognormal, single-level,
      with $F(\bsy) \,=\, G(u_h^s(\cdot,\bsy))$:
\begin{align*}
   \left| \partial^{\bsnu}F (\bsy) \right|
   &\,\lesssim\,
  |\bsnu|! \,
  \prod_{j=1}^s
  (\beta_j / \ln 2)^{\nu_j} \,
  \exp(\beta_j \, |y_j|)
  \,.
\end{align*}
\item Lognormal, multi-level, with $F_\ell(\bsy) \,=\,
    G((u_{h_\ell}^s-u_{h_{\ell-1}}^s)(\cdot,\bsy))$:
\begin{align*}
   \left| \partial^{\bsnu}F_\ell (\bsy) \right|
   &\,\lesssim\,
  (|\bsnu|+5)! \,
  \prod_{j=1}^s
  (2\overline{\beta}_j)^{\nu_j} \,
  \exp(9 \overline{\beta}_j \, |y_j|)
  \,.
\end{align*}
\end{itemize}
We note that the results in Section~\ref{sec:error} for the lognormal case
only hold for first order methods with all $\nu_j \le 1$.

We provide two Python scripts, \texttt{lat-cbc.py} and
\texttt{polylat-cbc.py} (as interfaces to the construction script
\texttt{spod-cbc.py}), to construct lattice rules and interlaced
polynomial lattice rules, respectively, in which the following generalized
bound on the mixed derivatives is assumed: for all $\bsnu \in \{0, ...,
\alpha\}^s$,
\begin{align}\label{eq:general-bound}
   \left| \partial^\bsnu F (\bsy) \right|
   &\,\lesssim\,
     \Bigl( (|\bsnu|+a_1)! \Bigr)^{d_1}
     \prod_{j=1}^s (a_2 \, \Bj)^{\nu_j}\,
                    \exp(a_3\,\Bj|y_j|)
  \,,
\end{align}
for some integers $\alpha \ge 1$ and $a_1 \ge 0$, real numbers $a_2 > 0$,
$a_3 \ge 0$ and $d_1 \ge 0$, and a sequence of positive numbers $\Bj$,
corresponding to the values of $b_j$, $\overline{b}_j$, $\beta_j$ or
$\overline{\beta}_j$, appropriate for the setting, see~\eqref{eq:bj-unif},
\eqref{eq:barbj-unif}, \eqref{eq:bj-logn} and \eqref{eq:barbj-logn}. An
overview of all parameters and their description for the
\texttt{lat-cbc.py} and the \texttt{polylat-cbc.py} scripts is given in
Table~\ref{tab:options}.

\begin{table}[t]
\caption{Options for the Python scripts \texttt{lat-cbc.py} and
\texttt{polylat-cbc.py}.}\label{tab:options}
{\hspace*{-1.3mm}\small\begin{tabular}{lll} \hline
\texttt{s} & & number of dimensions \\[1mm]
\texttt{m} & & number of points given by $2^m$ \\
           & & (or $p^m$ in case of optional argument \texttt{p} for lattice rules) \\[1mm]
\texttt{p} & optional & \mybullet for lattice rules $n=p^m$, prime, defaults to $p=2$ \\
           &          & \mybullet for polynomial lattice rules this is the \\
           &          & \hphantom{\mybullet} primitive modulus polynomial of degree $m$ \\
           &          & \hphantom{\mybullet} (the code uses a table of default polynomials) \\[1mm]
\texttt{alpha}
           & optional & \mybullet no effect for lattice rules, $\alpha = 1$ \\
           &          & \mybullet integer interlacing factor for polynomial lattice rules, $\alpha \ge 2$ \\
           &          & \hphantom{\mybullet} (defaults to $\alpha=2$) \\[1mm]
\texttt{a1} & optional, defaults to 0 & integer offset for factorial \\
            &                         & e.g., in our analysis, $a_1=0$ for single-level and \\
            &                         & $a_1=5$ for multi-level \\[1mm]
\texttt{a2} & optional, defaults to 1 & scaling in the product (can be Python expression) \\
            & & e.g., in our analysis, $a_2=1$ for uniform and \\
            & & $a_2=1/\ln 2$ for single-level lognormal and \\
            & & $a_2=2$ for multi-level lognormal \\[1mm]
\texttt{a3} & optional, defaults to 0 & boundary behaviour for the lognormal case \\
            & & ($a_3=0$ means we are in the uniform case) \\
            & & e.g., in our analysis, $a_3=1$ for single-level lognormal and \\
            & & $a_3=9$ for multi-level lognormal \\[1mm]
\texttt{d1} & optional, defaults to 1 & extra power on the factorial factor \\
            &                         & ($d_1=0$ implies product weights)
            \\[1mm]
\texttt{d2} & optional, defaults to 2 & decay of the $\Bj$ sequence, $d_2>1$ \\[1mm]
\texttt{b}  & optional, defaults to $cj^{-d_2}$ & the $\Bj$ sequence as a Python expression, see text \\
            &                                  & (alternatively as numerical values through a file with \verb+b_file+) \\[1mm]
\texttt{c}  & optional, defaults to 1 & in case no \texttt{b} and \texttt{b\_file} are given, $\Bj$ is set to $cj^{-d_2}$ \\[1mm]
\verb+b_file+ & optional & file name containing numerical values for the sequence $\Bj$ \\[1mm]
\texttt{out}  & optional & output directory to write results to \\
\hline
\end{tabular}}
\end{table}

The summary of bounds from the analysis of Section~\ref{sec:error}
corresponds to taking $d_1 = 1$. For randomly shifted lattice rules the
order of convergence is limited to~$1$ and thus $\alpha = 1$. For
interlaced polynomial lattice rules we need $\alpha \ge 2$. The uniform
case corresponds to taking $a_3 = 0$. The lognormal case corresponds to
taking $a_3 > 0$, with $a_3=1$ and $a_3=9$ for the single-level and
multi-level algorithms, respectively.
Our analysis lead to $a_1=0$ and $a_1=5$ for the single-level and
multi-level algorithms, respectively. (Following the proof arguments in
\cite{KSS15,DKLS} we could set $a_1=3$ for the multi-level algorithms in
the uniform case, but the sequence $\Bj$ is defined differently.) We have
$a_2 = 1$ in the uniform case, while in the lognormal case we have
$a_2=1/\ln 2$ and $a_2=2$ for the single-level and multi-level
algorithms, respectively. To cater for other potential integrands
which satisfy the generalized bound \eqref{eq:general-bound}, our scripts
can take general values of $a_1$, $a_2$, $a_3$ and $d_1$ as input.

To specify the sequence $\Bj$ the user has two main options. Either the
user provides a Python expression (with access to variables \texttt{j} and
\texttt{v}, to stand for the values of $j$ and $\nu_j$) as the argument to
command line option \texttt{b}, or the user provides the name of an input
file containing numerical values for each of the $\Bj$ by means of the
\texttt{b\_file} option. (Other possibilities are available, including
a configuration file, but are not discussed here for conciseness.)

We remark that the analysis in Section~\ref{sec:error} takes into account
the truncation from infinite dimensions to $s$~dimensions. Therefore it is
essential to have an idea of the summability of the infinite
sequence~$\Bj$. As before we are interested in the value of $p_* \in
(0,1)$ for which $\sum_{j=1}^\infty \Bj^{\;\,{p_*}} < \infty$ and we would
like $p_*$ to be as small as possible. Here it is more convenient to work
with the reciprocal value, denoted by $d_2>1$, and we call this the
``decay'' of the sequence~$\Bj$.

The theoretical QMC convergence rate in the context of PDE problems, with
the implied constant independent of $s$, is roughly of order
$n^{-\min(1,d_2-1/2)}$ for randomly shifted lattice rules, and
$n^{-\min(\alpha,d_2)}$ for interlaced polynomial lattice rules with
interlacing factor $\alpha\ge 2$.

The analysis in Section~\ref{sec:error} can be extended to handle the
bound~\eqref{eq:general-bound} with a general exponent $d_1$ on the
factorial factor, provided that $d_2 > d_1$ (to ensure that the implied
constant in the error estimate is bounded independently of $s$). The case
$d_1 = 0$ will lead to product weights $\gamma_\setu$, in which case the
CBC construction has a lower cost.

\subsection{Constructing randomly shifted lattice rules}

A synopsis of how to call the Python script to construct a randomly
shifted lattice rule with $n=2^m$ points is

{\small\begin{verbatim}
./lat-cbc.py --s={s} --m={m} [--a1={a1}] [--a2={a2}] [--a3={a3}] --d2={d2}     \
             [--b="{bound-function}" | --b_file={file_name} | --c={c}]
\end{verbatim}}

In particular, the default value of $a_3 = 0$ specifies the
uniform case, while any value of $a_3>0$ specifies the lognormal case.
These two cases correspond to the different function space settings in
Subsections~\ref{sec:sob} and~\ref{sec:Rs}.

The construction script will automatically choose the parameter $\lambda
\in (1/2,1]$ from Theorems~\ref{thm:QMC1} and~\ref{thm:QMC2} as (slightly
different from \eqref{eq:lambda-p0} and \eqref{eq:lambda-p1})
\begin{equation} \label{eq:lambda-delta}
  \frac1{2\lambda}
  \,=\,
  \begin{cases}
  1-\delta & \mbox{if } d_2\ge 3/2-\delta\,, \\
  d_2-\tfrac12 & \mbox{if } d_2 \le 3/2-\delta\,,
  \end{cases}
\end{equation}
with $\delta = 0.125$, which yields the theoretical convergence rate of
order $n^{-(1-\delta)}$ for $d_2\ge 3/2-\delta$ and $n^{-(d_2-1/2)}$ for
$d_2\le 3/2-\delta$. Therefore a correct value of $d_2$ should be
provided. (The value of $\delta$ can also be changed by a command
line argument.)

For the uniform case the script will use the weights (see
\eqref{eq:weight1} and \eqref{eq:weight2})
\[
  \gamma_\setu \,=\,
  \bigg(\bigg(\frac{(|\setu|+a_1)!}{a_1!}\bigg)^{d_1}\,
  \prod_{j\in\setu} \frac{a_2\Bj}{\sqrt{\rho(\lambda)}} \bigg)^{2/(1+\lambda)}\,.
\]
For the lognormal case the script will set the parameters $\alpha_j$ in
the weight functions \eqref{eq:psi} to be (see \eqref{eq:alpha3} and
\eqref{eq:alpha4})
\[
  \alpha_j \,=\, \frac{1}{2} \bigg(
  a_3\Bj + \sqrt{(a_3\Bj)^2+1-\frac{1}{2\lambda}}\bigg)\,,
\]
and use the weights (see \eqref{eq:weight3} and \eqref{eq:weight4})
\[
  \gamma_\setu =
  \bigg(\bigg(\frac{(|\setu|+a_1)!}{a_1!}\bigg)^{d_1}\, \prod_{j\in\setu} \frac{a_2\Bj}
  {2\,\exp((a_3\Bj)^2/2)\Phi(a_3\Bj)\, \sqrt{(\alpha_j-a_3\Bj)\,\rho_j(\lambda)}}
  \bigg)^{2/(1+\lambda)}\,.
\]

We give some examples on how to call the script:

{\small\begin{verbatim}
## uniform case, 100-dimensional rule, 2^10 points and with specified bounds b:
./lat-cbc.py --s=100 --m=10 --d2=3 --b="0.1 * j**-3 / log(j+1)"


## as above, but multi-level and with bounds from file:
./lat-cbc.py --s=100 --m=10 --a1=5 --d2=3 --b_file=bounds.txt

## lognormal case, 100-dimensional rule, 2^10 points and with algebraic decay:
./lat-cbc.py --s=100 --m=10 --a2="1/log(2)" --a3=1 --d2=3 --c=0.1

## as above, but multi-level and with bounds from file:
./lat-cbc.py --s=100 --m=10 --a1=5 --a2=2 --a3=9 --d2=3 --b_file=bounds.txt
\end{verbatim}}

This will produce several files in the output directory. The most
important one is \texttt{z.txt} which contains the generating vector. These
points then need to be randomly shifted for the theory to apply. In the
lognormal case, the randomly shifted points should be mapped to
$\bbR^s$ by applying the inverse of the cumulative normal distribution
function component-wise.

Codes are available in Python, Matlab/Octave and C++ to generate lattice
points. An example usage in Matlab is given below:

{\small\begin{verbatim}
load z.txt                     % load generating vector
latticeseq_b2('init0', z)      % initialize the procedural generator
Pa = latticeseq_b2(20, 512);   % first 512 20-dimensional points
Pb = latticeseq_b2(20, 512);   % next 512 20-dimensional points
\end{verbatim}}

With respect to the multi-level algorithm there are two important features
of these lattice rules: they are lattice sequences in terms of the number
of points, and they are constructed by a component-by-component
algorithm which allows a rule constructed for $s$~dimensions to be used
for a lower number of dimensions. This means the construction only has to
be done once for the maximum number of points $\max_{0\le\ell\le L}
n_\ell$ and the maximum number of dimensions $s_L$, since the parameters
in~\eqref{eq:general-bound} are the same for all levels.

As we have already mentioned, for the lognormal case there is work in
progress on the analysis with a different choice of weight functions
$\vpsi_j$, see Theorem~\ref{thm:QMC2b}, and this would yield a different
choice of weights $\bsgamma_\setu$. We may provide codes for this
alternative setting at a later time.

\subsection{Constructing interlaced polynomial lattice rules}

A synopsis of how to call the Python script to construct an interlaced
polynomial lattice rule with $n=2^m$ and interlacing factor $\alpha\ge 2$
is

{\small\begin{verbatim}
./polylat-cbc.py --s={s} --m={m} --alpha={alpha} [--a1={a1}] [--a2={a2}]       \
             [--b="{bound-function}" | --b_file={file_name} | --d2={d2} --c={c}]
\end{verbatim}}

The construction script will use the weights (see \eqref{eq:weight5} and
\eqref{eq:weight6})
\[
 \gamma_\setu
 \,=\,
 \sum_{\bsnu_\setu \in \{1:\alpha\}^{|\setu|}}
 \bigg(\frac{(|\bsnu_\setu|+a_1)!}{a_1!}\bigg)^{d_1}\,
 \prod_{j\in\setu} \big(2^{\delta(\nu_j,\alpha)} (a_2\Bj)^{\nu_j}\big)\,.
\]
We expect the theoretical convergence rate to be of order
close to $n^{-\min(\alpha,d_2)}$.

We give some examples on how to call the script:

{\small\begin{verbatim}
## 100-dimensional rule, 2^10 points, interlacing 3 and with specified bounds b:
./polylat-cbc.py --s=100 --m=10 --alpha=3 --b="0.1 * j**-3 / log(j+1)"

## as above, but multi-level and with bounds from file:
./polylat-cbc.py --s=100 --m=10 --alpha=3 --a1=5 --b_file=bounds.txt
\end{verbatim}}

Several files will be saved into the output directory.
The most important one is \texttt{Bs.col} which contains the
generating matrices of the interlaced polynomial lattice rule.
(Also the generating matrices of the non-interlaced polynomial
lattice rule are available in the file \texttt{Cs.col}.)

Codes are available in Python, Matlab/Octave and C++ to generate these
points. For interlacing to work correctly, the product $\alpha m$ should
be no more than the number of available bits. We note that the IEEE double
precision type only has 53~bits available and Matlab uses this type to do
its calculations. As a compromise (which comes with no guarantee) we load
the generating matrices truncated to 53~bits precision, which are
available in the file \texttt{Bs53.col}. Similarly, to cater for the
extended long double precision in C++ we provide a file \texttt{Bs64.col}.
For instance, with interlacing factor $4$ we can have up to $2^{13}$
points in Matlab and up to $2^{16}$ points in C++ using long double. The Python point
generator is implemented such that it can use arbitrary precision. One can
also change the C++ generator to use arbitrary precision.

Below we give an example to illustrate how to feed the output from the
construction script into the point generator. At the same time we
experiment on the effect of truncating the generating matrices. First
construct the generating matrices using the Python script and then save
these points (in long double precision) to the file \texttt{points.txt}
using the C++ example program:

{\small\begin{verbatim}
./polylat-cbc.py --s=10 --m=15 --alpha=4 --b="0.1 * j**-4" --out=.
./digitalseq_b2g <Bs64.col >points.txt
\end{verbatim}}

\noindent In this example of $2^{15}$ points with interlacing factor $4$,
we need $4\times 15 = 60$ bits of precision, which can be realized in full
in C++ using long double. Now we use the generating matrices
(truncated to 53~bits) in Matlab, make a plot, and compare these points to the
full precision points we just created on the command line:

{\small\begin{verbatim}
load Bs53.col                      % load generating matrices
s = 10; m = 15;
digitalseq_b2g('init0', Bs53)      % initialize the procedural generator
Pa = digitalseq_b2g(s, pow2(m-1)); % first half of the points
Pb = digitalseq_b2g(s, pow2(m-1)); % second half of the points
s1 = 2; s2 = 10;                   % pick some dimensions to show
plot(Pa(s1,:), Pa(s2,:), 'b.', Pb(s1,:), Pb(s2,:), 'r.')
axis([0 1 0 1]); axis square

load points.txt  % compare with the C++ generated points in long double
points = points';
Pc = points(:,1:pow2(m-1));
Pd = points(:,pow2(m-1)+(1:pow2(m-1)));
norm(Pc - Pa)    % this should be in the order of 1e-14
norm(Pd - Pb)    % and this as well...
\end{verbatim}}

\noindent Here the effect of truncating the generating matrices appears to
be empirically insignificant, but the higher order QMC convergence
theory no longer applies and there is no guarantee how well they would
perform in practice.

\subsection{Generating interlaced Sobol$'$ sequences}

We note that since the point generators operate using generating matrices,
they can be used to generate any other digital sequence, interlaced or
not. On the website we provide the generating matrices for an
implementation of the Sobol$'$ sequence from \cite{JK08} with $21201$
dimensions (as the file \texttt{sobol\_Cs.col}), as well as the
generating matrices for interlaced Sobol$'$ sequences with interlacing
factor $\alpha = 2,3,4,5$ (e.g., \texttt{sobol\_alpha3\_Bs53.col}).
An example usage in Matlab to generate the points is given below:

{\small\begin{verbatim}
load sobol_alpha3_Bs53.col                 % load generating matrices
digitalseq_b2g('init0', sobol_alpha3_Bs53) % initialize the procedural generator
Pa = digitalseq_b2g(10, 1024);             % first 1024 10-dimensional points
Pb = digitalseq_b2g(10, 1024);             % next 1024 10-dimensional points
\end{verbatim}}

\section{Concluding remarks}\label{sec:conc}

In this article we gave a survey of the results from
\cite{KSS12,KSS15,GKNSSS15,KSSSU,DKLNS14,DKLS} on the application of QMC
methods to PDEs with random coefficients, in a unified view. We outlined
three weighted function space settings for analyzing randomly shifted
lattice rules (first order) in both the uniform and lognormal cases, and
interlaced polynomial lattice rules (higher order) in the uniform case. At
present there is no QMC theory that can give higher order convergence for
the lognormal case.

We summarized the error analysis for single-level and multi-level
algorithms based on these QMC methods, in conjunction with FE methods and
dimension truncation. The key step of the analysis is to obtain bounds on
the appropriate weighted norm of the integrand, i.e., $G(u^s_h)$ for the
single-level algorithm and $G(u^s_{h_\ell}-u^s_{h_{\ell-1}})$ for the
multi-level algorithm. We discussed the strategy to obtain suitable
weights $\gamma_\setu$ for the function space setting and arrived at
weights of POD or SPOD form. These weights are to be fed into the CBC
construction of QMC points tailored to the PDE problems. This survey is
augmented with code to construct such tailored QMC points and we explained
how to do this in Section~\ref{sec:comp}.

The combination of a particular family of QMC methods with a specific
function space setting, and the careful designing of POD or SPOD weights
$\gamma_\setu$, means that we obtain QMC error bounds that are independent
of the truncation dimension, while optimizing on the theoretical
convergence rates under minimal assumptions on the PDE problems. Indeed,
we could consider other classes of QMC methods, or use the same QMC
methods but construct them with weights that are not as prescribed here,
however, we might not achieve the same theoretical error bounds. For
example, we recall from Subsection~\ref{sec:fo-sl-unif} that in the
uniform case with randomly shifted lattice rules we can achieve nearly
first order convergence if $p_0\le 2/3$, with $p_0$ from
Assumption~\ref{U5}. As pointed out in \cite{KSS12}, we could consider
randomly shifted lattice rules constructed with product weights (instead
of POD weights), or the deterministic lattice rules constructed following
\cite{Joe06}, or Niederreiter and Sobol$'$ sequences following the
analysis in \cite{Wan02}, but then to achieve nearly first order
convergence we would require, respectively, $p_0\le 1/2$, $p_0\le 1/2$,
and $p_0\le 1/3$, meaning that a stronger assumption on the PDE problem is
required to achieve the same convergence rate. More strikingly, we recall
from Subsection~\ref{sec:ho-sl-unif} that in the uniform case we can
construct interlaced polynomial lattice rules with interlacing factor $2$
to achieve first order convergence when $p_0<1$, which is a much weaker
condition on the PDE problem.

We already pointed out that the uniform framework can be extended to
general affine parametric operator equations, see
\cite{Sch12,DKLNS14,DKLS}. Thus the QMC strategy in this article applies
to a wide range of PDE problems including, e.g., stationary and
time-dependent diffusion in random media \cite{CDS10}, wave propagation
\cite{HS13}, parametric nonlinear PDEs \cite{CCS15}, and optimal control
problems for uncertain systems \cite{KunSch13}.

There may be additional properties of the functions $\psi_j$ in
\eqref{eq:axy-unif} or $\xi_j$ in \eqref{eq:axy-logn} that could be
exploited to improve the results. For example, in \cite{KSS15} a special
orthogonality property for multiresolution function systems was used to
obtain a better dimension truncation estimate for the multi-level error
analysis. In a different direction, \cite[Theorem~3.2]{DKLNS14} pointed
out that if we were to replace the $|\bsnu|! = (\sum_{j\ge1} \nu_j)!$
factor in Lemma~\ref{lem:regu1} by $\bsnu! = \prod_{j\ge 1} \nu_j!$, then
our analysis would yield product weights $\gamma_\setu$ rather than POD or
SPOD weights, and this would significantly reduce the CBC construction
cost.

We have assumed in this article that the cost for our single-level
algorithm is
\begin{equation} \label{eq:cost-kl}
  \calO(n\,s\,h^{-d})
\end{equation}
operations, based on $n$ instances of FE discretizations where the cost
for assembling the stiffness matrix is $\calO(s\,h^{-d})$ operations, with
the $\calO(s)$ factor coming from the number of KL terms in the
coefficient $a(\bsx,\bsy)$. An analogous cost model is assumed for the
multi-level algorithm, and we argued that the values of $n_\ell$,
$s_\ell$, $h_\ell$ should be chosen by minimizing the error for a fixed
cost. Understandably this latter optimization depends crucially on the
cost model. Below we discuss two related strategies to reduce cost.

We mentioned the ``circulant embedding'' strategy for the lognormal case,
see \cite{GKNSS11,GKNSS}, but did not go into any details in this article.
Roughly speaking, the idea is to sample the random field only at a
discrete set of grid points with respect to the covariance matrix, where
the number of grid points is of the same order as the number of FE nodes.
Then the problem of generating samples turns into a matrix factorization
problem which can be done in $\calO(h^{-d} \,(\log h^{-d}))$ operations
using FFT, by embedding the covariance matrix in a larger but circulant
matrix (some further padding may be required to ensure positive
definiteness). With this strategy the cost becomes
\[
  \calO(n\,h^{-d}\,(\log h^{-d}))
\]
operations, where we effectively replaced the $\calO(s)$ factor in
\eqref{eq:cost-kl} by $\calO(\log h^{-d})$.

The other strategy to reduce cost is the ``fast QMC matrix-vector
multiplication'', see \cite{DKLS15}, which exploits a certain structure in
the QMC point set. By choosing suitable QMC point sets, and by formulating
the QMC quadrature computation as a matrix-vector multiplication with a
circulant matrix obtained by indexing the QMC points in a certain way, the
cost becomes
\[
  \calO(n\,(\log n)\,h^{-d})
\]
operations by using FFT, which essentially means that the $\calO(s)$
factor in \eqref{eq:cost-kl} is replaced by $\calO(\log n)$.
Unfortunately, this strategy does \emph{not} work with randomly shifted
lattice rules or interlaced polynomial lattice rules considered in this
article, because both randomization and interlacing destroy the required
structure in the QMC point sets to yield a circulant matrix. However, this
strategy is compatible with ``tent transformed'' lattice rules or
polynomial lattice rules, which are deterministic QMC rules that can
achieve nearly first or second order convergence, see \cite{DNP14,GSY16}.
The error analysis of applying tent transformed lattice rules or
polynomial lattice rules for the PDE problems is work in progress.

The above three strategies for cost reduction (namely, multi-level
algorithms, circulant embedding, and fast QMC matrix-vector
multiplication) are \emph{not} mutually exclusive and they could
potentially be combined to have a compounding effect in reducing cost.
There is also a generalization of the multi-level concept called
``multi-index'' \cite{HNT16}, which is in some sense related to ``sparse
grid techniques'' \cite{BG04}. Note that each strategy has its
prerequisite: multi-level algorithm requires stronger regularity
assumptions on the PDE problem, circulant embedding requires stationary
covariance functions, while fast QMC matrix-vector multiplication requires
a certain structure in the QMC point set. It would be interesting to see
which strategy or which combination of strategies yields the most
effective reduction in cost under different scenarios.

Now we make some remarks on theory versus practice. Although the careful
tuning of weights $\gamma_\setu$ played a significant role in our analysis
and affected the theoretical QMC convergence rates, the outcomes from
numerical experiments so far have been inconclusive. We have seen that
some ``off the shelf'' lattice rules (i.e., lattice rules constructed with
product weights chosen to have some generic algebraic or geometric decay)
perform just as well as those lattice rules which are tailored to the PDE
problems. On the other hand, we have also seen that some badly tuned
lattice rules (e.g., when the weights $\gamma_\setu$ are badly scaled) can
perform poorly.

We also noted that the theoretical convergence rates are not always
reflected in the computations. In the lognormal computations in
\cite{GKNSSS15} we see that the convergence rates are not so much
influenced by the smoothness parameter $\nu$ of the Mat\'ern covariance
function as the theory predicted. Rather, it is the variance and the
correlation length that affect the empirical convergence rates.

The numerical experiments for the lognormal case in \cite{GKNSS11} were
obtained using randomly digitally shifted Sobol$'$ points with no theory,
yet the results were very encouraging. One could also try tent transformed
Sobol$'$ points or interlaced Sobol$'$ points. We suspect that they may
work reasonably well in practice, even though at present we are lacking a
strong supporting theory. A brief explanation on how to generate
interlaced Sobol$'$ sequences from the code can be found in
Section~\ref{sec:comp}.

Finally, if the story in this survey article sounds incomplete, we hope
the reader will understand that we are trying to tell a story that is
changing underneath us, even as we write. We live in interesting times!

\section{Appendix: selected proofs}\label{sec:proof}

In this section we provide the proofs for
Lemmas~\ref{lem:regu1}--\ref{lem:regu4c}. For simplicity of presentation,
in the proofs we will often omit the arguments $\bsx$ and $\bsy$ in our
notation. We start by collecting some identities and estimates that we
need for the proofs.

We will make repeated use of the Leibniz product rule
\begin{align} \label{eq:leibniz}
  \partial^\bsnu (AB) \,=\,
  \sum_{\bsm\le\bsnu} \sbinom{\bsnu}{\bsm} (\partial^{\bsm} A)\, (\partial^{\bsnu-\bsm} B)\,,
\end{align}
and the identity
\begin{align} \label{eq:nabla}
  \nabla\cdot (A\, \nabla B) \,=\, A\, \Delta B + \nabla A \cdot \nabla B\,.
\end{align}
We also need the combinatorial identity
\begin{align} \label{eq:ident0}
  \sum_{\satop{\bsm\le\bsnu}{|\bsm|=i}} \sbinom{\bsnu}{\bsm}
  \,=\, \sbinom{|\bsnu|}{i}\,,
\end{align}
which follows from a simple counting argument (i.e., consider the number
of ways to select $i$ distinct balls from some baskets containing a total
number of $|\bsnu|$ distinct balls). The identity \eqref{eq:ident0} is
used to establish the following identities
\begin{align}
 \sum_{\bsm\le\bsnu} \sbinom{\bsnu}{\bsm}\, |\bsm|!\,|\bsnu-\bsm|! &\,=\, (|\bsnu|+1)!\,, \label{eq:ident1}
 \\
 \sum_{\bsm\le\bsnu} \sbinom{\bsnu}{\bsm} |\bsm|! \,(|\bsnu-\bsm|+1)!
 &\,=\, \frac{(|\bsnu|+2)!}{2}\,, \label{eq:ident2}
 \\
 \sum_{\bsm\le\bsnu} \sbinom{\bsnu}{\bsm} \frac{(|\bsm|+2)!}{2}\,\frac{(|\bsnu-\bsm|+2)!}{2}
 &\,=\, \frac{(|\bsnu|+5)!}{120}\,. \label{eq:ident5}
\end{align}
Additionally, we need the recursive estimates in the next two lemmas. The
proofs can be found in \cite{DKLS} and \cite{KSSSU}, respectively.

\begin{lemma}
\label{lem:recur1}
Given a sequence of non-negative numbers $\bsb=(b_j)_{j\in\bbN}$, let
$(\bbA_\bsnu)_{\bsnu\in\indx}$ and $(\bbB_\bsnu)_{\bsnu\in\indx}$ be
non-negative numbers satisfying the inequality
\begin{align*}
  \bbA_\bsnu
  \,\le\, \sum_{j\in\supp(\bsnu)} \nu_j\,b_j\, \bbA_{\bsnu-\bse_j} + \bbB_\bsnu
  \quad\mbox{for any $\bsnu\in\indx$ \quad\textnormal{(}including $\bsnu=\bszero$\textnormal{)}}\,.
\end{align*}
Then
\begin{align*}
  \bbA_\bsnu
  \,\le\, \sum_{\bsm\le\bsnu} \sbinom{\bsnu}{\bsm}\, |\bsm|!\,
  \bsb^\bsm\,\bbB_{\bsnu-\bsm}
  \quad\mbox{for all $\bsnu\in\indx$}\,.
\end{align*}
The result holds also when both inequalities are replaced by equalities.
\end{lemma}

\begin{lemma}
\label{lem:recur2}
Given a sequence of non-negative numbers $\bsbeta=(\beta_j)_{j\in\bbN}$,
let $(\bbA_\bsnu)_{\bsnu\in\indx}$ and $(\bbB_\bsnu)_{\bsnu\in\indx}$ be
non-negative numbers satisfying the inequality
\begin{align*}
  \bbA_\bsnu
  \,\le\, \sum_{\satop{\bsm\le\bsnu}{\bsm\ne\bsnu}} \sbinom{\bsnu}{\bsm} \bsbeta^{\bsnu-\bsm} \bbA_\bsm + \bbB_\bsnu
  \quad\mbox{for any $\bsnu\in\indx$ \quad\textnormal{(}including $\bsnu=\bszero$\textnormal{)}}\,.
\end{align*}
Then
\begin{align*}
  \bbA_\bsnu
  \,\le\, \sum_{\bsk\le\bsnu} \sbinom{\bsnu}{\bsk} \Lambda_{|\bsk|}\, \bsbeta^\bsk\, \bbB_{\bsnu-\bsk}
  \quad\mbox{for all $\bsnu\in\indx$}\,,
\end{align*}
where the sequence $(\Lambda_n)_{n\ge 0}$ is defined recursively by
\begin{align} \label{eq:lambda}
  \Lambda_0 \,:=\, 1 \quad\mbox{and}\quad \Lambda_n \,:=\, \sum_{i=0}^{n-1} \sbinom{n}{i} \Lambda_i
  \quad\mbox{for all $n\ge 1$}\,.
\end{align}
The result holds also when both inequalities are replaced by equalities.
Moreover, we have
\begin{align} \label{eq:lambda-bound}
  \Lambda_n \,\le\, \frac{n!}{\alpha^n} \quad\mbox{for all \ $n\ge
  0$ \
  and \ $\alpha\le\ln 2=0.69...$}.
\end{align}
\end{lemma}


\subsection*{\textbf{Proof of Lemma~\ref{lem:regu1}}}

This result was proved in \cite{CDS10}. Since the same proof strategy is
used repeatedly in subsequent more complicated proofs, we include this
relatively simple proof as a first illustration.

 Let $f\in V^*$ and $\bsy\in U$. We prove this result by
induction on $|\bsnu|$. For $\bsnu=\bszero$, we take $v = u(\cdot,\bsy)$
in~\eqref{eq:weak} to obtain
\begin{align*}
 \int_D a(\bsx,\bsy)\, |\nabla u(\bsx,\bsy)|^2 \,\rd\bsx
  \,=\, \int_D f(\bsx)\,u(\bsx,\bsy)\,\rd\bsx\,,
\end{align*}
which leads to
\begin{align*}
 a_{\min}\, \|u(\cdot,\bsy)\|_V^2
  \,\le\, \|f\|_{V^*} \|u(\cdot,\bsy)\|_V
  \quad\implies\quad
  \|u(\cdot,\bsy)\|_V \,\le\, \frac{\|f\|_{V^*}}{a_{\min}}\,,
\end{align*}
as required (see also \eqref{eq:apriori-unif}).

Given any multi-index $\bsnu$ with $|\bsnu|\ge 1$, suppose that the result
holds for any multi-index of order $\le |\bsnu|-1$. Applying the mixed
derivative operators $\partial^\bsnu$ to the variational
formulation~\eqref{eq:weak}, recalling that $f$ is independent of $\bsy$,
and using the Leibniz product rule~\eqref{eq:leibniz}, we obtain the
identity (suppressing $\bsx$ and $\bsy$ in our notation)
\begin{align*}
  \int_D \bigg( \sum_{\bsm\le\bsnu} \sbinom{\bsnu}{\bsm}
  (\partial^{\bsm} a)\,
  \nabla (\partial^{\bsnu-\bsm} u) \cdot \nabla z \bigg)\,\rd\bsx
  \,=\, 0
  \qquad\mbox{for all}\quad
  z\in V\,.
\end{align*}
Observe that due to the linear dependence of $a(\bsx,\bsy)$ on the
parameters $\bsy$, the partial derivative $\partial^{\bsm}$ of $a$ with
respect to $\bsy$ satisfies
\begin{align} \label{eq:diff-unif}
  \partial^{\bsm} a(\bsx,\bsy)
  \,=\,
  \begin{cases}
  a(\bsx,\bsy) & \mbox{if } \bsm = \bszero, \\
  \psi_j(\bsx) & \mbox{if } \bsm = \bse_j, \\
  0 & \mbox{otherwise}.
  \end{cases}
\end{align}
Taking $z = \partial^\bsnu u(\cdot,\bsy)$ and separating out the $\bsm =
\bszero$ term, we obtain
\begin{align*}
  \int_D a\, |\nabla (\partial^\bsnu u)|^2 \,\rd\bsx
  \,=\, - \sum_{j\in\supp(\bsnu)} \nu_j\,
  \int_D \psi_j\,\nabla (\partial^{\bsnu-\bse_j}u) \cdot
  \nabla (\partial^{\bsnu} u)\,\rd\bsx\,,
\end{align*}
which yields
\begin{align*}
  a_{\min}\, \|\nabla (\partial^\bsnu u) \|_{L^2}^2
  &\,\le\, \sum_{j\ge 1}\nu_j\,\|\psi_j\|_{L^\infty}
  \|\nabla (\partial^{\bsnu-\bse_j} u)\|_{L^2}\,
  \|\nabla (\partial^{\bsnu} u)\|_{L^2},
\end{align*}
and therefore
\begin{align*}
  \|\nabla (\partial^\bsnu u) \|_{L^2}
  \,\le\, \sum_{j\ge 1}\nu_j\, b_j\, \|\nabla (\partial^{\bsnu-\bse_j} u)\|_{L^2}\,,
\end{align*}
where we used the definition of $b_j$ in \eqref{eq:bj-unif}. The induction
hypothesis then gives
\begin{align*}
  \|\nabla (\partial^\bsnu u) \|_{L^2}
  &\,\le\, \sum_{j\ge 1}\nu_j\, b_j\, |\bsnu-\bse_j|!\, \bsb^{\bsnu-\bse_j}\,\frac{\|f\|_{V^*}}{a_{\min}}
  \,=\, |\bsnu|!\,\bsb^\bsnu\,\frac{\|f\|_{V^*}}{a_{\min}}\,.
\end{align*}
This completes the proof. \qed


\subsection*{\textbf{Proof of Lemma~\ref{lem:regu2a}}}

This result corresponds to \cite[Theorem~6]{KSS15}. Here we take a more
direct route with the proof and the bound depends on the sequence
$\overline{\bsb}$ which is simpler than the sequence in \cite{KSS15}
(there the sequence depends on an additional parameter $\kappa\in (0,1]$
and other constants), at the expense of increasing the factorial factor
from $|\bsnu|!$ to $(|\bsnu|+1)!$. For simplicity we consider here the
case $f\in L^2(D)$, but the proof can be generalized to the case $f\in
H^{-1+t}(D)$ for $t\in [0,1]$ as in~\cite{KSS15}.

Let $f\in L^2(D)$ and $\bsy\in U$. For $\bsnu=\bszero$, we apply the
identity~\eqref{eq:nabla} to the strong formulation~\eqref{eq:strong} to
obtain (formally, at this stage, since we do not yet know that $\Delta
u(\cdot,\bsy)\in L^2(D)$)
\begin{align*}
  - a(\bsx,\bsy)\, \Delta u(\bsx,\bsy) \
  \,=\, \nabla a(\bsx,\bsy)\cdot \nabla u(\bsx,\bsy) + f(\bsx),
\end{align*}
which leads to
\begin{align*}
  a_{\min}\, \|\Delta u(\cdot,\bsy)\|_{L^2}
  \,\le\, \|\nabla a(\cdot,\bsy)\|_{L^\infty}\,\|\nabla u(\cdot,\bsy)\|_{L^2} + \|f\|_{L^2}\,.
\end{align*}
Combining this with~\eqref{eq:apriori-unif} gives
\begin{align} \label{eq:base2-unif}
  \|\Delta u(\cdot,\bsy)\|_{L^2}
  \,\le\, \frac{\sup_{\bsz\in U}\|\nabla a(\cdot,\bsz)\|_{L^\infty}}{a_{\min}}\,\frac{\|f\|_{V^*}}{a_{\min}}
  + \frac{\|f\|_{L^2}}{a_{\min}}
  \,\le\, C\,\|f\|_{L^2}\,,
\end{align}
where we could take
\begin{align*}
  C \,:=\, C_{\rm emb}\left(
  \frac{\sup_{\bsz\in U}\|\nabla a(\cdot,\bsz)\|_{L^\infty}}{a_{\min}^2} + \frac{1}{a_{\min}}\right)\,,
  \quad\mbox{with}\quad
  C_{\rm emb} \,:=\, \sup_{f\in L^2(D)} \frac{\|f\|_{V^*}}{\|f\|_{L^2}}\,.
\end{align*}
Thus the result holds for $\bsnu = \bszero$ (see also
\eqref{eq:delta-unif}).

For $\bsnu\ne \bszero$, we apply $\partial^\bsnu$ to the strong
formulation~\eqref{eq:strong} and use the Leibniz product
rule~\eqref{eq:leibniz} to obtain (suppressing $\bsx$ and $\bsy$)
\begin{align*}
  \nabla \cdot
  \bigg(\sum_{\bsm \le\bsnu} \sbinom{\bsnu}{\bsm}
  (\partial^{\bsm} a)\, \nabla (\partial^{\bsnu-\bsm} u) \bigg)
 \,=\, 0\,.
\end{align*}
Using again~\eqref{eq:diff-unif} and separating out the $\bsm = \bszero$
term yield the following identity
\begin{align*}
  \nabla\cdot(a\nabla(\partial^\bsnu u))
  \,=\, - \nabla \cdot \bigg(\sum_{j\in\supp(\bsnu)} \nu_j\,\psi_j(\bsx)\,
  \nabla (\partial^{\bsnu-\bse_j} u) \bigg)\,.
\end{align*}
Applying the identity~\eqref{eq:nabla} to both sides yields (formally)
\begin{align*}
  a\,\Delta (\partial^\bsnu u) + \nabla a \cdot \nabla(\partial^\bsnu u)
  \,=\, - \sum_{j\in\supp(\bsnu)} \nu_j\,
  \bigg(\psi_j\, \Delta(\partial^{\bsnu-\bse_j} u)
  \, + \, \nabla \psi_j\cdot \nabla (\partial^{\bsnu-\bse_j} u) \bigg)\,.
\end{align*}
In turn, we obtain
\begin{align*}
  a_{\min} \,\|\Delta (\partial^\bsnu u)\|_{L^2}
  &\,\le\,
  \sum_{j\in\supp(\bsnu)} \nu_j\,
  \bigg(\|\psi_j\|_{L^\infty}\, \|\Delta(\partial^{\bsnu-\bse_j} u)\|_{L^2}
   + \|\nabla \psi_j\|_{L^\infty} \, \|\nabla (\partial^{\bsnu-\bse_j} u)\|_{L^2} \bigg) \\
  &\qquad + \|\nabla a\|_{L^\infty} \,\|\nabla(\partial^\bsnu u)\|_{L^2}\,,
\end{align*}
which leads to
\begin{align*}
  \underbrace{\|\Delta (\partial^\bsnu u)\|_{L^2}}_{\bbA_\bsnu}
  &\,\le\,
  \sum_{j\in\supp(\bsnu)} \nu_j\,b_j\, \underbrace{\|\Delta(\partial^{\bsnu-\bse_j} u)\|_{L^2}}_{\bbA_{\bsnu-\bse_j}}
  \,+\, B_\bsnu\,,
\end{align*}
where we used the definition of $b_j$ in \eqref{eq:bj-unif}, and
\begin{align*}
  B_\bsnu
  \,:=\,
  \sum_{j\in\supp(\bsnu)} \nu_j\, \frac{\|\nabla \psi_j\|_{L^\infty}}{a_{\min}}
  \, \|\nabla (\partial^{\bsnu-\bse_j} u)\|_{L^2}
  + \frac{\|\nabla a\|_{L^\infty}}{a_{\min}} \,\|\nabla(\partial^\bsnu u)\|_{L^2}\,.
\end{align*}

Note that this formulation of $B_\bsnu$ cannot be used as $\bbB_\bsnu$ in
Lemma~\ref{lem:recur1} because the base step $\bbA_\bszero \le B_\bszero$
does not hold. {}From Lemma~\ref{lem:regu1} we can estimate
\begin{align*}
  B_\bsnu
  &\,\le\,
  \sum_{j\in\supp(\bsnu)} \nu_j\, \frac{\|\nabla \psi_j\|_{L^\infty}}{a_{\min}}
  \, |\bsnu-\bse_j|!\, \bsb^{\bsnu-\bse_j}\, \frac{\|f\|_{V^*}}{a_{\min}}
  + \frac{\sup_{\bsz\in U} \|\nabla a(\cdot,\bsz)\|_{L^\infty}}{a_{\min}}
  \,|\bsnu|!\, \bsb^{\bsnu}\, \frac{\|f\|_{V^*}}{a_{\min}} \\
  &\,\le\, C\,|\bsnu|!\, \overline\bsb^{\bsnu}\, \|f\|_{L^2(D)} \,=:\, \bbB_\bsnu\,,
\end{align*}
where we used the definition of $\overline{b}_j\ge b_j$ in
\eqref{eq:barbj-unif}. This definition of $\bbB_\bsnu$ ensures that the
base step $\bbA_\bszero\le \bbB_\bszero$ does hold,
see~\eqref{eq:base2-unif}. Now we apply Lemma~\ref{lem:recur1} to conclude
that
\begin{align*}
  \|\Delta (\partial^\bsnu u)\|_{L^2}
  &\,\le\,
  \sum_{\bsm\le\bsnu} \sbinom{\bsnu}{\bsm}\, |\bsm|!\, \bsb^\bsm\,
  C\,|\bsnu-\bsm|!\, \overline\bsb^{\bsnu-\bsm}\, \|f\|_{L^2}
  \,\le\, C\,(|\bsnu|+1)!\, \overline\bsb^{\bsnu}\, \|f\|_{L^2}\,,
\end{align*}
where we used the identity \eqref{eq:ident1}. This completes the proof.
\qed


\subsection*{\textbf{Proof of Lemma~\ref{lem:regu2b}}}

This result appeared as a technical step in the proof of
\cite[Theorem~7]{KSS15}, but only first derivatives were considered there,
i.e., $\nu_j\le 1$ for all $j$. Here we consider general derivatives, and
we make use of Lemma~\ref{lem:regu2a} instead of \cite[Theorem~6]{KSS15}
so that the sequence $\overline{\bsb}$ is different, the factorial factor
is larger, and we restrict to $f\in L^2(D)$.

Let $f\in L^2(D)$, $\bsy\in U$, and $\bsnu\in\indx$. Galerkin
orthogonality for the FE method yields
\begin{align}
\label{eq:GalerkinOrth}
 \scA(\bsy;u(\cdot,\bsy)-u_h(\cdot,\bsy), z_h)\,= \,0
 \qquad \text{for all} \quad z_h \in V_h\,,
\end{align}
Let $\calI:V\to V$ denote the identity operator and let $\calP_h =
\calP_h(\bsy) :V\to V_h$ denote the parametric FE projection onto $V_h$
which is defined, for arbitrary $w\in V$, by
\begin{align} \label{eq:Ph}
\scA(\bsy; \calP_h(\bsy) w - w, z_h) \,=\, 0
\qquad\text{for all} \quad z_h\in V_h \,.
\end{align}
In particular, we have $u_h = \calP_h  u \in V_h$ and
\begin{align}\label{eq:IdemP}
\calP_h^2(\bsy) \,\equiv\, \calP_h(\bsy) \quad\mbox{on}\quad V_h \,.
\end{align}
Moreover, since $\partial^\bsnu u_h \in V_h$ for every $\bsnu\in \indx$,
it follows from~\eqref{eq:IdemP} that
\begin{align}\label{eq:I-Ph}
(\calI - \calP_h(\bsy))(\partial^\bsnu u_h(\cdot,\bsy)) \,\equiv\, 0
\,.
\end{align}
Thus
\begin{align} \label{eq:new_1}
  \|\nabla \partial^\bsnu (u - u_{h})\|_{L^2}
  &\,=\, \|\nabla \calP_h \partial^\bsnu (u - u_{h}) + \nabla (\calI - \calP_h) \partial^\bsnu (u - u_{h}) \|_{L^2} \nonumber\\
  &\,\le\, \| \nabla \calP_h \partial^\bsnu (u - u_{h})\|_{L^2}
  \,+\, \|\nabla (\calI - \calP_h) \partial^\bsnu u \|_{L^2}
\,.
\end{align}
We stress here that, since the parametric FE projection $\calP_h(\bsy)$
depends on $\bsy$, in general we have $\partial^\bsnu (u(\cdot,\bsy) -
u_{h}(\cdot,\bsy)) \ne (\calI - \calP_h(\bsy)) (\partial^\bsnu
u(\cdot,\bsy))$; this is why we need the estimate \eqref{eq:new_1}.

Now, applying $\partial^\bsnu$ to~\eqref{eq:GalerkinOrth} and recalling
\eqref{eq:diff-unif}, we get for all $z_h\in V_h$,
\begin{align} \label{eq:new_2}
 \int_D a\, \nabla \partial^{\bsnu} (u-u_h)\cdot\nabla z_h \,\rd \bsx
 \,=\, - \sum_{j\in\supp(\bsnu)} \nu_j\,
 \int_D \psi_j \nabla\partial^{\bsnu-\bse_j}(u-u_h)\cdot\nabla z_h \,\rd \bsx\,.
\end{align}
Choosing $z_h = \calP_h \partial^\bsnu (u-u_h)$ and using the
definition~\eqref{eq:Ph} of $\calP_h$, the left-hand side
of~\eqref{eq:new_2} is equal to $\int_D a\, |\nabla \calP_h
\partial^{\bsnu} (u-u_h)|^2 \,\rd \bsx$. Using the Cauchy-Schwarz
inequality, we then obtain
\begin{align*}
 a_{\min}\, \|\nabla \calP_h \partial^{\bsnu} (u-u_h)\|_{L^2}^2 
 &\,\le\, \sum_{j\in\supp(\bsnu)} \nu_j\,
 \|\psi_j\|_{L^\infty} \|\nabla\partial^{\bsnu-\bse_j}(u-u_h)\|_{L^2}\,
 \|\nabla \calP_h \partial^{\bsnu} (u-u_h)\|_{L^2}\,.
\end{align*}
Canceling one common factor from both sides, we arrive at
\begin{align} \label{eq:new_3}
 \|\nabla \calP_h \partial^{\bsnu} (u-u_h)\|_{L^2}
 &\,\le\, \sum_{j\in\supp(\bsnu)} \nu_j\,b_j\, \|\nabla\partial^{\bsnu-\bse_j}(u-u_h)\|_{L^2}\,,
\end{align}
where we used the definition of $b_j$ in \eqref{eq:bj-unif}.
Substituting~\eqref{eq:new_3} into~\eqref{eq:new_1}, we then obtain
\begin{align*}
  \underbrace{\| \nabla \partial^\bsnu (u - u_{h})\|_{L^2}}_{\bbA_\bsnu}
  &\,\le\, \sum_{j\in\supp(\bsnu)} \nu_j\,b_j\,
  \underbrace{\|\nabla\partial^{\bsnu-\bse_j}(u-u_h)\|_{L^2}}_{\bbA_{\bsnu-\bse_j}}
  \,+\, \underbrace{\|\nabla (\calI - \calP_h) \partial^\bsnu u\|_{L^2}}_{\bbB_\bsnu} \,.
\end{align*}
Noting that $\bbA_\bszero = \bbB_\bszero$, we now apply
Lemma~\ref{lem:recur1} to obtain
\begin{align*}
  \| \nabla \partial^\bsnu (u - u_{h})\|_{L^2}
  &\,\le\, \sum_{\bsm\le\bsnu} \sbinom{\bsnu}{\bsm}\, |\bsnu|!\, \bsb^\bsnu\,
  \|\nabla (\calI - \calP_h) \partial^{\bsnu-\bsm} u\|_{L^2}\,.
\end{align*}
Next we use the FE estimate \eqref{eq:FEerror-unif} that for all $\bsy\in
U$ and $w\in H^2(D)$, we have $\|\nabla (\calI - \calP_h) w \|_{L^2}
\,\lesssim\, h\, \|\Delta w\|_{L^2}$. Hence from Lemma~\ref{lem:regu2a} we
obtain
\begin{align*}
  \| \nabla \partial^\bsnu (u - u_{h})\|_{L^2}
  &\,\lesssim\, \sum_{\bsm\le\bsnu} \sbinom{\bsnu}{\bsm}\, |\bsnu|!\, \bsb^\bsnu\,
  h\, \|\Delta (\partial^{\bsnu-\bsm} u)\|_{L^2} \nonumber \\
  &\,\lesssim\, \sum_{\bsm\le\bsnu} \sbinom{\bsnu}{\bsm}\, |\bsnu|!\, \bsb^\bsnu\,
  h\, (|\bsnu-\bsm|+1)!\, \overline\bsb^{\bsnu-\bsm}\, \|f\|_{L^2} \nonumber\\
  &\,\lesssim\, h\, \frac{(|\bsnu|+2)!}{2}\,\overline\bsb^\bsnu\, \|f\|_{L^2}\,,
\end{align*}
where we used the identity~\eqref{eq:ident2}. This completes the proof.
\qed


\subsection*{\textbf{Proof of Lemma~\ref{lem:regu2c}}}

This result generalizes \cite[Theorem~7]{KSS15} from first derivatives to
general derivatives. The proof is based on a duality argument since $G$ is
a bounded linear functional. It makes use of Lemma~\ref{lem:regu2b}, and
therefore the sequence $\overline{\bsb}$ is different, the factorial
factor is larger, and we restrict to $f,G\in L^2(D)$ here.

Let $f,G\in L^2(D)$ and $\bsy\in U$. We define $v^G(\cdot,\bsy)\in V$ and
$v_h^G(\cdot,\bsy)\in V_h$ via the adjoint problems
\begin{align}
  \scA(\bsy; w, v^G(\cdot,\bsy)) &\,=\, G( w )
  &&\text{for all}\quad w\in V\,,
  \label{eq:nitsche_1}
\\
  \scA(\bsy; w_h, v_h^G(\cdot,\bsy)) &\,=\, G ( w_h )
  &&\text{for all}\quad w_h\in V_h\,.
\label{eq:nitsche_2}
\end{align}
Due to Galerkin orthogonality \eqref{eq:GalerkinOrth} for the original
problem, by choosing the test function $w = u(\cdot,\bsy) -
u_h(\cdot,\bsy)$ in~\eqref{eq:nitsche_1}, we obtain
\begin{align} \label{eq:duality}
  G(u(\cdot,\bsy)-u_h(\cdot,\bsy))
 \,=\, \scA(\bsy; u(\cdot,\bsy)-u_h(\cdot,\bsy), v^G(\cdot,\bsy) - v_h^G(\cdot,\bsy))\,.
\end{align}

{}From the Leibniz product rule~\eqref{eq:leibniz}
and~\eqref{eq:diff-unif} we have for $\bsnu\in\indx$
\begin{align*}
  &\partial^\bsnu G(u - u_{h})
  \,=\,
  \int_D \partial^\bsnu \left(a \,\nabla (u
  - u_{h}) \cdot \nabla (v^G - v^G_h) \right)\, \rd\bsx \\
  &\,=\,
  \int_D \sum_{\bsm\le\bsnu} \sbinom{\bsnu}{\bsm} (\partial^{\bsm} a)\,
  \partial^{\bsnu-\bsm} \left( \nabla (u - u_{h}) \cdot \nabla (v^G - v^G_h) \right)\, \rd\bsx \\
  &\,=\,
  \int_D a\, \partial^{\bsnu} \left( \nabla (u - u_{h}) \cdot \nabla (v^G - v^G_h) \right)\, \rd\bsx \\
  &\qquad + \sum_{j\in\supp(\bsnu)} \nu_j\,
  \int_D \psi_j\,\partial^{\bsnu-\bse_j} \left( \nabla (u - u_{h}) \cdot \nabla (v^G - v^G_h) \right)\, \rd\bsx \\
  &\,=\,
  \int_D a\, \sum_{\bsk\le\bsnu} \sbinom{\bsnu}{\bsk} \nabla \partial^\bsk (u- u_{h}) \cdot
  \nabla \partial^{\bsnu-\bsk} (v^G - v^G_h)\, \rd\bsx \\
  &\qquad + \sum_{j\in\supp(\bsnu)} \nu_j\,
  \int_D \psi_j\,\sum_{\bsk\le\bsnu-\bse_j} \sbinom{\bsnu-\bse_j}{\bsk} \nabla \partial^\bsk (u- u_{h}) \cdot
  \nabla \partial^{\bsnu-\bse_j-\bsk} (v^G - v^G_h)\, \rd\bsx\,.
\end{align*}
The Cauchy-Schwarz inequality then yields
\begin{align} \label{eq:split}
  &|\partial^\bsnu G(u - u_{h})|
  \,\le\,
  a_{\max} \sum_{\bsk\le\bsnu} \sbinom{\bsnu}{\bsk} \|\nabla \partial^\bsk (u- u_{h})\|_{L^2} \,
  \|\nabla \partial^{\bsnu-\bsk} (v^G - v^G_h)\|_{L^2} \nonumber\\
  &\qquad + \sum_{j\in\supp(\bsnu)} \nu_j\,
  \|\psi_j\|_{L^\infty}\,\sum_{\bsk\le\bsnu-\bse_j} \sbinom{\bsnu-\bse_j}{\bsk} \|\nabla \partial^\bsk (u- u_{h})\|_{L^2}\,
  \|\nabla \partial^{\bsnu-\bse_j-\bsk} (v^G - v^G_h)\|_{L^2}\,.
\end{align}

We see from (the proof of) Lemma~\ref{lem:regu2b} that
\begin{align} \label{eq:nice_1}
  \| \nabla \partial^\bsk (u - u_{h})\|_{L^2}
  &\,\lesssim\, h\, \frac{(|\bsk|+2)!}{2}\,\overline\bsb^\bsk\, \|f\|_{L^2}\,.
\end{align}
Since the bilinear form $\scA(\bsy;\cdot,\cdot)$ is symmetric and since
the representer $g$ for the linear functional $G$ is in $L^2(D)$, all the
results hold verbatim also for the adjoint problem~\eqref{eq:nitsche_1}
and for its FE discretisation~\eqref{eq:nitsche_2}. Hence, as
in~\eqref{eq:nice_1}, we obtain
\begin{align} \label{eq:nice_2}
  \| \nabla \partial^{\bsnu-\bsk} (v^G - v^G_{h})\|_{L^2}
  &\,\lesssim\, h \,
 \frac{(|\bsnu-\bsk|+2)!}{2}\,  \overline\bsb^{\bsnu-\bsk}\, \|G\|_{L^2}\,.
\end{align}
Substituting~\eqref{eq:nice_1} and~\eqref{eq:nice_2}
into~\eqref{eq:split}, and using $\|\psi_j\|_{L^\infty} = a_{\min} b_j \le
a_{\max} \overline{b}_j$, yields
\begin{align*}
  &|\partial^\bsnu G(u - u_{h})| \\
  &\,\lesssim\,
  a_{\max}\,\|f\|_{L^2}\,\|G\|_{L^2}\, h^2\, \bigg(
  \sum_{\bsk\le\bsnu} \sbinom{\bsnu}{\bsk}\, \frac{(|\bsk|+2)!}{2}\,\overline\bsb^\bsk\,
  \frac{(|\bsnu-\bsk|+2)!}{2}\,\overline\bsb^{\bsnu-\bsk} \nonumber\\
  &\qquad +
  \sum_{j\in\supp(\bsnu)} \nu_j\,\overline{b}_j\,
  \sum_{\bsk\le\bsnu-\bse_j} \sbinom{\bsnu-\bse_j}{\bsk}\,\frac{(|\bsk|+2)!}{2}\,\overline\bsb^\bsk\,
  \frac{(|\bsnu-\bse_j-\bsk|+2)!}{2}\,\overline\bsb^{\bsnu-\bse_j-\bsk} \bigg) \\
  &\,\lesssim\,
  a_{\max}\,\|f\|_{L^2}\,\|G\|_{L^2}\, h^2\,
  \bigg(\frac{(|\bsnu|+5)!}{120}
   +
  \sum_{j\in\supp(\bsnu)} \nu_j\,\frac{(|\bsnu-\bse_j|+5)!}{120} \bigg)\overline\bsb^\bsnu \\
  &\,\lesssim\,
  a_{\max}\,\|f\|_{L^2}\,\|G\|_{L^2}\, h^2\,
  \frac{(|\bsnu|+5)!}{120}\, \overline\bsb^\bsnu\,,
\end{align*}
where we used the identity~\eqref{eq:ident5}. This completes the proof.
\qed


\subsection*{\textbf{Proof of Lemma~\ref{lem:regu3}}}

This result is \cite[Theorem~14]{GKNSSS15}. We include the proof here
since, unlike in the uniform case where we do induction directly for the
quantity $\|\nabla(\partial^{\bsnu}u(\cdot,\bsy))\|_{L^2}$, here we need
to work with
$\|a^{1/2}(\cdot,\bsy)\,\nabla(\partial^{\bsnu}u(\cdot,\bsy))\|_{L^2}$,
and this technical step is needed for the subsequent proof.

Let $f\in V^*$ and $\bsy\in U_\bsb$. We first prove by induction on
$|\bsnu|$ that
\begin{align} \label{eq:induct3}
 \|a^{1/2}(\cdot,\bsy)\,\nabla(\partial^{\bsnu}u(\cdot,\bsy))\|_{L^2}
 \,\le\, \Lambda_{|\bsnu|}\,\bsbeta^\bsnu\, \frac{\|f\|_{V^*}}{\sqrt{a_{\min}(\bsy)}}\,,
\end{align}
where the sequence $(\Lambda_n)_{n\ge 0}$ is defined recursively by
\eqref{eq:lambda} and satisfies \eqref{eq:lambda-bound}.

We take $v = u(\cdot,\bsy)$ in the weak form \eqref{eq:weak} to obtain
\begin{align*}
  \int_D a \,|\nabla u |^2\,\rd\bsx
  &\;\le\; \|f\|_{V^*}\,\|u(\cdot,\bsy)\|_{V}
  \;\le\; \frac{\|f\|_{V^*}}{\sqrt{a_{\min}(\bsy)}}
  \left(\int_D a \, |\nabla u |^2 \,\rd\bsx\right)^{1/2},
\end{align*}
and then cancel the common factor from both sides to obtain
\eqref{eq:induct3} for the case $\bsnu = \bszero$. Given any multi-index
$\bsnu$ with $|\bsnu| = n\ge 1$, we apply $\partial^\bsnu$ to
\eqref{eq:weak} to obtain
\[
  \int_D \bigg( \sum_{\bsm\le\bsnu} \sbinom{\bsnu}{\bsm} (\partial^{\bsnu-\bsm} a)
  \,\nabla (\partial^{\bsm} u) \cdot\nabla z\bigg) \rd\bsx \,=\, 0
  \qquad\mbox{for all}\quad z\in V\;.
\]
Taking $z = \partial^\bsnu u(\cdot,\bsy)$, separating out the $\bsm =
\bsnu$ term, dividing and multiplying by $a$, and using the Cauchy-Schwarz
inequality, we obtain
\begin{align*}
  &\int_D a\, |\nabla (\partial^{\bsnu} u) |^2 \rd\bsx
  \;=\; - \sum_{\satop{\bsm\le\bsnu}{\bsm\ne\bsnu}} \sbinom{\bsnu}{\bsm}
  \int_D (\partial^{\bsnu-\bsm} a) \,\nabla (\partial^{\bsm} u)
  \cdot\nabla (\partial^{\bsnu} u) \, \rd\bsx \\
  &\;\le\; \sum_{\satop{\bsm\le\bsnu}{\bsm\ne\bsnu}} \sbinom{\bsnu}{\bsm}
  \bigg\|\frac{\partial^{\bsnu-\bsm} a (\cdot,\bsy)}{a(\cdot,\bsy)}\bigg\|_{L^\infty}
  \left(\int_D a |\nabla (\partial^{\bsm} u) |^2\,\rd\bsx\right)^{1/2}
  \left(\int_D a |\nabla (\partial^{\bsnu} u) |^2\,\rd\bsx\right)^{1/2}\,.
\end{align*}
We observe from \eqref{eq:axy-logn}
that
\begin{equation} \label{eq:diff-logn-1}
  \partial^{\bsnu-\bsm} a \,=\, a\prod_{j\ge 1}
  (\sqrt{\mu_j}\,\xi_j)^{\nu_j-m_j} \quad \text{for all} \quad
  \bsnu\ne\bsm\,,
\end{equation}
and therefore
\begin{align} \label{eq:diff-logn-2}
  \bigg\|\frac{\partial^{\bsnu-\bsm} a (\cdot,\bsy)}{a(\cdot,\bsy)}\bigg\|_{L^\infty}
  \,=\, \bigg\|\prod_{j\ge 1} (\sqrt{\mu_j}\,\xi_j)^{\nu_j-m_j}\bigg\|_{L^\infty}
  \,\le\, \bsbeta^{\bsnu-\bsm}\,.
\end{align}
Thus we arrive at
\begin{align*}
  \left(\int_D a |\nabla (\partial^{\bsnu} u) |^2\,\rd\bsx\right)^{1/2}
  &\,\le\, \sum_{\satop{\bsm\le\bsnu}{\bsm\ne\bsnu}} \sbinom{\bsnu}{\bsm}\,
  \bsbeta^{\bsnu-\bsm}
  \left(\int_D a |\nabla (\partial^{\bsm} u)|^2\,\rd\bsx\right)^{1/2}\,.
\end{align*}
We now use the inductive hypothesis that \eqref{eq:induct3} holds when
$\vert \bsnu\vert \leq n-1$ in each of the terms on the right-hand side to
obtain
\begin{align*}
  \left(\int_D a |\nabla (\partial^{\bsnu} u) |^2\,\rd\bsx\right)^{1/2}
  &\,\le\, \sum_{i=0}^{n-1} \sum_{\satop{\bsm\le\bsnu}{|\bsm|=i}} \sbinom{\bsnu}{\bsm}\,
  \bsbeta^{\bsnu-\bsm}\, \Lambda_i\, \bsbeta^\bsm\,
  \frac{\|f\|_{V^*}}{\sqrt{a_{\min}(\bsy)}} \\
  &\,=\, \sum_{i=0}^{n-1} \sbinom{n}{i} \Lambda_i\, \bsbeta^\bsnu
  \frac{\|f\|_{V^*}}{\sqrt{a_{\min}(\bsy)}}
  \,=\, \Lambda_n\, \bsbeta^\bsnu
  \frac{\|f\|_{V^*}}{\sqrt{a_{\min}(\bsy)}}\,,
\end{align*}
where we used the identity~\eqref{eq:ident0}. This completes the induction
proof of \eqref{eq:induct3}.

The desired bound in the lemma is obtained by applying
\eqref{eq:lambda-bound} with $\alpha = \ln 2$ on the right-hand side of
\eqref{eq:induct3}, and by noting that the left-hand side of
\eqref{eq:induct3} can be bounded from below by
$\sqrt{a_{\min}(\bsy)}\,\|\partial^{\bsnu} u(\cdot,\bsy) \|_V$. The case
$\bsnu=\bszero$ corresponds to \eqref{eq:apriori-logn}. This completes the
proof. \qed


\subsection*{\textbf{Proof of Lemma~\ref{lem:regu4a}}}

This result was proved in \cite{KSSSU} based on an argument similar to the
proof of Lemma~\ref{lem:regu2a} in the uniform case. The tricky point of
the proof is in recognizing that for the recursion to work in the
lognormal case we need to multiply the expression by
$a^{-1/2}(\cdot,\bsy)$, which is not intuitive.

Let $f\in L^2(D)$ and $\bsy\in U_{\overline\bsbeta}$. For any multi-index
$\bsnu\ne\bszero$, we apply $\partial^\bsnu$ to \eqref{eq:strong} to
obtain (formally, at this stage)
\begin{align*}
 \nabla\cdot \partial^\bsnu (a\nabla u) \,=\,
\nabla \cdot
  \bigg(\sum_{\bsm \le\bsnu} \sbinom{\bsnu}{\bsm}
  (\partial^{\bsnu-\bsm} a)\, \nabla (\partial^\bsm u) \bigg)
 \,=\, 0\,.
\end{align*}
Separating out the $\bsm = \bsnu$ term yields the following identity
\begin{align*}
  g_\bsnu \,:=\,
  \nabla\cdot(a\nabla(\partial^\bsnu u))
  &\,=\, - \nabla \cdot \bigg(\sum_{\satop{\bsm\le\bsnu}{\bsm\ne\bsnu}} \sbinom{\bsnu}{\bsm}
  (\partial^{\bsnu-\bsm} a)\, \nabla (\partial^\bsm u) \bigg) \\
  &\,=\, - \sum_{\satop{\bsm\le\bsnu}{\bsm\ne\bsnu}} \sbinom{\bsnu}{\bsm}
  \nabla\cdot \bigg(\frac{\partial^{\bsnu-\bsm} a}{a}\, (a \nabla (\partial^\bsm u)) \bigg) \\
  &\,=\, - \sum_{\satop{\bsm\le\bsnu}{\bsm\ne\bsnu}} \sbinom{\bsnu}{\bsm}
  \bigg(\frac{\partial^{\bsnu-\bsm} a}{a}\, g_\bsm \, + \, \nabla
    \bigg(\frac{\partial^{\bsnu-\bsm} a}{a}\bigg)\cdot (a \nabla
    (\partial^\bsm u)) \bigg)\,,
\end{align*}
where we used the identity \eqref{eq:nabla}. Due to Assumption~\ref{L2} we
may multiply $g_\bsnu$ by $a^{-1/2}$ and obtain, for any $|\bsnu| > 0$,
the recursive bound
\begin{align} \label{eq:step1}
  \|a^{-1/2}g_\bsnu\|_{L^2}
  &\,\le\, \sum_{\satop{\bsm\le\bsnu}{\bsm\ne\bsnu}} \sbinom{\bsnu}{\bsm}
  \bigg(\bigg\|\frac{\partial^{\bsnu-\bsm}
          a}{a}\bigg\|_{L^\infty}\, \|a^{-1/2} g_\bsm\|_{L^2}
  \nonumber \\
  &\qquad\qquad\qquad + \;  \bigg\|\nabla \bigg(\frac{\partial^{\bsnu-\bsm}
    a}{a}\bigg)\bigg\|_{L^\infty}\, \|a^{1/2} \nabla (\partial^\bsm u)\|_{L^2}
\bigg)
\;.
\end{align}
By assumption, $-g_\bszero =f\in L^2(D)$, so that we obtain (by induction
with respect to $|\bsnu|$) from \eqref{eq:step1} that
$a^{-1/2}(\cdot,\bsy)\,g_\bsnu(\cdot,\bsy)\in L^2(D)$, and hence from
Assumption~\ref{L2} that $g_\bsnu(\cdot,\bsy)\in L^2(D)$ for every
$\bsnu\in \indx$. The above formal identities therefore hold in $L^2(D)$.

To complete the proof, it remains to bound the above $L^2$ norm. Applying
the product rule to \eqref{eq:diff-logn-1} we obtain
\[
  \nabla\bigg(\frac{\partial^{\bsnu-\bsm} a}{a}\bigg)
  \,=\, \sum_{k\ge 1} (\nu_k-m_k)(\sqrt{\mu_k}\,\xi_k)^{\nu_k-m_k-1} (\sqrt{\mu_k}\,\nabla\xi_k)
        \prod_{\satop{j\ge 1}{j\ne k}} (\sqrt{\mu_j}\,\xi_j)^{\nu_j-m_j}\,.
\]
Due to the definition of $\overline{\beta}_j$ in \eqref{eq:barbj-logn},
this implies, in a similar manner to \eqref{eq:diff-logn-2}, that
\begin{equation} \label{eq:diff-logn-3}
\bigg\|\nabla\bigg(\frac{\partial^{\bsnu-\bsm} a}{a}\bigg)\bigg\|_{L^\infty}
 \,\le\, |\bsnu-\bsm| \,\overline\bsbeta^{\bsnu-\bsm}.
\end{equation}
Substituting \eqref{eq:diff-logn-2} and \eqref{eq:diff-logn-3} into
\eqref{eq:step1}, we conclude that
\begin{align*}
 \underbrace{\|a^{-1/2}g_\bsnu\|_{L^2}}_{\bbA_\bsnu}  \,\le\,
 \sum_{\satop{\bsm\le\bsnu}{\bsm\ne\bsnu}} \sbinom{\bsnu}{\bsm}\,\bsbeta^{\bsnu-\bsm}\,
 \underbrace{\|a^{-1/2}g_\bsm \|_{L^2}}_{\bbA_\bsm}
 \,+\, B_\bsnu\,,
\end{align*}
where
\begin{align*}
  B_\bsnu
  \,:=&\,
  \sum_{\satop{\bsm\le\bsnu}{\bsm\ne\bsnu}} \sbinom{\bsnu}{\bsm}\,
  |\bsnu-\bsm| \,\overline\bsbeta^{\bsnu-\bsm}\,
  \|a^{1/2} \nabla (\partial^\bsm u)\|_{L^2} \\
 \,\le&\,
  \sum_{\satop{\bsm\le\bsnu}{\bsm\ne\bsnu}} \sbinom{\bsnu}{\bsm}\,
  |\bsnu-\bsm| \,\overline\bsbeta^{\bsnu-\bsm}\,
  \Lambda_{|\bsm|}\,\bsbeta^\bsm\,
 \frac{\|f\|_{V^*}}{\sqrt{a_{\min}(\bsy)}}
 \,\le\,
  \overline\Lambda_{|\bsnu|}\,\overline\bsbeta^\bsnu
  \frac{\|f\|_{V^*}}{\sqrt{a_{\min}(\bsy)}}
  \,,
\end{align*}
where we used \eqref{eq:induct3} and again the identity \eqref{eq:ident0}
to write, with $n = |\bsnu|$,
\[
  \sum_{\satop{\bsm\le\bsnu}{\bsm\ne\bsnu}}
  \sbinom{\bsnu}{\bsm}\,
  |\bsnu-\bsm|\, \Lambda_{|\bsm|}
  \,=\, \sum_{i=0}^{n-1} \sum_{\satop{\bsm\le\bsnu}{|\bsm|=i}}
  \sbinom{\bsnu}{\bsm}\, (n-i)\, \Lambda_i
  \,=\, \sum_{i=0}^{n-1} \sbinom{n}{i}\, (n-i)\,
    \Lambda_i
  \,=:\, \overline\Lambda_n\,.
\]
Since $\bbA_\bszero = \|a^{-1/2} f\|_{L^2} \le
\|f\|_{L^2}/\sqrt{a_{\min}(\bsy)}$, we now define
\[
  \bbB_\bsnu \,:=\, C_{\rm emb}\, \overline\Lambda_{|\bsnu|}\,\overline\bsbeta^\bsnu
  \frac{\|f\|_{L^2}}{\sqrt{a_{\min}(\bsy)}}\,,
\]
so that $\bbA_\bszero\le \bbB_\bszero$ and $B_\bsnu\le \bbB_\bsnu$ for all
$\bsnu$. We may now apply Lemma~\ref{lem:recur2} to obtain
\begin{align}\label{eq:step2}
 \|a^{-1/2}g_\bsnu\|_{L^2}
 &\,\le\,
 \sum_{\bsk\le\bsnu} \sbinom{\bsnu}{\bsk}\Lambda_{|\bsk|}\,\bsbeta^\bsk\,
 C_{\rm emb}\, \overline\Lambda_{|\bsnu-\bsk|}\,\overline\bsbeta^{\bsnu-\bsk}
  \frac{\|f\|_{L^2}}{\sqrt{a_{\min}(\bsy)}}\,.
\end{align}
Note the extra factor $n-i$ in the definition of $\overline\Lambda_n$
compared to $\Lambda_n$ in \eqref{eq:lambda} so that
$\Lambda_n\le\overline\Lambda_n$. Using the bound in
\eqref{eq:lambda-bound} with $\alpha\le\ln 2$, we have
\begin{align*}
  \overline\Lambda_n
  \,\le\, \sum_{i=0}^{n-1} \sbinom{n}{i} (n-i)\,\frac{i!}{\alpha^i}
  \,=\, \frac{n!}{\alpha^n}\,\alpha\,\sum_{i=0}^{n-1} \frac{\alpha^{n-i-1}}{(n-i-1)!}
  \,=\, \frac{n!}{\alpha^n}\,\alpha\,\sum_{k=0}^{n-1} \frac{\alpha^k}{k!}
  \,\le\, \frac{n!}{\alpha^n}\,\alpha\, e^\alpha
  \,\le\, \frac{n!}{\alpha^n},
\end{align*}
where the final step is valid provided that $\alpha\,e^\alpha \le 1$. Thus
it suffices to choose $\alpha \le 0.567\cdots$. For convenience we take
$\alpha = 0.5$ to bound \eqref{eq:step2}. This together with the identity
\eqref{eq:ident1} gives
\begin{equation} \label{eq:step3}
 \|a^{-1/2}g_\bsnu\|_{L^2}
 \,\le\, C_{\rm emb}\,
 (|\bsnu|+1)!\,2^{|\bsnu|}\,\overline\bsbeta^\bsnu\,\frac{\|f\|_{L^2}}{\sqrt{a_{\min}(\bsy)}}\,.
\end{equation}

Since $a^{-1/2}g_\bsnu = a^{-1/2}\nabla\cdot(a\nabla(\partial^\bsnu u)) =
a^{1/2}\Delta(\partial^\bsnu u) + a^{-1/2}\,(\nabla a\cdot
\nabla(\partial^\bsnu u))$ by applying \eqref{eq:nabla}, we have
\begin{align*}
 \|a^{1/2}\Delta(\partial^\bsnu u)\|_{L^2}
 \,\le\, \|a^{-1/2}g_\bsnu\|_{L^2} \,+\, \|a^{-1/2}\,(\nabla a\cdot\nabla(\partial^\bsnu u))\|_{L^2}\,,
\end{align*}
which yields
\begin{align*}
 \sqrt{a_{\min}(\bsy)}\,\|\Delta(\partial^\bsnu u)\|_{L^2}
 \,\le\, \|a^{-1/2}g_\bsnu\|_{L^2}
 \,+\, \frac{\|\nabla a(\cdot,\bsy)\|_{L^\infty} }{a_{\min}(\bsy)} \|a^{1/2}\nabla(\partial^\bsnu u)\|_{L^2}\,,
\end{align*}
and in turn
\begin{align} \label{eq:step4}
 \|\Delta(\partial^\bsnu u)\|_{L^2}
 \,\le\, \frac{\|a^{-1/2}g_\bsnu\|_{L^2}}{\sqrt{a_{\min}(\bsy)}}
 \,+\, \frac{\|\nabla a(\cdot,\bsy)\|_{L^\infty} }{a_{\min}(\bsy)}\,
 \frac{\|a^{1/2}\nabla(\partial^\bsnu u)\|_{L^2}}{\sqrt{a_{\min}(\bsy)}}\,.
\end{align}
Substituting \eqref{eq:step3} and \eqref{eq:induct3} into
\eqref{eq:step4}, and using $\Lambda_{|\bsnu|}\le 2^{|\bsnu|}|\bsnu|!$ and
$\bsbeta^\bsnu\le \overline\bsbeta^\bsnu$, we conclude that
\begin{align*}
 \|\Delta(\partial^\bsnu u)\|_{L^2}
 &\,\le\, C_{\rm emb}\, \bigg( \frac{1}{a_{\min}(\bsy)} + \frac{\|\nabla a(\cdot,\bsy)\|_{L^\infty}}{a_{\min}^2(\bsy)}\bigg)
 (|\bsnu|+1)!\,2^{|\bsnu|}\,\overline\bsbeta^\bsnu\,\|f\|_{L^2}\,.
\end{align*}
This completes the proof. \qed


\subsection*{\textbf{Proof of Lemma~\ref{lem:regu4b}}}

This result was proved in \cite{KSSSU}. We include the proof here to
provide a complete unified view of the proof techniques discussed in this
survey.

Let $f\in L^2(D)$ and $\bsy\in U_{\overline{\bsbeta}}$. Following
\eqref{eq:Ph}--\eqref{eq:new_1} in the uniform case, we can write in the
lognormal case
\begin{align} \label{eq:new_1-logn}
  \| a^{1/2} \nabla \partial^\bsnu (u - u_{h})\|_{L^2}
  \,\le\, \| a^{1/2} \nabla \calP_h \partial^\bsnu (u - u_{h})\|_{L^2}
  \,+\, \|a^{1/2} \nabla (\calI - \calP_h) \partial^\bsnu u \|_{L^2}\,.
\end{align}
Now, applying $\partial^\bsnu$ to \eqref{eq:GalerkinOrth} and separating
out the $\bsm = \bsnu$ term, we get for all $z_h\in V_h$ in the lognormal
case that
\begin{align} \label{eq:new_2-logn}
 \int_D a\, \nabla \partial^{\bsnu} (u-u_h)\cdot\nabla z_h \,\rd \bsx
 \,=\, - \sum_{\satop{\bsm\le\bsnu}{\bsm\ne\bsnu}} \sbinom{\bsnu}{\bsm}
 \int_D (\partial^{\bsnu-\bsm} a)
  \nabla\partial^{\bsm}(u-u_h)\cdot\nabla z_h \,\rd \bsx\,.
\end{align}
Choosing $z_h = \calP_h \partial^\bsnu (u-u_h)$ and using the definition
\eqref{eq:Ph} of $\calP_h$, the left-hand side of \eqref{eq:new_2-logn} is
equal to $\int_D a\, |\nabla \calP_h
\partial^{\bsnu} (u-u_h)|^2 \,\rd \bsx$. Dividing and multiplying the
right-hand side of \eqref{eq:new_2-logn} by $a$, and using the
Cauchy-Schwarz inequality, we then obtain
\begin{align*}
 &\int_D a\, |\nabla \calP_h \partial^{\bsnu} (u-u_h)|^2 \,\rd \bsx \\
 &\ \ \le\, \sum_{\satop{\bsm\le\bsnu}{\bsm\ne\bsnu}} \sbinom{\bsnu}{\bsm}
 \left\|\frac{\partial^{\bsnu-\bsm} a}{a}\right\|_{L^\infty}
 \left(\int_D a\, |\nabla\partial^{\bsm}(u-u_h)|^2\,\rd\bsx \right)^{\frac{1}{2}}
 \left(\int_D a\, |\nabla \calP_h \partial^{\bsnu} (u-u_h)|^2\,\rd\bsx \right)^{\frac{1}{2}}.
\end{align*}
Canceling one common factor from both sides and using
\eqref{eq:diff-logn-2}, we arrive at
\begin{align} \label{eq:new_3-logn}
 \|a^{1/2} \nabla \calP_h \partial^{\bsnu} (u-u_h)\|_{L^2}
 &\,\le\, \sum_{\satop{\bsm\le\bsnu}{\bsm\ne\bsnu}} \sbinom{\bsnu}{\bsm}\bsbeta^{\bsnu-\bsm}\,
 \|a^{1/2} \nabla\partial^{\bsm}(u-u_h)\|_{L^2}.
\end{align}
Substituting \eqref{eq:new_3-logn} into \eqref{eq:new_1-logn}, we then
obtain
\begin{align*}
  \underbrace{\| a^{1/2} \nabla \partial^\bsnu (u - u_{h})\|_{L^2}}_{\bbA_\bsnu} 
  &\,\le\, \sum_{\satop{\bsm\le\bsnu}{\bsm\ne\bsnu}} \sbinom{\bsnu}{\bsm}\bsbeta^{\bsnu-\bsm}
  \underbrace{\|a^{1/2} \nabla\partial^{\bsm}(u-u_h)\|_{L^2}}_{\bbA_\bsm}
  \,+\, \underbrace{\|a^{1/2} \nabla (\calI - \calP_h) \partial^\bsnu u
    \|_{L^2}}_{\bbB_\bsnu}\,.
\end{align*}
Note that we have $\bbA_\bszero = \bbB_\bszero$. Now applying
Lemma~\ref{lem:recur2} with $\alpha=0.5$, together with
\eqref{eq:FEerror-logn}, Lemma~\ref{lem:regu4a} and \eqref{eq:ident2}, we
conclude that
\begin{align*}
  &\| a^{1/2} \nabla \partial^\bsnu (u - u_{h})\|_{L^2}
  \,\le\, \sum_{\bsm\le\bsnu} \sbinom{\bsnu}{\bsm} \Lambda_{|\bsm|}\,\bsbeta^\bsm\,
  \|a^{1/2} \nabla (\calI - \calP_h) \partial^{\bsnu-\bsm} u
  \|_{L^2} \\
  &\,\lesssim\, h\, a_{\max}^{1/2}(\bsy)\, \sum_{\bsm\le\bsnu} \sbinom{\bsnu}{\bsm} \Lambda_{|\bsm|}\, \bsbeta^\bsm\,
  \|\Delta (\partial^{\bsnu-\bsm} u)\|_{L^2} \\
  &\,\lesssim\, h\, T(\bsy)\,a_{\max}^{1/2}(\bsy)\,
  \sum_{\bsm\le\bsnu} \sbinom{\bsnu}{\bsm} \, |\bsm|!\,2^{|\bsm|}\, \bsbeta^\bsm\,
  (|\bsnu-\bsm|+1)!\,2^{|\bsnu-\bsm|}\,\overline\bsbeta^{\bsnu-\bsm}\, \|f\|_{L^2} \\
  &\,\lesssim\, h T(\bsy)\,\, a_{\max}^{1/2}(\bsy)\,
  \frac{(|\bsnu|+2)!}{2}\,2^{|\bsnu|}\, \overline\bsbeta^{\bsnu}\, \|f\|_{L^2} \,,
\end{align*}
where $T(\bsy)$ is defined in \eqref{eq:defT}. This completes the proof.
\qed


\subsection*{\textbf{Proof of Lemma~\ref{lem:regu4c}}}

This result was proved in \cite{KSSSU}. Again we include the proof here to
provide a complete unified view of the proof techniques discussed in this
survey.

 Let $f,G\in L^2(D)$ and $\bsy\in U_{\overline{\bsbeta}}$.
Following \eqref{eq:nitsche_1}--\eqref{eq:duality} in the uniform case,
and using the Leibniz product rule~\eqref{eq:leibniz}, we have for the
lognormal case that
\begin{align*}
  &\partial^\bsnu G(u - u_{h})
  \,=\,
  \int_D \partial^\bsnu \left(a \,\nabla (u
  - u_{h}) \cdot \nabla (v^G - v^G_h) \right)\, \rd\bsx \\
  &\,=\,
  \int_D \sum_{\bsm\le\bsnu} \sbinom{\bsnu}{\bsm} (\partial^{\bsnu-\bsm} a)\,
  \partial^\bsm \left( \nabla (u - u_{h}) \cdot \nabla (v^G - v^G_h) \right)\, \rd\bsx \\
  &\,=\,
  \int_D \sum_{\bsm\le\bsnu} \sbinom{\bsnu}{\bsm} (\partial^{\bsnu-\bsm} a)\,
  \sum_{\bsk\le\bsm} \sbinom{\bsm}{\bsk} \nabla \partial^\bsk (u- u_{h}) \cdot
  \nabla \partial^{\bsm-\bsk} (v^G - v^G_h)\, \rd\bsx \\
  &\,=\,
  \int_D \sum_{\bsm\le\bsnu} \sbinom{\bsnu}{\bsm} \frac{\partial^{\bsnu-\bsm} a}{a}
  \sum_{\bsk\le\bsm} \sbinom{\bsm}{\bsk}
  \left(a^{1/2}\nabla \partial^\bsk (u - u_{h})\right) \cdot
  \left(a^{1/2}\nabla \partial^{\bsm-\bsk} (v^G - v^G_h)\right)\, \rd\bsx.
\end{align*}
Using the Cauchy-Schwarz inequality and \eqref{eq:diff-logn-2}, we obtain
\begin{align} \label{eq:split-logn}
 &|\partial^\bsnu G(u - u_{h}) | \nonumber\\
 &\,\le\,
 \sum_{\bsm\le\bsnu} \sbinom{\bsnu}{\bsm}
 \bsbeta^{\bsnu-\bsm}
 \sum_{\bsk\le\bsm} \sbinom{\bsm}{\bsk}
 \| a^{1/2} \nabla \partial^\bsk (u - u_{h})\|_{L^2}\,
 \| a^{1/2} \nabla \partial^{\bsm-\bsk} (v^G - v^G_h)\|_{L^2}\,.
\end{align}

We have from Lemma~\ref{lem:regu4b} that
\begin{align} \label{eq:nice_1-logn}
  \| a^{1/2} \nabla \partial^\bsk (u - u_{h})\|_{L^2}
  \,\lesssim\, h \, T(\bsy)\,a_{\max}^{1/2}(\bsy)\,
   \frac{(|\bsk|+2)!}{2}\, 2^{|\bsk|}\, \overline\bsbeta^{\bsk}\,\|f\|_{L^2}\,.
\end{align}
Since the bilinear form $\scA(\bsy;\cdot,\cdot)$ is symmetric and since
the representer $g$ for the linear functional $G$ is in $L^2$, all the
results hold verbatim also for the adjoint problem \eqref{eq:nitsche_1}
and for its FE discretisation \eqref{eq:nitsche_2}. Hence, as in
\eqref{eq:nice_1-logn}, we obtain
\begin{align} \label{eq:nice_2-logn}
  \| a^{1/2} \nabla \partial^{\bsm-\bsk} (v^G - v^G_{h})\|_{L^2}
  &\lesssim h T(\bsy)\,\,a_{\max}^{1/2}(\bsy)\,
  \frac{(|\bsm-\bsk|+2)!}{2}\,2^{|\bsm-\bsk|}\,\overline\bsbeta^{\bsm-\bsk}\,\|G\|_{L^2}\,.
\end{align}
Substituting \eqref{eq:nice_1-logn} and \eqref{eq:nice_2-logn} into
\eqref{eq:split-logn}, and using the identity \eqref{eq:ident5}, we obtain
\begin{align*}
 |\partial^\bsnu G(u - u_{h}) |
 &\,\lesssim\,
 h^2\, T^2(\bsy)\,a_{\max}(\bsy)\,
 \sum_{\bsm\le\bsnu} \sbinom{\bsnu}{\bsm}
 \frac{(|\bsm|+5)!}{120}\,2^{|\bsm|}\,\overline\bsbeta^{\bsnu}\,
 \|f\|_{L^2} \,\|G\|_{L^2}\,.
\end{align*}
Using again \eqref{eq:ident0}, with $n=|\bsnu|$ we have
\begin{align*}
 &\sum_{\bsm\le\bsnu} \sbinom{\bsnu}{\bsm} 2^{|\bsm|}\,\frac{(|\bsm|+5)!}{120} \\
 &\,=\, \sum_{i=0}^{n} \sbinom{n}{i} 2^{i}\,\frac{(i+5)!}{120}
\,=\, n!\,\sum_{i=0}^{n} \frac{(i+1)(i+2)(i+3)(i+4)(i+5) 2^{i}}{120(n-i)!}
 \,\le\, \frac{(n+5)!}{120} 2^n e\,.
\end{align*}
This yields the required bound in the lemma. This completes the proof.
\qed

\begin{acknowledgements}
We graciously acknowledge many insightful discussions and valuable
comments from our collaborators Josef Dick, Mahadevan Ganesh, Thong Le
Gia, Alexander Gilbert, Ivan Graham, Yoshihito Kazashi, James Nichols,
Pieterjan Robbe, Robert Scheichl, Christoph Schwab, and Ian Sloan. We
especially thank Mahadevan Ganesh for suggesting an alternative proof
strategy to improve some existing estimates. We are also grateful for the
financial supports from the Australian Research Council (FT130100655 and
DP150101770) and the KU Leuven research fund (OT:3E130287 and
C3:3E150478).
\end{acknowledgements}


\end{document}